\documentclass[12pt]{amsart}
\usepackage{amsmath}
\usepackage{amscd}
\usepackage{amsfonts}
\usepackage{amssymb}
\usepackage{amsthm}

  =4.mm
\newcommand{\rar}{\rightarrow}
\newcommand{\calg}{\mathcal}
\newcommand{\hb}{\widehat{B}}
\newcommand{\hc}{\widehat{C}}
\newcommand{\ilim}{\text{lim } }
\newcommand{\rhh}{\text{RHom}}
\newcommand{\rhhh}{\text{R}\mathcal{H}om}
\newcommand{\hhh}{\mathcal{H}om}
\newcommand{\chh}{{\mathcal{H}om}^{\text{cont}}}
\newcommand{\dd}{\text{D}_{poly}}
\newcommand{\strc}{\mathcal{O}_{X}}
\newcommand{\strcc}[1]{\mathcal{O}_{#1}}

\newcommand{\tenp}{\otimes_{\mathcal{O}_{X}}}

\newcommand{\dcat}{\text{D}^{b}(X)}
\newcommand{\edcat}{\widetilde{\text{D}^{b}(X)}}

\newcommand{\dcatb}{\text{D}^{+}(X)}
\newcommand{\bbcc}{\text{Ch}^{+}({{\mathcal O}_X}-{\text{mod}})}
\newcommand{\bcc}{\text{Ch}^{b}({{\mathcal O}_X}-{\text{mod}})}
\newcommand{\bacc}{\text{Ch}^{-}({{\mathcal O}_X}-{\text{mod}})}

\newcommand{\exal}[2]{ \wedge^{#2} \Omega_{#1} [#2]}
\newcommand{\texal}[2]{ \wedge^{#2} T_{#1} [-#2]}
\newcommand{\hkr}{\text{I}_{HKR}}
\newcommand{\thkr}{\widehat{\text{I}_{HKR}}}
\newcommand{\pbw}{\text{I}_{PBW}}
\newcommand{\bhkr}{\text{J}}

\newcommand{\sss}{{\textbf{S}}^{\bullet}}
\newcommand{\ttt}{\text{t}}
\newcommand{\enn}{\text{End}}
\newcommand{\ennnn}{{\mathcal E}\text{nd}}
\newcommand{\tot}{\text{Tot}}
\newcommand{\tdd}{\triangledown}
\newcommand{\chn}{\text{ch}}
\newcommand{\caph}{\text{H}}
\newcommand{\exx}{\text{exp}}
\newcommand{\adj}{\text{ad}}
\newcommand{\id}{\text{id}}
\newcommand{\pll}{\Pi \text{LIE}}
\newcommand{\stdp}{\text{END}_{T[-1]}}

\newcommand{\frel}[1]{{\mathcal L}(#1)}
\newcommand{\ufrel}[1]{\mathcal{T}(#1)}
\newcommand{\smap}[1]{\text{J}(#1)}
\newcommand{\bsmap}[1]{\text{B}(#1)}
\newcommand{\shuf}[2]{\text{Sh}_{#1,#2}}
\newcommand{\fld}{\mathbb K}

\newtheorem{thm}{Theorem}
\newtheorem*{thm1}{Theorem 1}
\newtheorem*{thm2p}{Theorem 2'}
\newtheorem*{thm3}{Theorem 3}
\newtheorem*{thm4}{Theorem 4}
\newtheorem*{thm5}{Theorem 5}
\newtheorem*{thmrr}{the relative Riemann-Roch theorem}
\newtheorem{prop}{Proposition}
\newtheorem{cor}{Corollary}
\newtheorem{lem}{Lemma}

\title{ The relative Riemann-Roch theorem from Hochschild homology}

\author{Ajay C. Ramadoss}
\address{Department of Mathematics, University of Oklahoma \newline
\text{        }       621 Elm Avenue, Norman,  OK-73019}
\email{aramadoss@math.ou.edu}

\begin{document}

\maketitle

\begin{abstract}

This paper attempts to clarify a preprint of Markarian [2]. The
preprint by Markarian [2] proves the relative Riemann-Roch theorem
using a result describing how the HKR map fails to respect
comultiplication. This paper elaborates on the core computations in
[2]. These computations show that the HKR map twisted by the square
root of the Todd genus "almost preserves" the Mukai pairing. This
settles a part of a conjecture of Caldararu [3]. The relative
Riemann-Roch theorem follows from this and a result of Caldararu
[4].

\end{abstract}

\section*{Introduction}

The purpose of this paper is to explain in detail an alternative
approach to the relative Riemann-Roch theorem which first appeared
in a very interesting but cryptic preprint [2] of Markarian. This
approach leads to a proof of the relative Riemann-Roch theorem by a
direct computation of the pairing on Hodge cohomology to which the
Mukai pairing on Hochschild homology defined by Caldararu [4]
descends via the Hochschild-Kostant-Rosenberg map multiplied by the
square root of the Todd genus. In the framework of this approach,
the fundamental reason for the appearance of the Todd genus in the
Riemann-Roch theorem is the failure of the
Hochschild-Kostant-Rosenberg map $\hkr$ introduced in Section 2 to
respect co-multiplication. One of the main theorems in [2] (Theorem
1 of [2]) describes the Duflo like error term that measures by how
much $\hkr$ fails to respect co-multiplication. The proof supplied
in [2] however, has a nontrivial error. A "dual" version of this
result equivalent to the original theorem has since been correctly
proven by the author in [1]. A correct proof of another version of
this result has been outlined by Markarian himself in a revised
version [6] of [2]. Theorem 1 of [2] appears as Theorem 2' in this
paper. \\

Let $X$ be a smooth proper scheme over a field $\fld$ of
characteristic $0$. We use Theorem 2' of this paper to prove the
main result (Theorem 1) of the current paper. This explicitly
interrelates the Hochschild-Kostant-Rosenberg map, the twisted HKR
map introduced in Section 1 and a map which we call the duality map
between $\text{RHom}_X(\strc, \Delta^* \strcc{\Delta})$ and
$\text{RHom}_{X}(\Delta^* \strcc{\Delta},S_X)$. Here, $S_X$ is the
shifted line bundle on $X$ tensoring with which is the Serre duality
functor on $\dcat$. Theorem 1 of this paper is equivalent to a
corrected version of an erroneous result (Theorem 8 of [2]) that appears in [2]. \\

Theorem 1 enables us to compute the pairing on Hodge cohomology to
which the Mukai pairing on Hochschild homology defined by Caldararu
[4] descends via the HKR-map twisted by the square root of the Todd
genus . Given Theorem 1, doing this computation is fairly easy. This
pairing on Hodge cohomology (see Proposition 5 of this paper) is
very similar to the generalized Mukai pairing Caldararu defined in
[3]. In particular, it satisfies the adjointness property one
expects from the Mukai pairing. Moreover, it coincides with the
Mukai pairing defined by Caldararu [3] on Mukai vectors (see [3]) of
elements of $\dcat$. However, the pairing obtained in Proposition 5
is not exactly the same as Caldararu's Mukai pairing on Hodge
co-homology. This settles a part of Caldararu's conjecture in [3]
regarding the equivalence between the Hochschild and Hodge
structures of a smooth proper complex variety - to be able to say
that the HKR map twisted by the square root of the Todd genus
preserves the Mukai pairing, one has to replace the Mukai pairing
defined by Caldararu [3] with the pairing that shows up in
Proposition 5. Since pairing in Proposition 5 does not in general
coincide with the Mukai pairing on Hodge cohomology for K3-surfaces,
we stay
contented by saying that the HKR map twisted by the square root of the Todd genus "almost preserves" the Mukai pairing. \\

The relative Riemann-Roch theorem follows from Proposition 5 and the adjointness property of the Mukai pairing on Hochschild homology. The
adjointness property of the Mukai pairing on Hochschild homology was proven in a paper [4] of Caldararu. \\

In order to prove Theorem 1 , we elaborate upon the core computations in [2]. These computations appeared in [2] in a very cryptic way. Some of
these computations do not appear in [6], whose approach differs in some details from [2]. In particular, unlike [2],[6] does not contain a
result equivalent to Theorem 1 and does not compute the Mukai pairing on Hochschild homology at the level of Hodge co-homology.  \\

The key steps in this computational approach are covered by Theorem
2' and Lemmas 2,3 and 4 of this paper. The most crucial
computations, Theorem 2' and Lemma 4 are related to very familiar
computations in elementary Lie theory. Theorem 2' is related to
computing the pull-back of a left invariant $1$-form on a Lie group
$G$ via the exponential map. Similarly, Lemma 4 is related computing
the pull-back of a left invariant volume form on $G$ via the the map
$\bar{\text{exp}}$ where $\bar{\text{exp}}(Z) = \text{exp}(-Z)$ for
any element $Z$ of the Lie algebra $\mathfrak{g}$ of $G$. We aim to
make these relations transparent by developing a "dictionary" in
this paper in three separate subsections containing remarks meant
for this purpose only. \\

\subsection*{The layout of this paper}

Section 1 begins by introducing the notations and conventions that shall be used for the rest of this paper. It then goes on to state Theorem 1
after defining the various maps involved in Theorem 1. The pairing on Hodge cohomology to which the Mukai pairing on Hochschild homology
descends ( via the HKR map twisted by the square root of the Todd genus ) is then computed. Finally, Section 1 uses this computation to prove
the relative Riemann-Roch theorem. The remaining sections of this paper are
devoted to proving Theorem 1.\\

Section 2 introduces two "connections" on the complex of completed
Hochschild chains of a smooth scheme $X$. Their properties are
proven in various propositions in this section. This section also
proves Theorem 2'. Sections 2.3 and 2.5 develop the "dictionary"
making the analogy between Theorem 2' and its
counterpart in elementary Lie theory more transparent. \\

Section 3 consists of a number of definitions, technical
propositions and two lemmas (Lemma 1 and Lemma 2) pertaining mainly
to linear algebra. These are used in later sections at various
points. The definitions of this section are important to understand
later sections. Proofs of propositions in later sections time and
again
refer to propositions in this section .\\

The key result of Section 4 is Lemma 3. This in turn follows from Lemma 4 and Lemma 2. Besides proving Lemma 3 and Lemma 4, Section 4 has a
subsection (Section 4.3) which explains the analogy between Lemma 4 and its counterpart in basic Lie theory. Section 4.3 is the last of the
three sections (2.3,2.5 and 4.3) developing the "dictionary" in this paper. \\

Section 5 undertakes the final computations leading to the proof of Theorem 1. \\

\subsection*{Acknowledgements} I thank Prof. Victor Ginzburg for
introducing me to the works of Markarian and Caldararu and for
useful discussions. I also thank Prof. Madhav Nori for useful
discussions and Prof. Andy Magid and Dr. Victor Protsak for useful
suggestions regarding the presentation of this paper.

\section{The main theorem , the Mukai pairing and the relative Riemann-Roch theorem}

We begin this section by clarifying some notation and conventions that shall be followed throughout this paper. Immediately after that, in
Section 1.1, we state the main theorem (Theorem 1) of this paper after describing the maps involved. This section then goes on to explain in
detail why
Theorem 1 implies the relative Riemann-Roch theorem. This is done in Sections 1.2 and 1.3. \\

\subsection*{Notation and conventions}
Let $X$ be a smooth proper scheme over a field $\fld$ of
characteristic $0$. All schemes and complex varieties that we
encounter in this paper are assumed to be proper. Let $\Delta: X
\rar X \times X $ denote the diagonal embedding. Let $p_1$ and $p_2$
denote the projections from $X \times X$ onto the first and second
factors respectively. As usual , $\strc$ denotes
the structure sheaf of $X$. \\

 $\bcc$ (resp. $\bacc$ and $\bbcc$) denotes the category of
  bounded (resp. bounded above and bounded
below) chain complexes of $\strc$-modules with coherent co-homology.
$\dcat$ denotes the bounded derived category of complexes of
$\strc$-modules with coherent co-homology. Similarly,
$\text{D}^{b}(X \times X)$ denotes the bounded derived category of
complexes
of $\strcc{X \times X}$-modules with coherent co-homology. \\

Whenever $f:X \rar Y$ is a morphism of schemes, $f^*: \text{D}^{b}(Y) \rar \dcat$ denotes the left derived functor of the pullback via $f$.
Similarly, $f_*:\dcat \rar \text{D}^{b}(Y)$ denotes the  right derived functor of the push-forward via $f$. $\strcc{\Delta}$ will denote
$\Delta_* \strc$. Similarly, when we refer to a tensor product, we will mean the corresponding left-derived functor unless stated otherwise
explicitly. Also, if $\calg E$ is
an object of $\dcat$ and $\varphi$ is a morphism in $\dcat$ ,\\
$\calg E \otimes \varphi$ shall denote the morphism $\textbf{1}_{\calg E} \otimes \varphi$. At times, $\calg E$ shall be used to denote the
morphism $\textbf{1}_{\calg E}$.
\\

If $\calg E $ and $\calg F$ are objects of $\dcat$,
$\rhh_{\dcat}(\calg E,\calg F)$ shall be denoted by $\rhh_X(\calg
E,\calg F)$. Similarly, $\rhh_{\text{D}^{b}(X \times X)}(-,-)$ shall
be denoted by $\rhh_{X \times X}(-,-)$. \\

$\Omega_X$ denotes the cotangent bundle of $X$. $\sss(\Omega_X[1])$ will denote the symmetric algebra generated over $\strc$ by $\Omega_X$
concentrated in degree $-1$. Note that $\sss(\Omega_X[1]) = \oplus_i \exal{X}{i}$. We shall often denote $\wedge^i \Omega_X$ by $\Omega_X^i$.
From Section 2 onwards, $\Omega_X$ and $\sss(\Omega_X[1])$ shall be denoted by $\Omega$ and $\sss(\Omega[1])$ respectively. The tangent bundle
of $X$ shall be denoted by
$T_X$ in this section and $T$ from Section 2 onwards. \\

Where convenient,we shall denote the Hodge co-homology $\oplus_{p,q}
\text{H}^q(X,\Omega_X^p)$ by $\text{H}^*(X)$ and $
\text{H}^q(X,\Omega_X^p)$ by $\text{H}^{p,q}(X)$. Note that
$\text{H}^*(X) = \rhh_X(\strc,\sss(\Omega_X[1]))$. The product on
$\text{H}^*(X)$ induced by the wedge product on $\sss(\Omega_X[1])$
will be denoted by $\wedge$. However, we shall often suppress
the $\wedge$ : If $a,b \in \text{H}^*(X)$, $ab$ should be understood to mean $a \wedge b \in \text{H}^*(X)$.    \\

Despite our attempts to minimize abuse of notation, it does happen
at times. There are many situations in this paper where we encounter
maps between tensor products of various objects in $\dcat$ that
rearrange factors. Very often, such maps are denoted by the symbol
$\tau$. In each
such situation, we specify what $\tau$ means unless we feel it is obvious to the reader. \\

\subsection{The crux of this paper}

Recall from Yekutieli [5] that the completed complex of Hochschild chains $\hc^{\bullet}(X) \in \bacc$ is a complex of flat $\strc$-modules that
represents $\Delta^* \strcc{\Delta}$ in $\dcat$. Upper indexing is used here to convert what was originally
a chain complex into a co-chain complex. \\

The Hochschild-Kostant-Rosenberg (HKR) map $\hkr$ from
$\hc^{\bullet}(X)$ to $\sss(\Omega_X[1])$ is a map of complexes of
$\strc$ modules. We describe this in greater detail in Section 2. We
identify $\Delta^* \strcc{\Delta}$ with $\hc^{\bullet}(X)$. Thus,
the HKR map $\hkr$ can be thought of as a map in $\dcat$ from
$\Delta^*
\strcc{\Delta}$ to $\sss(\Omega_X[1])$. \\

Let $S_X$ denote the object $\Omega^n[n]$ in $\dcat$. Let $\pi_n$
denote the projection from $\sss(\Omega_X[1])$ to the direct
summmand $\Omega^n[n]$. Consider the pairing
$$\langle -,- \rangle: \sss(\Omega_X[1]) \otimes \sss(\Omega_X[1]) \rar
S_X$$ given by the composite $$
\begin{CD} \sss(\Omega_X[1]) \otimes \sss(\Omega_X[1]) @>(- \wedge -)>> \sss(\Omega_X[1]) @>\pi_n>> S_X \end{CD} $$ of morphisms in $\dcat$.
 One also has a
$\textbf{ twisted HKR map}$ from $\sss(\Omega_X[1])$ to
\\ $\rhhh_X(\Delta^* \strcc{\Delta}, S_X)$. This arises out of the
composite

\begin{equation} \label{1}  \begin{CD} \Delta^*{\strcc{\Delta}}
\otimes \sss(\Omega[1]) @> \hkr \otimes
{\textbf{1}_{\sss(\Omega_X[1])}}
>> \sss(\Omega_X[1]) \otimes \sss(\Omega_X[1]) @>{\langle,\rangle}>>
S_X \\ \end{CD} \end{equation}

of morphisms in $\dcat$. We denote the twisted HKR map by $\thkr$.
\\

$\textbf{The duality map:}$ The material in this paragraph is
recalled from Caldararu [4]. Recall that $\Delta^*:\text{D}^{b}(X
\times X) \rar \dcat$ is the left adjoint of $\Delta_*:\dcat \rar
\text{D}^{b} (X \times X)$. Also recall that the functor of
tensoring by the shifted line bundle $S_X$ is the Serre duality
functor on $\dcat$. Similarly, tensoring by the shifted line bundle
$S_{X \times X}$ is the Serre duality functor on $\text{D}^{b}(X
\times X)$. We denote the functor given by tensoring by a shifted
line bundle $\calg L$ by $\calg L$ itself. The left adjoint
$\Delta_!:\dcat \rar
\text{D}^{b}(X \times X)$ of $\Delta^*$ is given by ${S_{X \times X}}^{-1}\Delta_* S_X$.   \\

Since $\Delta_!$ is the left adjoint of $\Delta^*$ we have an
isomorphism \begin{equation} \label{2} \calg I:
\rhh_X(\strc,\Delta^* \strcc{\Delta}) \simeq \rhh_{X \times
X}(\Delta_! \strc, \strcc{\Delta}) \text{ . } \end{equation}  Now,
$\Delta_! \strc = {S_{X \times X}}^{-1}\Delta_* S_X \simeq \Delta_*
S_X^{-1}$. We also have an isomorphism

\begin{equation} \label{3} {\calg T}: \rhh_{X \times X}(\Delta_!
\strc, \strcc{\Delta}) \simeq \rhh_{X \times X}(\strcc{\Delta},
\Delta_* S_X) \end{equation}

 given by tensoring an element of
$\rhh_{X \times X}(\Delta_! \strc, \strcc{\Delta})$ on the right by the shifted line bundle $p_2^* S_X$ and making the obvious identifications.
Now,since $\Delta^*$ is the left adjoint of $\Delta_*$ we have an isomorphism
 \begin{equation} \label{4} \calg
J: \rhh_{X \times X}(\strcc{\Delta}, \Delta_* S_X) \simeq
\rhh_X(\Delta^* \strcc{\Delta}, S_X) \text{ . }
\end{equation}

Let $D_{\Delta}:\rhh_X(\strc,\Delta^* \strcc{\Delta}) \rar
\rhh_X(\Delta^* \strcc{\Delta},S_X)$ denote the composite $\calg J
\circ {\calg T} \circ \calg I$. We refer to $D_{\Delta}$ as the
$\textbf{duality map}$.\\

\subsection*{The main theorem} The main theorem of this paper relates
the HKR, twisted HKR and duality maps. This is a corrected version
of Theorem 8 of Markarian's preprint [2]. \\

Note that $\hkr$ induces a map $$\hkr: \rhh_X(\strc, \Delta^*
\strcc{\Delta}) \rar \rhh_X(\strc, \sss(\Omega_X[1])) =
  \text{H}^*(X) \text{ . }$$ Similarly, $\thkr$ induces a map
$$\thkr:\text{H}^*(X) = \rhh_X(\strc,
 \sss(\Omega_X[1])) \rar \rhh_X(\Delta^* \strcc{\Delta}, S_X)
\text{ . }$$ Let $J$ denote the endomorphism on $\sss(\Omega_X[1])$
that multiplies $\exal{X}{i}$ by ${(-1)}^i$. $J$ induces an
endomorphism on $\text{H}^q(X,\Omega_X^p)$. Let $\text{td}(T_X) \in
\oplus_i \text{H}^i(X, \Omega_X^i)$ denote the Todd genus of the
tangent bundle of $X$. Recall that the wedge product $(\text{ - }
\wedge \text{ - }): \sss(\Omega_X[1]) ^{\otimes 2} \rar
\sss(\Omega_X[1])$ induces a product on $\text{H}^*(X)$. We are now
in a position to state the main theorem.

\begin{thm}

The following diagram commutes: \\

$$\begin{CD}
\rhh_X(\strc, \Delta^* \strcc{\Delta}) @>D_{\Delta}>>
\rhh_X(\Delta^* \strcc{\Delta}, S_X) \\
@VV{\hkr}V   @A{\thkr}AA \\
\text{H}^*(X) @>{ (- \wedge \text{td}(T_X)) \circ
J}>> \text{H}^*(X) \\
\end{CD} $$

The map in the bottom row of the above diagram takes an element
$\alpha \in \text{H}^*(X)$ to \\ $J(\alpha) \wedge \text{td}(T_X) $.

\end{thm}

Theorem 1 can be thought of as an explicit computation of the
duality map.

\subsection{The Mukai pairing}
We now try to understand how Theorem 1 leads to the relative
Riemann-Roch theorem. It is in this attempt that we see how Theorem
1 helps us calculate what the Mukai pairing on Hochschild homology
[4] descends to in Hodge co-homology. This settles a part of a
conjecture by Caldararu in [3]. \\

Let $\text{HH}_{i}(X)$ denote $\text{Hom}_{\dcat}(\strc, \Delta^* \strcc{\Delta}[i])$. $\text{HH}_i(X)$ is called the $i$-th {\it Hochschild
homology} of $X$. Let $\calg I$, ${\calg T}$ and $\calg J$ be as in \eqref{2}, \eqref{3}, and \eqref{4} respectively. Let $\text{tr}_X$ and
$\text{tr}_{X \times X}$ denote the canonical identifications of $\text{Hom}_{\dcat}(\strc,S_X)$ and $\text{Hom}_{\text{D}^b(X \times
X)}(\Delta_*S_X^{-1}, \Delta_* S_X)$ with $\fld$ respectively. We recall that Caldararu [4] defined a Mukai pairing on Hochschild homology. This
was a pairing

\begin{equation} \label{5} \text{HH}_{i}(X) \times \text{HH}_{-i}(X)
\rar \fld \end{equation}

\begin{equation*} (v,w) \leadsto \text{tr}_{X \times X}({\calg T}(\calg
I(v)) \circ \calg I(w)) \text{ . } \end{equation*}

On the other hand we can consider the pairing

\begin{equation} \label{6} \text{HH}_{i}(X) \times \text{HH}_{-i}(X)
\rar \fld \end{equation}

\begin{equation*} (v,w) \leadsto \text{tr}_{X}(D_{\Delta}(v) \circ w) \text{ . } \end{equation*}

\begin{prop} The pairings on Hochschild homology defined in \eqref{5}
and \eqref{6} are identical. \end{prop}

\begin{proof} By definition, $D_{\Delta}(v) = \calg J ({\calg T}(\calg
I(v))) $. The proposition would follow if we can check that
\begin{equation} \label{cal} \text{tr}_{X} (\calg J(\alpha) \circ
\beta) = \text{tr}_{X \times X}(\alpha \circ \calg I(\beta))
\end{equation} for any $\alpha \in \text{Hom}_{\text{D}^{b}(X \times
X)}(\strcc{\Delta}, \Delta_* S_X[i])$ and $\beta \in
\text{Hom}_{\dcat}(\strc[i], \Delta^* \strcc{\Delta})$. This is just
saying that $\calg I$ is the map "dual"
to the map $\calg J$ in \eqref{4}.  \\

We remark here that the assertion \eqref{cal} is similar to the last part of Proposition 3.1 of Caldararu's paper [4]. Proposition 3.1 of [4]
describes the construction of a right adjoint to a functor from $\dcat$ to $\text{D}^{b}(Y)$ given a left adjoint (via Serre duality). In our
situation, $\Delta_!$ is a left adjoint to $\Delta^*: \text{D}^{b}(X \times X) \rar \dcat$ . $\Delta_!$ was constructed in [4] using the right
adjoint $\Delta_*$
of $\Delta^*$ and Serre duality. \\

\end{proof}

Moreover, let $\int_X:  \text{H}^*(X) \rar \fld$ denote the linear functional that is $0$ on $\text{H}^{p,q}(X)$ if $(p,q) \neq (n,n)$ and
coincides with the identification of $\text{H}^n(X,\Omega_X^n) = \text{Hom}_{\dcat}(\strc,S_X)$ with $\fld$ on $\text{H}^{n,n}(X)$. Recall the
definition of the twisted HKR map $\thkr$ . The following proposition is immediate from the definition of $\thkr$.

\begin{prop} If $a \in \text{H}^*(X)$ and $b
\in \rhh_X(\strc, \Delta^* \strcc{\Delta})$, then
$$ \text{tr}_X(\thkr(a) \circ b) = \int_X \hkr(b) \wedge a \text{ . }$$
\end{prop}

Let $J$ be as in Theorem 1. Let $\langle, \rangle$ denote the Mukai
pairing on Hochschild homology {\it in this subsection only}. The
following proposition is immediate from Proposition 1,Proposition 2
and Theorem 1.

\begin{prop}
If $a \in \text{HH}_{i}(X)$ and $b \in \text{HH}_{-i}(X)$ then
$$\langle a,b \rangle = \int_X \hkr(b) \wedge J(\hkr(a)) \wedge \text{td}(T_X) \text{ . } $$
\end{prop}

Note that the product on $\text{H}^*(X)$ is graded commutative. Also note that $\int_X$ is non-vanishing only on $\text{H}^{2n}(X)$. Therefore,
$$\int_X v \wedge w = \int_X \bar{w} \wedge v $$ where $\bar{w}$ is obtained from $w$ by multiplying its component in $\text{H}^k(X)$ by
$(-1)^k$. Note that if $w \in \text{H}^*(X)$, $\bar{J(w)} = K(w)$ where $K$ is the endomorphism on $\text{H}^*(X)$ multiplying
$\text{H}^q(X,\Omega^p_X)$ by ${(-1)}^q$. Since $\text{td}(T_X) \in \oplus_i \text{H}^{2i}(X)$, $\text{td}(T_X)$ commutes with every element of
$\text{H}^*(X)$. Proposition 3 may therefore, be rewritten as  \\

\begin{equation} \label{pr3'}
\langle a,b \rangle = \int_X K(\hkr(a)) \wedge \hkr(b) \wedge
\text{td}(T_X) \text{ . }
\end{equation}

\subsection*{A Mukai like pairing on Hodge co-homology}
Now suppose that $X$ is a \\ {\it smooth complex variety}. Recall that a generalized Mukai pairing $\langle,\rangle_{C}$ has been defined by
Caldararu [3] on the Hodge cohomology $ \text{H}^*(X)$. Let $\omega_X = \Omega_X^n$ and let $\tau$ denote the endomorphism on $\text{H}^*(X)$
that is multiplication by $ {\sqrt {(-1)}}^{p+q} $ on $\text{H}^q(X,\Omega_X^p)$. Let $\text{ch}: \dcat \rar \text{H}^*(X)$ denote the Chern
character. Recall from [3] that $\sqrt{\text{ch}(\omega_X)}$ is a well-defined element of $\text{H}^*(X)$ . Then, if $v,w \in \text{H}^*(X)$,

\begin{equation} \langle v,w \rangle_{C} = \int_X
\frac{\tau(v)}{\sqrt{\text{ch}(\omega_X)}} \wedge w \text{ . }
\end{equation}

Let $\bar{\tau}$ denote the endomorphism on $\text{H}^*(X)$ that is multiplication by $ {\sqrt {(-1)}}^{q-p} $ on $\text{H}^q(X,\Omega_X^p)$.
Then, $K = \tau \circ \bar{\tau} = \bar{\tau} \circ \tau$. Define a pairing $\langle, \rangle_M$ on $\text{H}^*(X)$ by setting

\begin{equation} \langle v,w \rangle_M = \langle \bar{\tau}(v),w
\rangle_{C} = \int_X \frac{K(v)}{\sqrt{\text{ch}(\omega_X)}} \wedge
w \text{ . }
\end{equation}

\begin{prop}
If $f:X \rar Y$ is a proper morphism of smooth complex varieties,
then
$$ \langle f^*(v),w \rangle_M = \langle v,f_*(w) \rangle_M $$
for all $v \in \text{H}^*(Y)$ and $w \in \text{H}^*(X)$.
\end{prop}

\begin{proof} We recall from Caldararu [3] that $\langle f^*(v),w
\rangle_C = \langle v, f_*(w) \rangle_C$ for all $v \in  \text{H}^*(Y)$ and $w \in \text{H}^*(X)$ . Now, $f^*(\bar{\tau} v) =
\bar{\tau}(f^*(v))$ for any $v \in \text{H}^*(Y)$. Thus,
$$ \langle f^*(v),w \rangle_M = \langle \bar{\tau}(f^*(v)),w
\rangle_C= \langle f^*(\bar{\tau}(v)),w \rangle_C = \langle
\bar{\tau}(v),f_*(w) \rangle_C = \langle v,f_*(w) \rangle_M
$$ for all $v \in \text{H}^*(Y)$ and $w \in
\text{H}^*(X)$. \\
\end{proof}

\begin{prop}
If $a \in \text{HH}_{i}(X)$ and $b \in \text{HH}_{-i}(X)$ then
$$\langle a,b \rangle = \langle \text{ }
\hkr(a) \wedge \sqrt{\text{td}(T_X)},\hkr(b) \wedge
\sqrt{\text{td}(T_X)} \text{ }\rangle_M \text{ . }
$$
\end{prop}

\begin{proof} Since $\sqrt{\text{td}(T_X)}$ is a linear combination
of elements in $\text{H}^i(X, \Omega_X^i)$, it commutes with other
elements in $\text{H}^*(X)$. The RHS of the equation in this
proposition is therefore,
$$ \int_X \frac{K(\hkr(a) \wedge \sqrt{\text{td}(T_X)})}{\sqrt{\text{ch}(\omega_X)}} \wedge
\sqrt{\text{td}(T_X)} \wedge \hkr(b) \text{ . } $$ But $K$ is a ring
endomorphism of $\text{H}^*(X)$. Thus
$$  K(\hkr(a) \wedge \sqrt{\text{td}(T_X)}) =  K(\hkr(a)) \wedge
K(\sqrt{\text{td}(T_X)}) \text{ . }$$ But $K(\sqrt{\text{td}(T_X)})
= \tau(\sqrt{\text{td}(T_X)})$ since both $K$ and $\tau$ are
multiplication by ${(-1)}^i$ on $\text{H}^i(X, \Omega_X^i)$. It has
also been shown in Caldararu [3] that
$$\frac{\tau(\sqrt{\text{td}(T_X)})}{\sqrt{\text{ch}(\omega_X)}} =
\sqrt{\text{td}(T_X)} \text{ . }$$ It follows that
$$ \int_X \frac{K(\hkr(a) \wedge \sqrt{\text{td}(T_X)})}{\sqrt{\text{ch}(\omega_X)}} \wedge
\sqrt{\text{td}(T_X)} \wedge \hkr(b) = \int_X  K(\hkr(a)) \wedge
\text{td}(T_X) \wedge \hkr(b)$$ $$ = \int_X K(\hkr(a)) \wedge
\hkr(b) \wedge \text{td}(T_X) \text{ . }$$ The desired proposition
now follows from \eqref{pr3'}.

\end{proof}

\subsection*{Remark 1} Recall that if $\calg E \in \dcat$, then $\text{ch}(\calg E).\sqrt{\text{td}(T_X)}$ is
called the {\it Mukai vector} of $\calg E$ (Caldararu [3]). The
pairing $\langle,\rangle_M$ is slightly different from the
generalized Mukai pairing $\langle,\rangle_C$ defined by Caldararu
[3]. However, if $v$ and $w$ are Mukai vectors of elements of
$\dcat$, then $\langle v,w \rangle_M = \langle v,w \rangle_C$.

\subsection*{Remark 2} Let $X$ and $Y$ be smooth complex varieties.
Recall the definition of an integral transform $\Phi: \dcat \rar
\text{D}^{b}(Y)$ from Caldararu [3],[4]. An integral transform
$\Phi$ induces a map $\Phi_*: \text{H}^*(X) \rar \text{H}^*(Y)$. We
remark that the pairing $\langle,\rangle_M$ satisfies the {\it
Adjointness} one expects from a Mukai pairing. More precisely, if
$\Phi: \dcat \rar \text{D}^{b}(Y) $ and $\Psi:\text{D}^{b}(Y) \rar
\dcat$ are integral transforms such that $\Psi$ is a left adjoint of
$\Phi$, then
$$\langle \Psi_* v,w \rangle_M = \langle v,\Phi_*w \rangle_M $$ for
all $v \in \text{D}^{b}(Y)$ and $w \in \dcat$. This follows from the analogous property for the pairing $\langle,\rangle_C$ and the fact (see
Caldararu [3]) that integral transforms preserve the columns of the Hodge diamond. We are thus justified when we refer to the pairing
$\langle,\rangle_M$ as a {\textbf{Mukai like pairing}.

\subsection*{Remark 3} A part of the main conjecture of Caldararu
[3] was that the HKR map twisted by the square root of the Todd
genus of X preserves the Mukai pairing. However, instead of taking
the Mukai pairing on Hochschild homology to $\langle,\rangle_C$ , it
takes it to $\langle,\rangle_M$ by Proposition 5. The latter pairing
is itself a Mukai like pairing and is very similar to the former
pairing. However, $\langle,\rangle_M$ does not coincide with the
Mukai pairing on the Hodge cohomology of a K3-surface in general. In
particular, if $v \in \text{H}^{2,0}(X)$ and $w \in
\text{H}^{0,2}(X)$ then $\langle v,w \rangle_M \neq \langle v,w
\rangle_C$. This is why we do not go so far as to call
$\langle,\rangle_M$ a generalized Mukai pairing. We can however,
justifiably say that the HKR map twisted by the square root of the
Todd genus of X "almost preserves" the Mukai pairing.

\subsection{The relative Riemann-Roch theorem}
The relative Riemann-Roch theorem follows from Proposition 5. Recall
that Caldararu [4] defined a Chern character $$\text{Ch}: \dcat \rar
\text{HH}_0(X) \text{ . }$$ He also showed in [4] that if $f:X \rar
Y$ is a proper morphism of smooth schemes, then
\begin{equation} \label{9} \text{Ch}(f_* \calg E) = f_* \text{Ch}(\calg E)
\end{equation}
for any $\calg E \in \dcat$. Also, in [3], it was shown that $$\hkr
\circ \text{Ch}(\calg E) = \text{ch}(\calg E) \text{ . }
$$
We now have
\begin{thmrr}
Let $f:X \rar Y$ be a proper morphism of smooth proper complex
varieties. Then, if $\calg E$ is a vector bundle on $X$, $$\int_X
f^*(l) \text{ch}(\calg E) \text{td}(T_X) = \int_Y l \text{ch}(f_*
\calg E) \text{td}(T_Y) $$ for any $l \in \text{H}^*(Y)$.
\end{thmrr}

\begin{proof}
Note that $\hkr: \text{HH}_{*}(X) \rar \text{H}^*(X)$ is an
isomorphism of complex vector spaces . Let \\ $a = \hkr^{-1}(K(l))
\in \text{HH}_*(Y)$. Then,
\begin{equation} \label{10} \langle f^*a, \text{Ch}(\calg E) \rangle = \langle a, f_*
\text{Ch}(\calg E) \rangle = \langle a,  \text{Ch}(f_* \calg E)
\rangle \text{ . } \end{equation}

The first equality in \eqref{10} is due to the Adjointness property
of the Mukai pairing (see Caldararu [4]). The second equality in
\eqref{10} is due to \eqref{9}. By Proposition 5 and the fact that
$\text{td}(T_X)$ commutes with other elements of $\text{H}^*(X)$ ,
\begin{equation} \label{11} \langle f^*a, \text{Ch}(\calg E) \rangle
= \int_X K(\hkr(f^*a)) \hkr(\text{Ch}(\calg E)) \text{td}(T_X)
\end{equation}

$$\langle a,  \text{Ch}(f_* \calg E)\rangle = \int_Y K(\hkr(a)) \hkr(\text{Ch}(f_* \calg
E)) \text{td}(T_Y) \text{ . }$$

Now note that $\hkr \circ f^* = f^* \circ \hkr $ and that $K \circ f^* = f^* \circ K$.
 Also $\hkr \circ \text{Ch} = \text{ch}$. Now, applying
these facts to \eqref{11} and using \eqref{10}, the desired theorem follows.

\end{proof}

\subsection*{Remark} Note that the Chern character to Hochschild
homology actually commutes with push-forwards as shown in Caldararu
[4]. The Todd genus in the relative Riemann-Roch theorem thus occurs
as a consequence of the fact that the Mukai pairing on Hochschild
homology does not correspond to a Mukai like pairing on Hodge
cohomology under $\hkr$. For the Mukai pairing on Hochschild
homology to be "preserved" in any sense, one has to twist $\hkr$ by
$\sqrt{\text{td}(T_X)}$.

\section{Two "connections" on the Hochschild chain complex $\hc^{\bullet}(X)$}

\subsection{The completed bar and Hochschild chain complexes}

 Let $U = \text{Spec }R$ be an open affine subscheme of
$X$. The restriction of $\strcc{\Delta}$ to $U \times U$ has a free
$R
\otimes R$-module resolution given by the {\it Bar resolution} : \\

$$B^{-n}(R) = R^{\otimes n+2} $$

$$d(r_0 \otimes ... \otimes r_{n+1}) = r_0r_1 \otimes .... \otimes
r_{n+1} -r_0 \otimes r_1r_2 \otimes ... \otimes r_{n+1} +..... +
{(-1)}^{n} r_0 \otimes ... \otimes r_nr_{n+1} \text{    } r_i \in
R$$.

The $R \otimes R$-module structure is given by multiplication with
the extreme factors. Let $I_n$ denote the kernel of the $n+2$-fold
multiplication $R^{\otimes n+2} \rar R$. Let $$\hb^{-n}(R) = \ilim_k
\frac{B^{-n}(R)}{I_n^k} \text{ . } $$ Note that each summand of the
differential $d$ takes $I_n$ to $I_{n-1}$. It follows that the
differential on $B^{\bullet}(R)$ extends to yield a differential on
$\hb^{\bullet}(R)$. Yekutieli [5] shows that completing the bar
resolution in this manner yields a complex $\hb^{\bullet}(X)$ of
coherent sheaves on $X \times X$. He also shows that
$\hb^{\bullet}(X)$ is a resolution of $\strcc{\Delta}$ by flat
$\strcc{X \times X}$-modules. \\

It follows that $\Delta^* \strcc{\Delta}$ is represented by the
complex $$\hc^{\bullet}(X):=  \strc \otimes_{\Delta^{-1} \strcc{X
\times X}} \Delta^{-1} \hb^{\bullet}(X) \text{ . }
$$ $\hc^{\bullet}(X)$ is called the complex of completed {\it Hochschild chains} on $X$.\\

On an open subscheme $U = \text{Spec }R$ of $X$ before completion,
$$C^{-n}(R) = R^{\otimes n+1} $$
$$d(r_0 \otimes .... \otimes r_n) = r_0r_1 \otimes .... \otimes r_n
-r_0 \otimes r_1r_2 \otimes .... \otimes r_n + ... $$
$$...+{(-1)}^{n-1}r_0 \otimes ... \otimes r_{n-1}r_n +{(-1)}^n
r_nr_0 \otimes .... \otimes r_{n-1} $$
$$\hc^{-n}(R) = \ilim_k \frac{B^{-n}(R)}{I_n^k}
\otimes_{R^{\otimes n+2}} C^{-n}(R) \text{ . } $$

Yekutieli [5] also showed that $\text{RD}(\hc^{\bullet}(X))$ is
represented in $\dcat$ by the complex $\dd^{\bullet}(X)$ of
poly-differential operators on $X$ equipped with Hochschild
co-boundary. \\

Let us describe some operations on $\hc^{\bullet}(X)$ that endow it
with the structure of a Hopf algebra in $\bacc$ ( and therefore in $\dcat$ ) . \\

\subsubsection*{Product on $\hc^{\bullet}(X)$ :} The product $m:\hc^{\bullet}(X)
\otimes_{\strc} \hc^{\bullet}(X) \rar \hc^{\bullet}(X)$ is given by
the {\it signed shuffle product}. Recall that a $(p,q)$-shuffle
$\sigma$ is a permutation of $\{1,....,p+q\}$ such that
$\sigma(1)<...<\sigma(p)$ and $\sigma(p+1) <...< \sigma(p+q)$.
Denote the set of $(p,q)$-shuffles by $\shuf{p}{q}$.  On an open
subscheme $U= \text{Spec } R$ of $X$ before completion, this product
is given by

$$(r_0 \otimes r_1 \otimes ....\otimes r_p ) \otimes_R (r'_0 \otimes
r_{p+1} \otimes ... \otimes r_{p+q}) \leadsto \sum_{\sigma \in
\shuf{p}{q} } \text{sgn}(\sigma) r_0r'_0 \otimes r_{\sigma^{-1}(1)}
\otimes .... \otimes r_{\sigma^{-1}(p+q)} $$.

This is easily seen to be a (graded) commutative product.

\subsubsection*{Co-product on $\hc^{\bullet}(X)$ } The co-product $\hc^{\bullet}(X) \rar  \hc^{\bullet}(X)
\otimes_{\strc} \hc^{\bullet}(X)$ is given by the {\it cut co-product }. Contrary to the usual practise, we denote the co-product by
$\textbf{C}$ to avoid confusion with $\Delta$ which denotes the diagonal map $X \rar X \times X$  in this paper. On an open subscheme $U =
\text{Spec } R $ of $X$ before completion,
$$\textbf{C}(r_0 \otimes ... \otimes r_n) = \sum_{p+q=n} r_0 \otimes
r _1 \otimes .. \otimes r_p \otimes_R 1 \otimes r_{p+1} \otimes ...
\otimes r_n \text{ . }$$

\subsubsection*{Unit for $\hc^{\bullet}(X)$ :} There is a unit map
$\epsilon: \strc \rar \hc^{\bullet}(X)$. On an open subscheme \\ $U
= \text{Spec }R$ of $X$ before completion , $\epsilon$ is given by
the composite
$$ R \simeq C^0(R) \hookrightarrow C^{\bullet}(R) \text{ . }$$

\subsubsection*{Co-unit for $\hc^{\bullet}(X)$ : } There is a co-unit
$\eta: \hc^{\bullet}(X) \rar \strc$. On an open subscheme  $U=
\text{Spec }R$ of $X$ before completion, this is given by the
projection from $C^{\bullet}(R)$ to $C^0(R)$.

\begin{prop} $M,\textbf{C},\epsilon$ and $\eta$ make
$\hc^{\bullet}(X)$ a Hopf algebra in $\bacc$.
\end{prop}

\begin{proof} It has been proven by Yekutieli[5] that applying the functor
$\hhh_{\strc}^{\text{cont}}(-,\strc)$ (with $\strc$ given the
discrete topology ) to $\hc^{\bullet}(X) \in \bacc$ yields
\\ $\dd^{\bullet}(X) \in \bbcc$. It is easy to verify that this
functor takes the product,co-product,unit and co-unit of
$\hc^{\bullet}(X)$ to the co-product, product,co-unit and unit of
$\dd^{\bullet}(X)$ respectively. The desired proposition then
follows from Proposition
2 of [1]. \\

\end{proof}

\subsubsection*{Antipode on $\hc^{\bullet}(X)$ : }
$\hc^{\bullet}(X)$ also comes equipped with an {\it Antipode } map.
We will denote this map by $S$. On $U = \text{Spec } R$ before
completion,
$$S(r_0 \otimes ... \otimes r_n)= {(-1)}^{\frac{n(n+1)}{2}}r_0
\otimes r_n \otimes r_{n-1} \otimes ..... \otimes r_1 \text{ . }$$

\subsection*{The Hochschild-Kostant-Rosenberg (HKR) map :} There is
a quasi-isomorphism (see Yekutieli [5]) $$ \hkr: \hc^{\bullet}(X)
\rar \sss(\Omega[1])$$ of complexes of $\strc$-modules. On an open
subscheme $U = \text{Spec }R$ of $X$ before completion, $$\hkr(r_0
\otimes .... \otimes r_n ) = \frac{1}{n!} r_0 dr_1 \wedge ... \wedge
dr_n \text{ . }$$

\subsection{Two connections on $\hc^{\bullet}(X)$ :}
Let $\pi_k:\sss(\Omega[1]) \rar \Omega^k[k]$ denote the natural
projection. Denote by $\alpha_R$ the composite
$$\begin{CD} \hc^{\bullet}(X) @>\textbf{C}>> \hc^{\bullet}(X)
\otimes \hc^{\bullet}(X) @>(\textbf{1}_{\hc^{\bullet}(X)} \otimes
\pi_1 \circ \hkr) >> \hc^{\bullet}(X) \otimes \Omega[1] \\ \end{CD}
\text{ . }$$

More concretely, on an open subscheme $U= \text{Spec } R$ of $X$
before completion,

$$ \alpha_R(r_0 \otimes .... \otimes r_n) = r_0 \otimes ... \otimes
r_{n-1} \otimes dr_n \text{ . } $$

Let $\alpha_L: \hc^{\bullet}(X) \rar \hc^{\bullet}(X) \otimes
\Omega[1]$ be the map such that

$$ \alpha_L(r_0 \otimes .... \otimes r_n) = {(-1)}^{n-1} r_0 \otimes r_2 \otimes
..... \otimes r_n \otimes dr_1 $$ on any open subscheme $U=
\text{Spec } R$ of $X$ before completion.

Then, $\alpha_L = - (S \otimes \Omega[1]) \circ  \alpha_R \circ S$. Let $\alpha_R \otimes \hc^{\bullet}(X)$ denote the composite $$ \begin{CD}
\hc^{\bullet}(X) \otimes \hc^{\bullet}(X) @> \alpha_R \otimes \hc^{\bullet}(X) >> \hc^{\bullet}(X) \otimes \Omega[1] \otimes \hc^{\bullet}(X)
@>\hc^{\bullet}(X) \otimes \tau >> \hc^{\bullet}(X) \otimes \hc^{\bullet}(X) \otimes \Omega[1] \end{CD} $$ where $\tau: \Omega[1] \otimes
\hc^{\bullet}(X) \rar \hc^{\bullet}(X) \otimes \Omega[1]$ is the map that swaps factors. Similarly, let $\alpha_L \otimes \hc^{\bullet}(X)$
denote the composite $$
\begin{CD} \hc^{\bullet}(X) \otimes \hc^{\bullet}(X) @> \alpha_L \otimes \hc^{\bullet}(X) >> \hc^{\bullet}(X) \otimes \Omega[1] \otimes
\hc^{\bullet}(X) @>\hc^{\bullet}(X) \otimes \tau >> \hc^{\bullet}(X)
\otimes \hc^{\bullet}(X) \otimes \Omega[1] \end{CD} \text{ . } $$
 We now have the
following proposition: \\

\begin{prop} The following diagrams commute in $\bacc$:
\begin{equation} \label{bbb} \begin{CD} \hc^{\bullet}(X) \otimes \hc^{\bullet}(X)
@>\alpha_R \otimes \hc^{\bullet}(X) + \hc^{\bullet}(X) \otimes
\alpha_R >>
\hc^{\bullet}(X) \otimes \hc^{\bullet}(X) \otimes \Omega[1] \\
@VVmV                         @Vm \otimes \Omega[1]VV \\
\hc^{\bullet}(X) @> \alpha_R >> \hc^{\bullet}(X) \otimes \Omega[1]
\\
\end{CD} \end{equation}

\begin{equation} \label{cc} \begin{CD} \hc^{\bullet}(X) \otimes \hc^{\bullet}(X)
@>\alpha_L \otimes \hc^{\bullet}(X) + \hc^{\bullet}(X) \otimes
\alpha_L >>
\hc^{\bullet}(X) \otimes \hc^{\bullet}(X) \otimes \Omega[1] \\
@VVmV                         @Vm \otimes \Omega[1]VV \\
\hc^{\bullet}(X) @> \alpha_L >> \hc^{\bullet}(X) \otimes \Omega[1]
\\
\end{CD} \end{equation}
\end{prop}

\begin{proof}
This proposition is proven by a combinatorial argument. On an open
subscheme $U = \text{Spec } R$ of $X$ before completion,

\begin{equation} \label{aa}  (\alpha_R \otimes \hc^{\bullet}(X) + \hc^{\bullet}(X) \otimes
\alpha_R)((r_0 \otimes r_1 \otimes ... \otimes r_p) \otimes_R (1
\otimes r_{p+1} \otimes ... \otimes r_{p+q})) = $$ $${(-1)}^q (r_0
\otimes r_1 \otimes ... \otimes r_{p-1}) \otimes_R (1 \otimes
r_{p+1} \otimes ... \otimes r_{p+q}) \otimes dr_p + $$ $$ (r_0
\otimes r_1 \otimes ... \otimes r_p) \otimes_R (1 \otimes r_{p+1}
\otimes ... \otimes r_{p+q-1}) \otimes dr_{p+q} \text{ . } \end{equation}\\

Note that if $\sigma$ is a $(p,q)$-shuffle, then $\sigma^{-1}(p+q)=
p$ or $\sigma^{-1}(p+q) = p+q $. Let $\shuf{p}{q}^1$ denote the set
of all $(p,q)$-shuffles $\sigma$ such that $\sigma^{-1}(p+q) = p$.
Let $\shuf{p}{q}^2$ denote $\shuf{p}{q} \setminus \shuf{p}{q}^1$. \\

Also note that there is a sign preserving bijection from
$\shuf{p}{q-1}$ to $\shuf{p}{q}^2$ the inverse of which takes an
element $\sigma$ of $\shuf{p}{q}^2$ to its restriction to the set
$\{1,...,p+q-1\}$.
Denote this bijection by $I:\shuf{p}{q-1} \rar \shuf{p}{q}^2 $. \\

 We also have a bijection from $\shuf{p}{q}^1$ to
$\shuf{p-1}{q}$. Let \\ $\psi:\{1,..,p+q-1\} \rar
\{1,...,p-1,p+1,.....,p+q-1\}$ be the unique order-preserving
map.The permutation in $\shuf{p-1}{q}$ corresponding to a
permutation $\sigma$ in $\shuf{p}{q}^1$ is given by the composite
$$\begin{CD} \{1,..,p+q-1\} @>\psi>> \{1,..,p-1,p+1,...,p+q\}
@>\sigma>> \{1,..,p+q-1\} \end{CD} \text{ . }$$ This bijection from
$\shuf{p}{q}^1$ to $\shuf{p-1}{q}$ however, changes the sign by
${(-1)}^q$. Denote the inverse of this bijection by $J:\shuf{p-1}{q}
\rar
\shuf{p}{q}^1$. \\

Then, $$m((r_0 \otimes r_1 \otimes ... \otimes r_p) \otimes_R (1
\otimes r_{p+1} \otimes ... \otimes r_{p+q})) = \sum_{\sigma \in
\shuf{p}{q}} \text{sgn}(\sigma) r_0 \otimes r_{\sigma^{-1}(1)}
\otimes ... r_{\sigma^{-1}(p+q)} $$ $$ = \sum_{\sigma \in
\shuf{p}{q}^1} \text{sgn}(\sigma) r_0 \otimes r_{\sigma^{-1}(1)}
\otimes ... r_{\sigma^{-1}(p+q)} + \sum_{\sigma \in
\shuf{p}{q}^2}\text{sgn}(\sigma)  r_0 \otimes r_{\sigma^{-1}(1)}
\otimes ... r_{\sigma^{-1}(p+q)} $$  $$= \sum_{\sigma \in
\shuf{p-1}{q}} {(-1)}^q \text{sgn}(\sigma) r_0 \otimes
r_{J(\sigma)^{-1}(1)} \otimes .... \otimes r_{J(\sigma)^{-1}(p+q-1)}
\otimes r_p  $$ $$ +  \sum_{\sigma \in \shuf{p}{q-1}}
\text{sgn}(\sigma)r_0 \otimes r_{I(\sigma)^{-1}(1)} \otimes ....
\otimes r_{I(\sigma)^{-1}(p+q-1)} \otimes r_{p+q} $$ $$ ={(-1)}^q
m((r_0 \otimes r_1 \otimes ... \otimes r_{p-1}) \otimes_R (1 \otimes
r_{p+1} \otimes ... \otimes r_{p+q})) \otimes r_p + $$ $$ m((r_0
\otimes r_1 \otimes ... \otimes r_p) \otimes_R (1 \otimes r_{p+1}
\otimes ... \otimes r_{p+q-1})) \otimes r_{p+q} \text{ . } $$

It follows that $$\alpha_R(m((r_0 \otimes r_1 \otimes ... \otimes
r_p) \otimes_R (1 \otimes r_{p+1} \otimes ... \otimes r_{p+q}))) =
$$ $$ ={(-1)}^q
m((r_0 \otimes r_1 \otimes ... \otimes r_{p-1}) \otimes_R (1 \otimes
r_{p+1} \otimes ... \otimes r_{p+q})) \otimes dr_p + $$ $$ m((r_0
\otimes r_1 \otimes ... \otimes r_p) \otimes_R (1 \otimes r_{p+1}
\otimes ... \otimes r_{p+q-1})) \otimes dr_{p+q} \text{ . } $$

It follows from \eqref{aa} that this is precisely $$(m \otimes
\Omega[1]) ((\alpha_R \otimes \hc^{\bullet}(X) + \hc^{\bullet}(X)
\otimes \alpha_R)) ((r_0 \otimes r_1 \otimes ... \otimes r_p)
\otimes_R (1 \otimes r_{p+1} \otimes ... \otimes r_{p+q})) \text{ .
}
$$

This proves that the diagram \eqref{bbb} commutes. Proving that the diagram \eqref{cc}
commutes is very similar and left to the reader. \\

\end{proof}

Let $\alpha_R^{\circ i}$ denote the composite $$ \begin{CD}
\hc^{\bullet}(X)@>\alpha_R>> \hc^{\bullet}(X) \otimes \Omega[1]
@>\alpha_R \otimes \Omega[1] >> \hc^{\bullet}(X) \otimes
\Omega[1]^{\otimes 2} @>>> ... @>\alpha_R \otimes \Omega[1]^{\otimes
i-1} >> \hc^{\bullet}(X) \otimes \Omega[1]^{\otimes i} \end{CD} $$.

Let $p:\Omega[1]^{\otimes i} \rar \Omega^i[i]$ be the standard
projection. On an open subscheme $U = \text{Spec } R$, of $X$,
$$p(r_0dr_1 \otimes ... \otimes dr_i) = r_0 dr_1 \wedge .... \wedge
dr_i \text{ . }$$

Let $\alpha_R^i$ denote the composite
$$\begin{CD} \hc^{\bullet}(X) @> \alpha_R^{\circ i} >>
\hc^{\bullet}(X) \otimes \Omega[1]^{\otimes i} @> \hc^{\bullet}(X)
\otimes p >> \hc^{\bullet}(X) \otimes \Omega^i[i] \end{CD} \text{ .
}
$$

Let $\text{exp}(\alpha_R)$ denote the sum $$\sum_i \frac{1}{i!}
\alpha_R^i : \hc^{\bullet}(X) \rar \hc^{\bullet}(X) \otimes
\sss(\Omega[1]) \text{ . } $$

We then have the following proposition:
\begin{prop}
$$ (\hc^{\bullet}(X) \otimes \hkr) \circ \textbf{C} =
\text{exp}(\alpha_R) \text{ . }$$
\end{prop}

\begin{proof}
On an open subscheme $U = \text{Spec } R$ before completion,
$$\text{exp}(\alpha_R)(r_0 \otimes ... \otimes r_n) =\frac{1}{i!} \sum_i r_0
\otimes ... \otimes r_{n-i} \otimes_R dr_{n-i+1} \wedge .. \wedge
dr_n $$ $$ = \sum_i r_0 \otimes ... \otimes r_{n-i} \otimes_R \hkr(1
\otimes r_{n-i+1} \otimes ... \otimes r_n) $$ $$= (\hc^{\bullet}(X)
\otimes \hkr) ( \sum_i r_0 \otimes ... \otimes r_{n-i} \otimes_R 1
\otimes r_{n-i+1} \otimes ... \otimes r_n ) $$ $$= (\hc^{\bullet}
\otimes \hkr) \circ \textbf{C} (r_0 \otimes .... \otimes r_n) \text{
. }$$ This verifies the desired proposition.

\end{proof}

Let $\tau: \Omega[1] \otimes \Omega[1] \rar \Omega[1] \otimes \Omega[1]$ denote the swap map. The following proposition
tells us that $\alpha_L$ and $\alpha_R$ "commute" with each other. \\

\begin{prop}
$$(\alpha_R \otimes \Omega[1]) \circ \alpha_L - (\hc^{\bullet}(X) \otimes \tau) \circ
(\alpha_L \otimes \Omega[1]) \circ \alpha_R = 0 \text{ . }$$
\end{prop}

\begin{proof}
On an open subscheme $U = \text{Spec }R$ before completion,
$$(\alpha_R \otimes \Omega[1]) \circ \alpha_L(r_0 \otimes ..... \otimes r_n) =
{(-1)}^{n-1} (r_0 \otimes r_2 \otimes ... \otimes r_{n-1}) \otimes_R
dr_n \otimes_R dr_1 $$
$$(\alpha_L \otimes \Omega[1]) \circ \alpha_R (r_0 \otimes ..... \otimes r_n) =
{(-1)}^{n-2} (r_0 \otimes r_2 \otimes ... \otimes r_{n-1}) \otimes_R
dr_1 \otimes_R dr_n $$
$$(\hc^{\bullet}(X) \otimes \tau)((r_0 \otimes r_2 \otimes ... \otimes r_{n-1}) \otimes_R
dr_1 \otimes_R dr_n ) = -(r_0 \otimes r_2 \otimes ... \otimes
r_{n-1}) \otimes_R dr_n \otimes_R dr_1 $$. The desired proposition
is
now immediate. \\
\end{proof}

\subsection{Remark - the beginning of a dictionary} For reasons that
will become clear later in this section, the reader should think of
$\hc^{\bullet}(X)$ as analogous to the ring of functions on an open
"symmetric" neighborhood $U_G$ of the identity of a Lie group $G$.
By "symmetric" we mean that if $g \in U_G$ then $g^{-1} \in U_G$.
$T[-1]$ is the analog of the Lie algebra $\mathfrak{g}$ of the Lie
group $G$. Thus, $\Omega[1]$ is the analog of $\mathfrak{g}^*$.
Proposition 7 says that in this picture, both $\alpha_L$ and
$\alpha_R$ are analogs of "connections" on the
ring of functions of $G$. \\

In the same picture, $\sss(\Omega[1])$ is to be thought of as analogous to the ring of functions on a neighborhood $\calg V$ of $0$ in
$\mathfrak{g}$ that is diffeomorphic to $U_G$ via the exponential map . The
Hochschild-Kostant-Rosenberg map is then the analog of the pull back by the exponential map $\text{exp}^*$. \\

The antipode map $S$ is the analog of the pull-back by the map which takes an element of $G$ to its inverse. \\

\subsection{More on the maps $\alpha_R$ and $\alpha_L$ }

The question that arises at this stage is "Can $\alpha_L$ and
$\alpha_R$ be described by explicit formulae as maps in $\dcat$ from
$\sss(\Omega[1])$ to $\sss(\Omega[1]) \otimes \Omega[1]$ ?".
Markarian sought to answer this in Theorem 1 of [2]. The proof there
was however
 erroneous. A result dual to what we want is available in [1]. A later version [6] of
 [2] also contains a result equivalent to Theorem 1 of [2]. \\

Recall from Kapranov [7] that $T[-1]$ is a Lie algebra object in
$\dcat$. The Lie bracket of $T[-1]$ is given by the Atiyah class
$$\text{At}_T:T[-1] \otimes T[-1] \rar T[-1] $$ of the tangent bundle
of $X$. It is also known that the universal enveloping algebra of
$T[-1]$ in $\dcat$ is represented by the complex $\dd^{\bullet}(X)$.
This was proven in [1]. Equivalent results have been proven using
methods different from that in [1] by Markarian [6] and Roberts and Willerton [8]. \\

Let $\mu$ denote the wedge product on $\sss(T[-1])$. Let \\
 $\delta: \sss(T[-1]) \otimes T[-1] \rar  \sss(T[-1]) \otimes T[-1] \otimes T[-1]$
be the map
$$\delta(v_1 \wedge... \wedge v_k \otimes y)= \sum_{i=1}^{i=k}
{(-1)} ^{k-i} \widehat{v_1 \wedge .... i ... \wedge v_k} \otimes v_i
\otimes y $$ for sections $v_1,..,v_k,y$ of $T$ over an open
subscheme $U= \text{Spec } R$ of $X$. \\

We have a map
$$\bar{\omega}:\sss(T[-1]) \otimes T[-1] \rar \sss(T[-1]) \otimes
T[-1] $$
$$\bar{\omega} = (\textbf{1}_{\sss(T[-1])} \otimes \text{At}_{T} ) \circ
\delta \text{ . }$$ Note that $\mu \circ \bar{\omega}$ yields the
right
adjoint action of the Lie algebra object $T[-1]$ on $\sss(T[-1])$.\\

The Hochschild-Kostant-Rosenberg map $\hkr:\sss(T[-1]) \rar
\dd^{\bullet}(X) $ is the "dual" of the HKR map
$\hkr:\hc^{\bullet}(X) \rar \sss(\Omega[1])$. The following theorem,
which figures as Corollary 1 in [1] and Theorem 2 in [6], describes
the Duflo-like error term that measures how the map
$\hkr:\sss(T[-1]) \rar \dd^{\bullet}(X)$ fails to commute with
multiplication. \\

\begin{thm}
(Recalled from [1]. [6] has a similar result) The following diagram commutes in $\dcat$ :\\
$$\begin{CD} \dd^{\bullet}(X) \otimes \dd^{\bullet}(X) @> m >> \dd^{\bullet}(X) \\
@AA \hkr \otimes \hkr A            @A \hkr AA \\
\sss(T[-1]) \otimes T_[-1] @> \mu \circ
\frac{\bar{\omega}}{1-e^{-\bar{\omega}}} >> \sss(T[-1])  \\
\end{CD}$$ \\
\end{thm}

Note that applying the functor $\text{RD}$ to $\bar{\omega}$ gives
us a morphism in $\dcat$ $$\bar{\omega}: \sss(\Omega[1]) \otimes
\Omega[1] \rar \sss(\Omega[1]) \otimes \Omega[1] \text{ . }$$
Further, denote the co-multiplication on $\sss(\Omega[1])$ by
$\textbf{C}_{\Omega}$. Denote the map $$
(\sss(\Omega[1]) \otimes \pi_1) \circ \textbf{C}_{\Omega} : \sss(\Omega[1]) \rar \sss(\Omega[1]) \otimes \Omega[1] $$ by $\bar{\textbf{C}}$.\\

We denote the map $$\frac{\bar{\omega}}{e^{\bar{\omega}}-1} \circ \bar{\textbf{C}}: \sss(\Omega[1]) \rar \sss(\Omega[1]) \otimes \Omega[1]$$ by
$\Phi_L$. The map $$\frac{\bar{\omega}}{1-e^{-\bar{\omega}}} \circ \bar{\textbf{C}}: \sss(\Omega[1]) \rar \sss(\Omega[1]) \otimes \Omega[1]$$
will be denoted by $\Phi_R$.

We can now state the following theorem. This can be thought of as the starting point for the computations leading to Theorem 1.

\begin{thm2p}
The following diagrams commute in $\dcat$ : \\
\begin{equation}
\label{dd}
\begin{CD}
\hc^{\bullet}(X) @>\alpha_R>> \hc^{\bullet}(X) \otimes \Omega[1] \\
@VV \hkr V                      @V \hkr \otimes \Omega[1] VV \\
\sss(\Omega[1]) @> \Phi_R >> \sss(\Omega[1]) \otimes \Omega[1] \\
\end{CD}
\end{equation}

\begin{equation}
\label{ee}
\begin{CD}
\hc^{\bullet}(X) @>\alpha_L>> \hc^{\bullet}(X) \otimes \Omega[1] \\
@VV \hkr V                      @V \hkr \otimes \Omega[1] VV \\
\sss(\Omega[1]) @> \Phi_L >> \sss(\Omega[1]) \otimes \Omega[1] \\
\end{CD}
\end{equation}
\end{thm2p}

\begin{proof}
The fact that the diagram \eqref{dd} commutes in $\dcat$ is obtained
by applying the functor $\text{RD}$ to the diagram in Theorem 2. \\

Note that since $\alpha_L = - (S \otimes \Omega[1]) \circ  \alpha_R
 \circ S$ the following diagram
commutes in \\ $\bacc$ (and therefore in $\dcat$ ): \\

\begin{equation}
\label{fa}
\begin{CD}
\hc^{\bullet}(X) @> \alpha_R>> \hc^{\bullet}(X) \otimes \Omega[1] \\
@VV S V    @V S \otimes \Omega[1] VV \\
\hc^{\bullet}(X) @>-\alpha_L>> \hc^{\bullet}(X) \otimes \Omega[1] \\
\end{CD}
\text{ . }
\end{equation}

Recall that $J$ is the endomorphism of $\sss(\Omega[1])$ multiplying
$\exal{ }{i}$ by ${(-1)}^i$. Then,\begin{equation} \label{ga} \hkr
\circ S = J \circ \hkr \text{ . } \end{equation}

To see \eqref{ga}, note that on an open subscheme $U = \text{Spec }
R$ of $X$ before completion, $$\hkr \circ S(r_0 \otimes .... \otimes
r_n) = {(-1)}^{\frac{n(n+1)}{2}} \frac{1}{n!} r_0 dr_n \wedge ...
\wedge dr_1 =
$$ $$ {(-1)}^{\frac{n(n+1)}{2}} {(-1)}^{\frac{n(n-1)}{2}} \frac{1}{n!} r_0 dr_1 \wedge
... \wedge dr_n = {(-1)}^n \hkr (r_0 \otimes ... \otimes r_n) \text{
. }
$$

Further, note that applying the functor $\text{RD}$ to $J$ yields
the endomorphism $I$ of $\sss(T[-1])$ that multiplies $\texal{X}{i}$
by ${(-1)}^i$. For sections $v_1,...,v_k,y$ of $T$ over an open
subscheme $U = \text{Spec } R$ of $X$,

$$\delta(I(v_1 \wedge .... \wedge v_k) \otimes y) = {(-1)}^ k \sum_{i=1}^{i=k}
{(-1)} ^{k-i} \widehat{v_1 \wedge .... i ... \wedge v_k} \otimes v_i
\otimes y $$ $$ ={(-1)}^k {(-1)}^{k-1} \sum_{i=1}^{i=k} {(-1)}
^{k-i} I(\widehat{v_1 \wedge .... i ... \wedge v_k}) \otimes v_i
\otimes y \text{ . }$$

It follows from the above computation and the fact that $\bar{\omega} = (\sss(T[-1]) \otimes \text{At}_T) \circ \delta$ that the following
 diagram commutes in $\dcat$: \\

\begin{equation} \label{ha}
\begin{CD}
\sss(T[-1]) \otimes T[-1] @>{\bar{\omega}}>> \sss(T[-1]) \otimes
T[-1]\\
@VV I \otimes T[-1] V    @VV I \otimes T[-1] V \\
\sss(T[-1]) \otimes T[-1] @>-\bar{\omega}>> \sss(T[-1]) \otimes
T[-1]\\
\end{CD}
\end{equation}

Applying the functor $\text{RD}$ to the diagram \eqref{ha} we obtain
the following diagram :

\begin{equation} \label{ka}
\begin{CD}
\sss(\Omega[1]) \otimes \Omega[1] @>{\bar{\omega}}>> \sss(\Omega[1])
\otimes
\Omega[1]\\
@VV J \otimes \Omega[1] V   @V J \otimes \Omega[1] VV \\
\sss(\Omega[1]) \otimes \Omega[1] @>-{\bar{\omega}}>>
\sss(\Omega[1]) \otimes
\Omega[1]\\
\end{CD}
\end{equation}

A calculation similar to the one made to verify \eqref{ha} also
shows that the following diagram commutes : \\

\begin{equation} \label{la}
\begin{CD}
\sss(\Omega[1])  @>\bar{\textbf{C}}>> \sss(\Omega[1]) \otimes
\Omega[1]\\
@VV J  V   @V J \otimes \Omega[1] VV \\
\sss(\Omega[1])  @>-\bar{\textbf{C}}>> \sss(\Omega[1]) \otimes
\Omega[1]\\
\end{CD}
\end{equation}

Combining \eqref{la} and \eqref{ka} we obtain the following
commutative diagram :\\

\begin{equation} \label{ma}
\begin{CD}
\sss(\Omega[1])  @>\frac{\bar{\omega}}{1-e^{-\bar{\omega}}} \circ
\bar{\textbf{C}}>> \sss(\Omega[1]) \otimes
\Omega[1]\\
@VV J  V   @V J \otimes \Omega[1] VV \\
\sss(\Omega[1])  @>-\frac{\bar{\omega}}{e^{\bar{\omega}}-1} \circ
\bar{\textbf{C}}>> \sss(\Omega[1]) \otimes
\Omega[1]\\
\end{CD}
\end{equation}

It follows from \eqref{dd} and \eqref{ma} that all squares in the
diagram below commute in $\dcat$. \\

\begin{equation} \label{na} \begin{CD}
\hc^{\bullet}(X) @>\alpha_R>> \hc^{\bullet}(X) \otimes \Omega[1] \\
@VV \hkr V                      @V \hkr \otimes \Omega[1] VV \\
\sss(\Omega[1]) @> \frac{\bar{\omega}}{1-e^{-\bar{\omega}}} \circ
\bar{\textbf{C}} >> \sss(\Omega[1]) \otimes \Omega[1] \\
@VV J V   @V J \otimes \Omega[1] VV \\
\sss(\Omega[1]) @>-\frac{\bar{\omega}}{e^{\bar{\omega}}-1} \circ
\bar{\textbf{C}}>> \sss(\Omega[1]) \otimes \Omega[1] \\
\end{CD} \end{equation}

The diagram \eqref{na} says that $$-\Phi_L \circ J \circ \hkr = [(J
\circ \hkr) \otimes \Omega[1] ] \circ \alpha_R \text{ . }$$ But by
\eqref{ga}
$$J \circ \hkr = \hkr \circ S \implies -\Phi_L \circ \hkr \circ S = (\hkr \otimes \Omega[1]) \circ (S \otimes \Omega[1]) \circ \alpha_R \text{ .}$$
Since $\alpha_L = -(S \otimes \Omega[1]) \circ \alpha_R \circ S$ and
$S \circ S = \textbf{1}_{\hc^{\bullet}}(X)$, $$\Phi_L \circ \hkr =
(\hkr \otimes \Omega[1]) \circ \alpha_L  \text{ . }$$

This proves that the diagram \eqref{ee} commutes.\\
\end{proof}

\subsection{A long remark - enlarging the dictionary }
This subsection is a continuation of Section 2.3. Recall that
$\hc^{\bullet}(X)$ should be thought of as analogous to the ring of
functions on a "symmetric" neighborhood $U_G$ of the identity in a
Lie group $G$. $T[-1]$ is analogous to the Lie algebra
$\mathfrak{g}$ of $G$. $\sss(\Omega[1])$ is analogous to the ring of
functions on a neighborhood $\calg V$ of $0$ in $\mathfrak{g}$ that
is diffeomorphic to $U_G$ via the exponential map .
$\hkr:\hc^{\bullet}(X) \rar \sss(\Omega[1])$ is analogous to the
pull-back of functions by the map $\text{exp}:\mathfrak{g} \rar
G$. \\

\subsubsection{The classical picture- I} Let us look at the classical picture for now. Choose a basis $\{X_1,..,X_n\}$ of $\mathfrak{g}$ and a
basis $\{Y_1,...,Y_n\}$ of $\mathfrak{g}^*$ dual to
$\{X_1,...,X_n\}$. Then, $\sum_{i=1}^{i=n} X_i \otimes Y_i$ yields
an element $\mathfrak{g} \otimes \mathfrak{g}^*$. Denote the ring of
functions on $U_G$ by $C(G)$. We identify elements of $\mathfrak{g}$
(resp. elements of $\mathfrak{g}^*$) with left-invariant vector
fields (resp. left-invariant $1$-forms) on $U_G$ whenever necessary.
Given any element of $\mathfrak{g} \otimes \mathfrak{g}^*$ , letting
$\mathfrak{g}$ act on $C(G)$ as a differential operator yields us a
map from $C(G)$ to $C(G) \otimes \mathfrak{g}^*$ satisfying the
Leibniz rule. Therefore, any element of $\mathfrak{g} \otimes
\mathfrak{g}^*$ yields a global connection on $C(G)$. It is easy to
verify that the connection on $C(G)$ yielded by $\sum_{i=1}^{i=n}
X_i \otimes Y_i$ is precisely $d_G$, the De-Rham differential from
$C(G)$ to the sections
of the co-tangent bundle of $U_G$. \\

Denote the ring of functions on $\calg V$ by $C(\mathfrak{g})$. Note
that replacing $C(G)$ by $C(\mathfrak{g})$ in the previous paragraph
enables us to conclude that $\sum_{i=1}^{i=n} X_i \otimes Y_i$
yields the connection $d_{\mathfrak{g}}$ on $C(\mathfrak{g})$ where
$d_{\mathfrak{g}}$ is
precisely the De-Rham differential from $C(\mathfrak{g})$ to the sections of the co-tangent bundle of $\calg V$.\\

Note that elements of $C(G) \otimes \mathfrak{g}^*$ are smooth
functions from $U_G$ to $\mathfrak{g}^*$ i.e, sections of the
co-tangent bundle of $U_G$. Similarly, elements of $C(\mathfrak{g})
\otimes \mathfrak{g}^*$ are smooth functions from $\calg V$ to
$\mathfrak{g}^*$ i.e, sections of the co-tangent bundle of $\calg
V$. \\

Given a smooth function $A$ from $\calg V$ to $\enn(\mathfrak{g})$
and a smooth function $h$ from $\calg V$ to $\mathfrak{g}$ (resp.
$\mathfrak{g}^*$), one can obtain a smooth function $A(h)$ from
$\calg V$ to $\mathfrak{g}$ (resp. $\mathfrak{g}^*$) by setting
$$A(h)(Z) = A(Z)h(Z) $$ for every $Z \in \calg V$. Denote the
smooth function $Z \leadsto \text{ad}(Z)$ from $\calg V$ to
$\enn(\mathfrak{g})$ by $\bar{\text{ad}}$.\\

Consider the connection $\Phi$ on $C(\mathfrak{g})$ such that the following diagram commutes: \\

$$\begin{CD}
C(G) @>d_G>> C(G) \otimes \mathfrak{g}^* \\
@VV{\text{exp}^*}V    @V{\text{exp}^* \otimes \mathfrak{g}^*}VV \\
C(\mathfrak{g}) @>\Phi>> C(\mathfrak{g}) \otimes \mathfrak{g}^* \\
\end{CD} $$

We are interested in comparing $\Phi$ with $d_{\mathfrak{g}}$. This
is done as follows: \\

Let $f$ be any function on $U_G$. Then, $d_G(f)$ is a $1$-form on
$U_G$. The pull-back of the $1$-form $d_G(f)$ via the exponential
map is precisely the $1$-form $d_{\mathfrak{g}}(\text{exp}^*(f))$.
Recall the formula $$d(\text{exp})_Z
=\frac{1-\text{e}^{-\text{ad}(Z)}}{\text{ad}(Z)} $$ $$\implies
(\text{exp}^* \otimes d(\text{exp})^*) =
\frac{1-\text{e}^{-\bar{\text{ad}}}}{\bar{\text{ad}}} \circ
(\text{exp}^* \otimes \mathfrak{g}^*) \text{ . } $$

By definition , $$\Phi(\text{exp}^*(f)) = (\text{exp}^* \otimes
\mathfrak{g}^*) d_G(f) \text{ . } $$ The fact that
$d_{\mathfrak{g}}(\text{exp}^*(f))$ is the pull-back of $d_G(f)$ via
the exponential map implies that $$
d_{\mathfrak{g}}(\text{exp}^*(f)) = (\text{exp}^* \otimes
d(\text{exp})^*)(d_G(f)) =
\frac{1-\text{e}^{-\bar{\text{ad}}}}{\bar{\text{ad}}} \circ
(\text{exp}^* \otimes \mathfrak{g}^*)(d_G(f)) =
\frac{1-\text{e}^{-\bar{\text{ad}}}}{\bar{\text{ad}}} \circ
\Phi(\text{exp}^*(f)) \text{ .} $$

Since any smooth function on $\calg V$ is of the form
$\text{exp}^*(f)$,

\begin{equation} \label{concp1} d_{\mathfrak{g}} =
\frac{1-\text{e}^{-\bar{\text{ad}}}}{\bar{\text{ad}}} \circ \Phi
\end{equation}
$$ \implies \Phi =
\frac{\bar{\text{ad}}}{1-\text{e}^{-\bar{\text{ad}}}} \circ
d_{\mathfrak{g}} \text{ . }$$

Now assume that $U_G$ and $\calg V$ are chosen so that $\calg V$ is
a sufficiently small open disc in $\mathfrak{g}$ containing $0$.
Also assume that $G$ is not $1$-dimensional. Then, any closed
$1$-form on $U_G$ is also exact. Let $Y \in \mathfrak{g}^*$. Since
the left-invariant $1$-form $1 \otimes Y \in C(G) \otimes
\mathfrak{g}^*$ is closed, it is exact as well. Thus, there is a
function $f_Y$ on $U_G$ such that $d_G(f_Y) = 1 \otimes Y$. Then,
$\Phi(\text{exp}^*(f_Y)) = 1 \otimes Y $. On the other hand,
$d_{\mathfrak{g}}(\text{exp}^*(f_Y))$ is the pull-back of the
$1$-form $1 \otimes Y$ via the exponential map.  The formula
$\eqref{concp1}$ thus amounts to the formula for the pullback of a
left-invariant
$1$-form on $U_G$ via the exponential map.\\

Note that the bracket $\text{ad}: \mathfrak{g} \otimes \mathfrak{g}
\rar \mathfrak{g}$ induces a map $\text{ad}:\mathfrak{g}^* \rar
\mathfrak{g}^* \otimes \mathfrak{g}^*$. Let $\mu:C(\mathfrak{g})
\otimes \mathfrak{g}^* \rar C(\mathfrak{g})$ denote the product
taken after treating an element of $\mathfrak{g}^*$ as a function on
$\calg V$. We now claim that as an endomorphism in the space of
smooth sections of $\calg V \times \mathfrak{g}^*$ ,
$\bar{\text{ad}}$ is given by the following composite map
$$\begin{CD} C(\mathfrak{g}) \otimes \mathfrak{g}^* @>C(\mathfrak{g}) \otimes \text{ad} >> C(\mathfrak{g}) \otimes \mathfrak{g}^* \otimes
\mathfrak{g}^* @>\mu \otimes \mathfrak{g}^* >> C(\mathfrak{g})
\otimes \mathfrak{g}^* \end{CD} \text{ . }$$

{\it Denote the above composite map by $\omega_{\mathfrak{g}}$} .\\

To verify this, choose a basis $\{X_i\}$ of $\mathfrak{g}$ and a basis $\{Y_i\}$ of $\mathfrak{g}^*$ dual to $\{X_i\}$. Suppose that $[X_i,X_j]
= \sum_k C_{ij}^k X_k$. Then, $\text{ad}(Y_k) = \sum_{i,j} C_{ij}^k Y_i \otimes Y_j $. Therefore, if $f \in C(\mathfrak{g})$, then,

$$(\mu \otimes \mathfrak{g}^*) \circ (C(\mathfrak{g}) \otimes \text{ad})(f \otimes Y_k)(\sum_i a_i X_i) = f(\sum_i a_i X_i). \sum_{i,j} a_i
C_{ij}^k Y_j \text{ . } $$

On the other hand, $$\bar{\text{ad}}(f \otimes Y_k)(\sum_i a_i X_i)
= f(\sum_i a_i X_i) \text{ad}(\sum_i a_i X_i)(Y_k) $$ But
$\text{ad}(X_i)(Y_k) = \sum_{j} C_{ij}^k Y_j$. Therefore,
$$\bar{\text{ad}}(f \otimes Y_k)(\sum_i a_i X_i)  = f(\sum_i a_iX_i)
\sum_{i,j} a_iC_{ij}^k Y_j \text{ . } $$.

It follows that
$$\frac{\bar{\text{ad}}}{1-\text{e}^{-\bar{\text{ad}}}} =
\frac{\omega_{\mathfrak{g}}}{1-\text{e}^{-\omega_{\mathfrak{g}}}} $$
and that $$\Phi =
\frac{\omega_{\mathfrak{g}}}{1-\text{e}^{-\omega_{\mathfrak{g}}}}
\circ d_{\mathfrak{g}} \text{ . } $$

\subsubsection{The classical picture- II} Let $\bar{\text{exp}}$ denote the map from $\calg V$ to $U_G$
such that $Z \leadsto \text{exp}(Z)^{-1} = \text{exp}(-Z)$. The
discussion in Section 2.5.1 with $\text{exp}$ replaced by
$\bar{\text{exp}}$ together with the fact that
$$d(\bar{\text{exp}})_Z = -\frac{1-\text{e}^{-\text{ad}(-Z)}}{\text{ad}(-Z)} = -\frac{\text{e}^{\text{ad}(Z)}-1}{\text{ad}(Z)} $$
tells us that if $$\Psi =
-\frac{\omega_{\mathfrak{g}}}{\text{e}^{\omega_{\mathfrak{g}}}-1}
\circ d_{\mathfrak{g}}$$ then the following diagram commutes :

$$\begin{CD}
C(G) \otimes \mathfrak{g}^* @>d_G>> C(G) \otimes \mathfrak{g}^* \\
@VV{\bar{\text{exp}}^*}V        @V{\bar{\text{exp}}^* \otimes \mathfrak{g}^*}VV\\
C(\mathfrak{g}) @>\Psi>> C(\mathfrak{g}) \otimes \mathfrak{g}^* \\
\end{CD} $$

This is equivalent to the formula , $$d_{\mathfrak{g}} =
-\frac{\text{e}^{\bar{\text{ad}}}-1}{\bar{\text{ad}}} \circ \Psi
\text{ . }$$ The above formula is equivalent to giving a formula for
the pull back by $\bar{\text{exp}}$
of a left-invariant $1$-form on $U_G$. \\

\subsubsection{Enlarging our dictionary}

The discussion in the previous two subsections helps us understand Theorem 2' better. The discussion in Section 2.5.1 says that if $$\Phi =
\frac{\omega_{\mathfrak{g}}}{1-\text{e}^{-\omega_{\mathfrak{g}}}} \circ d_{\mathfrak{g}}: C(\mathfrak{g}) \otimes \mathfrak{g}^* \rar
C(\mathfrak{g}) \otimes \mathfrak{g}^*$$ then the diagram
$$\begin{CD}
C(G) @>d_G>> C(G) \otimes \mathfrak{g}^* \\
@VV{\text{exp}^*}V    @V{\text{exp}^* \otimes \mathfrak{g}^*}VV \\
C(\mathfrak{g}) @>\Phi>> C(\mathfrak{g}) \otimes \mathfrak{g}^* \\
\end{CD} $$
commutes. The analogy between this and the diagram \eqref{dd} of
Theorem 2' is now fairly explicit. $\hc^{\bullet}(X)$ is analogous
to $C(G)$. The map $\alpha_R$ is analogous to the connection $d_G$.
$\hkr$ is analogous to $\text{exp}^*$. $T[-1]$ is analogous to
$\mathfrak{g}$. $\Omega[1]$ is analogous to $\mathfrak{g}^*$.
$\sss(\Omega[1])$ is analogous to $C(\mathfrak{g})$. The map
$\bar{\textbf{C}}$ is analogous to $d_{\mathfrak{g}}$. The map
$\bar{\omega}:\sss(\Omega[1]) \otimes  \Omega[1] \rar
\sss(\Omega[1]) \otimes \Omega[1]$ is analogous to
$\omega_{\mathfrak{g}}$, and the map $\Phi_R$ is analogous to
$\Phi$. {\it In short, the diagram in \eqref{dd} is analogous to the
computation of the Duflo-like term describing the correction that
needs to be applied to $d_{\mathfrak{g}}$ to yield the map $\Phi$.
This in
turn is equivalent to the formula for the pull-back of a left-invariant $1$-form on $U_G$ via the exponential map. } \\

The map $S:\hc^{\bullet}(X) \rar \hc^{\bullet}(X)$ is analogous to
the map $S_G:C(G) \rar C(G)$ given by pulling back an element of
$C(G)$ by the map $g \leadsto g^{-1}$. The map $-\alpha_L = (S
\otimes \Omega[1]) \circ \alpha_R \circ S$ is therefore analogous to
$(S_G \otimes \mathfrak{g}^*) \circ d_G \circ S_G$. Moreover,
$\text{exp}^* \circ S_G = \bar{\text{exp}}^*$. Therefore $\hkr \circ
S$ is analogous to $\bar{\text{exp}}^*$.It follows that $\Phi_L$ is
analogous to the map $-\Psi$ of the previous subsection. {\it In
short, the diagram \eqref{ee} is analogous to the computation of the
term describing the correction that has to be applied to
$d_{\mathfrak{g}}$ to describe the map $- \Psi$. This in turn is
equivalent to the formula for the pull-back of a left-invariant
$1$-form on $U_G$
via the map $\bar{\text{exp}}$}. \\

Finally we recall from [1] that the Universal enveloping algebra in
$\dcat$ of the Lie algebra object $T[-1]$ was shown to be
represented by the complex $\dd^{\bullet}(X)$ of poly-differential
operators with Hochschild co-boundary. Results equivalent to this
statement were obtained using different methods by Roberts and
Willerton [8] and Markarian [6] as well. Yekutieli [5] showed that
the functor $\text{RD}$ applied to $\dd^{\bullet}(X)$ yields
$\hc^{\bullet}(X)$. It follows that the Hopf-algebra object
$\hc^{\bullet}(X)$ of $\dcat$ is the "dual" in $\dcat$ of
the universal enveloping algebra in $\dcat$ of the Lie algebra object $T[-1]$ of $\dcat$.\\

We warn the reader that the material in the remaining part of this subsection is rather hazy. It can be verified that $\hc^{\bullet}(X)$
satisfies all the formal properties that a ring of functions on a Lie group is required to satisfy. It might therefore be possible to say that
$\hc^{\bullet}(X)$ corresponds to a "Lie group object" in $\dcat$ . \\

One has to be very careful here. Since the concept of a geometric object like a manifold in $\dcat$ does not make sense by itself, the best we
can do is to try to define such a notion by attempting to define a {\it ring of functions on a manifold} in $\dcat$. This has to be a
commutative algebra object in
$\dcat$. If such a definition is possible , $T[-1]$ thought of as a manifold in $\dcat$ should correspond to the algebra object $\sss(\Omega[1])$ of $\dcat$.\\

 The Lie algebra of the  Lie group object $\hc^{\bullet}(X)$ should be $T[-1]$ . The diagrams
\eqref{dd} and \eqref{ee} in Theorem 2' could then be thought of as
being equivalent to "computing the pull-back of a $1$-form on on the
Lie group object $\hc^{\bullet}(X)$ via the maps $\text{exp}$ and
$\bar{\text{exp}}$ respectively". Of course, $\text{exp}$ and
$\bar{\text{exp}}$ are defined solely by what they are as maps from
$\hc^{\bullet}(X)$ to $\sss(\Omega[1])$. However, for this to make any sense, one should
be able to define the notion of a differential form on a manifold in $\dcat$. \\

\section{Some essential linear algebra}
The first subsection of this section proves some propositions and a
lemma in linear algebra. The second subsection describes their
extensions pertaining to the object $\sss(\Omega[1])$ of $\dcat$. We
remind the reader that maps between tensor products of graded
$\fld$-vector spaces (or graded $\strc$-modules) that rearrange
factors are assumed to take the appropriate signs into account. If
$x$ is a homogenous
element of a graded $\fld$-vector space $W$, $|x|$ will denote the degree of $x$. \\

\subsection{Some propositions and a lemma of linear algebra}
In this subsection, we will work with differential graded vector spaces over a field $\fld$ of characteristic $0$. Almost every dg-vector space
in this section has $0$ differential. We shall assume that the differential on a graded vector space is $0$ unless we explicitly say otherwise.
Let $V$ be a finite dimensional
vector space over $\fld$. Denote the dual of $V$ by $V^*$. {\it Denote the dimension of $V$ in this subsection by $m$}.\\

As usual $\sss(V[1])$ denotes the symmetric algebra $\oplus_i \wedge^i V[i]$ generated over $\fld$ by $V$ concentrated in degree $-1$.
Similarly, $\sss(V^*[-1])$ denotes the symmetric algebra $\oplus_i \wedge^iV^*[-i]$ generated over $\fld$ by $V^*$ concentrated in degree $1$.
Note that the products on $\sss(V[1])$ and $\sss(V^*[-1])$
are wedge products.  \\

There is a map $$\textbf{i}_V: \sss(V[1]) \rar \enn(\sss(V[1])) $$
$$ Z \leadsto (W \leadsto W \wedge Z) \text{ . } $$
$\textbf{i}_V$ takes an element $Z$ of $\sss(V[1])$ to the
endomorphism of $\sss(V[1])$ given by multiplication by $Z$ on the
right. \\

We will often denote the related map $$\sss(V[1]) \otimes \sss(V[1]) \rar \sss(V[1]) $$ $$ W \otimes Z \leadsto W \wedge Z $$ by $(\text{ - }
\wedge \text{ - })_V$. The subscript $V$ may be dropped at times when it is obvious. \\

Choose a basis $\{x_1,....,x_m\}$ of $V$. Let $\{y_1,...,y_m\}$ be a
basis of $V^*$ dual to $\{x_1,....,x_m\}$. Let
$\textbf{j}_{V^*}(y_i)$ be the endomorphism of $\sss(V[1])$ given by
$$ \textbf{j}_{V^*}(y_i)(x_j) = \delta_{ij}$$ $$ \textbf{j}_{V^*}(y_i)(1) = 0  $$
\begin{equation} \label{leib} \textbf{j}_{V^*}(y_i)(x_{i_1} \wedge .... \wedge x_{i_k}) =
\delta_{ii_k} x_{i_1} \wedge .... \wedge x_{i_{k-1}} -
\textbf{j}_{V^*}(y_i)(x_{i_1} \wedge .... \wedge x_{i_{k-1}}) \wedge
x_{i_k} \text{ . } \end{equation}\\

Note that $\textbf{j}_{V^*}$ extends by linearity to a map from $V^*[-1]$ to $\enn(\sss(V[1]))$. Let \\
$\circ : \enn(\sss(V[1]))^{\otimes 2}
\rar \enn(\sss(V[1]))$ denote the composition map. Extend $\textbf{j}_{V^*}$ to a map
$$\textbf{j}_{V^*}:\sss(V^*[-1]) \rar \enn(\sss(V[1]))$$ by setting
\begin{equation} \label{hodop} \textbf{j}_{V^*}(Y_1 \wedge Y_2) = \textbf{j}_{V^*}(Y_2) \circ
\textbf{j}_{V^*}(Y_1) \text{   } \forall \text{  } Y_1,Y_2 \in
\sss(V^*[-1]) \text{ . } \end{equation}

We will often denote the related map $$\sss(V[1]) \otimes \sss(V^*[-1]) \rar \sss(V[1]) $$ $$ W \otimes Y \leadsto \textbf{j}_{V^*}(Y)(W) $$ by
$(\text{- }
\bullet \text{ - })_V$. The subscript $V$ may be dropped at times when it is obvious. \\

\subsection*{Remark 1 - A geometric analogy:} One can think of $\sss(V[1])$ as the ring of
functions on an odd super-manifold $M_V$. For an element $Z$ of
$\sss(V[1])$, $\textbf{i}_V(Z)$ is just the operator given by
"multiplication on the right by $Z$". All operators on $\sss(V[1])$
act on the right in our viewpoint. From this viewpoint, elements of
$V^*[-1]$ yield {\it vector fields} on $M_V$. The map
$\textbf{j}_{V^*}$ takes an element $Y$ of $V^*[-1]$ to the constant
vector field on $M_V$ associated with $Y$. The reader can observe
that the {\it Leibniz rule is part of the definition} of the map
$\textbf{j}_{V^*}:V^*[-1] \rar \enn(\sss(V[1]))$ ( see equation
\eqref{leib}). The equation \eqref{hodop} in the definition of
$\textbf{j}_{V^*}$
 just says that the map $\textbf{j}_{V^*}$ takes an
element $Y$ of $\sss(V^*[-1])$ to the constant differential operator
on $M_V$ associated with $Y$. \\

\noindent Note that on $\enn(\sss(V[1]))$ we have a natural product given by the composition $\circ$ . Let $\circ^{\text{op}}$ denote the
product on $\enn(\sss(V[1]))^{\text{op}}$. If $a,b \in \enn(\sss(V[1]))$ then $a \circ^{\text{op}} b = b \circ a$.
We have the following proposition \\

\begin{prop}
The composite
$$\begin{CD} \sss(V^*[-1]) \otimes \sss(V[1]) @>\textbf{j}_{V^*} \otimes
\textbf{i}_{V} >> \enn(\sss(V[1])) \otimes \enn(\sss(V[1]) @> \circ^{\text{op}} >> \enn(\sss(V[1])) \end{CD} $$ is an isomorphism of graded
$\fld$ vector spaces with $0$ differential.
\end{prop}

\begin{proof}
Denote the given composite by $G_r$. \\

Since all vector spaces involved in this proposition have $0$
differential, it suffices to check that $G_r$ is an isomorphism of
graded $\fld$-vector spaces. Further, it is easy to check that
$\textbf{i}_V$, $\textbf{j}_{V^*}$ and $\circ^{\text{op}}$ are
degree preserving. It therefore, suffices to check that $G_r$ is an
isomorphism of $\fld$-vector spaces. All the $\fld$-vector spaces
involved in this proposition are finite-dimensional. Further, $
\sss(V^*[-1]) \otimes \sss(V[1])$ and $ \enn(\sss(V[1]))$ have the
same dimension as $\fld$-vector spaces. It
therefore, suffices to check that $G_r$ injective. \\

Choose a basis $\{x_1,....,x_m\}$ of $V$. Let $\{y_1,...,y_m\}$ be a
basis of $V^*$ dual to $\{x_1,....,x_m\}$. By convention, we list
the elements of any subset of $\{1,...,m\}$ in {\it ascending
order}. We can then define an ordering $\prec$ on the set of subsets
of $\{1,...,m\}$ by setting
 $$ S \prec T \text{ if } |S| < |T| $$
$$ \{i_1,...,i_k\} \prec \{j_1,...,j_k\} \text{ if } (i_1,..,i_k)
\prec (j_1,...,j_k) \text{ in the lexicographic order . } $$

For $S=\{i_1,...,i_k\}$, let $x_S$ denote $x_{i_1} \wedge ... \wedge x_{i_k} $. Similarly , $y_S$ will denote $y_{i_1} \wedge .... \wedge
y_{i_k} $. It is easy to verify that
\begin{equation} \label{prel} \textbf{j}_{V^*}(y_S)(x_S) = \pm 1
\end{equation}
$$\textbf{j}_{V^*}(y_S)(x_T) = 0 \text{ if } T \prec S \text{ . }$$

An element of $\sss(V^*[-1]) \otimes \sss(V[1])$ is given by an
expression of the form

$$ \sum_{S \subset \{1,...,m\}} y_S \otimes a_S \text{ , } a_S \in
\sss(V[1]) \text{ . } $$

Choose $S_0 \subset \{1,...,m\}$ to be the {\it least} subset (under
the ordering $\prec$ ) of $\{1,...,m\}$ such that $a_{S_0} \neq 0$.
Then,

$$G_r(\sum_{S \subset \{1,...,m\}} y_S \otimes a_S)(x_{S_0}) = \pm a_{S_0} $$

by \eqref{prel}. It follows that $G_r$ is injective. This proves the
desired proposition. \\

\end{proof}

\subsection*{Remark - 2 } Let $M_V$ be the super-manifold of whose
ring of functions is $\sss(V[1])$. Under our convention that all
operators on $\sss(V[1])$ act on the right, $G_r$ identifies
$\wedge^iT[-i] \otimes \sss(V[1])   $ with the space of principal
symbols of differential operators of order $i$ on $M_V$. Proposition
10 says that every endomorphism of $\sss(V[1])$ is given by a
differential operator on
$M_V$. \\

\subsection*{Notation} For the rest of this paper, $G_r$ shall denote
the isomorphism in Proposition 10 and $F_r$ shall denote its inverse. In addition, if \\
$\tau: \sss(V^*[-1]) \otimes \sss(V[1]) \rar \sss(V[1])
\otimes \sss(V^*[-1])$ denotes the swap map the composite
$$
\begin{CD}\sss(V^*[-1]) \otimes \sss(V[1]) @> \tau>> \sss(V[1])
\otimes \sss(V^*[-1]) @>\textbf{i}_V \circ^{\text{op}} \textbf{j}_{V^*}>> \enn(\sss(V[1])) \end{CD} $$ shall be denoted by $G_l$ .
 The inverse of $G_l$ will be denoted by $F_l$. \\

Let $\pi_k:\sss(V[1]) \rar \wedge^kV[k]$ denote the natural
projection. Note that we have a non-degenerate pairing
$\langle,\rangle$ on
 $\sss(V[1])$. This is given by the composite
 $$\begin{CD} \sss(V[1]) \otimes \sss(V[1]) @>\wedge>> \sss(V[1])
 @>\pi_m>> \wedge^mV[m] \\ \end{CD} $$ where $m$ is the dimension of
 $V$.

\subsection*{Definition :} The {\it adjoint }  $\Phi^+ \in \enn(\sss(V[1]))$ of a homogenous
element $\Phi$ of $\enn(\sss(V[1]))$ is the unique element of
$\enn(\sss(V[1]))$ satisfying
$$ \langle \Phi(a),b \rangle = {(-1)}^{|\Phi||b|}\langle a,\Phi^+(b)
\rangle \text{  } \forall \text{ homogenous } a,b \in \sss(V[1])
\text{ . }
$$ The adjoint of an arbitrary element of $\enn(\sss(V[1])$ is the
sum
of the adjoints of its homogenous components.\\

Let $$\text{ev}_V: \sss(V[1]) \otimes \enn(\sss(V[1])) $$ be the map
$$ W \otimes \Phi \leadsto \Phi(W) \text{ . }$$ Denote the map which takes an
element of $\enn(\sss(V[1]))$ to its adjoint by $\textbf{A}_V$. Let $\langle \text{ },\text{ev}^+ \rangle_V$ denote the composite
$$\begin{CD}
\sss(V[1]) \otimes \sss(V[1]) \otimes \enn(\sss(V[1])\\
@VV{\sss(V[1]) \otimes \sss(V[1]) \otimes \textbf{A}_V}V\\
\sss(V[1]) \otimes \sss(V[1]) \otimes \enn(\sss(V[1]) @>\sss(V[1]) \otimes \text{ev}_V>> \sss(V[1]) \otimes \sss(V[1])
 @>\langle \text{ },\text{ } \rangle>> \wedge^mV[m] \end{CD} $$.

  Then, the following diagram commutes: \\

$$\begin{CD}
\sss(V[1]) \otimes \enn(\sss(V[1])) \otimes \sss(V[1]) @> \sss(V[1]) \otimes \tau >> \sss(V[1]) \otimes \sss(V[1]) \otimes \enn(\sss(V[1])) \\
@VV{\text{ev}_V \otimes \sss(V[1])}V     @V{\langle \text{ },\text{ev}^+ \rangle_V}VV\\
\sss(V[1]) \otimes \sss(V[1])    @> \langle,\rangle>> \wedge^mV[m] \\
\end{CD} $$

The map $\tau: \enn(\sss(V[1])) \otimes \sss(V[1]) \rar \sss(V[1]) \otimes \enn(\sss(V[1]))$ in the top row of the above diagram swaps factors.\\

The following proposition
describes the basic properties of the adjoint.  \\

\begin{prop}
1. Let $L$ be an element of $\enn(\sss(V[1]))$. Then,
\begin{equation*} L^{++} = L \text{ . }
\end{equation*}
2. If $L_1$ and $L_2$ are homogenous elements of $\enn(\sss(V[1]))$,
then $$ (L_1 \circ L_2)^+ = {(-1)}^{|L_1||L_2|} (L_2)^+ \circ
(L_1)^+ \text{ . }$$
3. If $Z \in \sss(V[1])$ is a homogenous
element, then
 $$\textbf{i}_V(Z)^+ = \textbf{i}_V(Z) \text{ . } $$
4. If $Y \in \sss(V^*[-1])$ is a homogenous element, then
$$\textbf{j}_{V^*}(Y)^+ = {(-1)}^{|Y|} \textbf{j}_{V^*}(Y) \text{ . }$$
5. If $Z \in \sss(V[1])$ and $Y \in \sss(V^*[-1])$ are homogenous
elements, then
$$G_l( Y \otimes Z)^+ = {(-1)}^{|Y|}
G_r(Y \otimes Z) \text{ . }$$
\end{prop}

\begin{proof}
\noindent  Observe that if $L \in \enn(\sss(V[1]))$ is homogenous, then $|L^+| = |L|$. Also note that if $a,b \in \sss(V[1])$ are homogenous,
then $\langle a,b \rangle = {(-1)}^{|a||b|} \langle b,a \rangle$. Also recall that the
pairing $\langle,\rangle$ is non-degenerate.\\

\noindent If $a,b \in \sss(V[1])$ are homogenous elements and if $L
\in \enn(\sss(V[1]))$ is homogenous  , note that
$$ \langle L^{++}(a),b \rangle = {(-1)}^{|b||a|+|b||L|} \langle
b,L^{++}(a) \rangle = {(-1)}^{|b||a| +|b||L|+|a||L|} \langle
L^+(b),a \rangle $$ $$ =  {(-1)}^{|b||a|
+|b||L|+|a||L|+|a||b|+|a||L|} \langle a,L^+(b) \rangle = \langle
L(a),b \rangle \text{ . }$$ Part 1 of this proposition now follows
immediately
from this calculation. \\

\noindent For the rest of this proof $a$ and $b$ shall be homogenous
elements
of $\sss(V[1])$. \\

\noindent If $L_1$ and $L_2$ are homogenous elements of
$\enn(\sss(V[1]))$ then, $$ {(-1)}^{|L_1||L_2|} \langle a, (L_2)^+
\circ (L_1)^+ (b) \rangle = {(-1)}^{|L_1||L_2|+ |L_2|(|b|+|L_1|)}
\langle L_2(a), L_1^+(b) \rangle $$ $$ = {(-1)}^{|L_1||L_2|+
|L_2|(|b|+|L_1|)+ |L_1||b|} \langle L_1 \circ L_2(a),b \rangle =
{(-1)}^{|b|(|L_1|+|L_2|)} \langle L_1 \circ L_2(a),b \rangle \text{
. }$$ Part 2 of this proposition now follows from the observation
that $|L_1
\circ L_2| = |L_1| +|L_2|$. \\

\noindent Part 3 of this proposition is immediate from the relevant
definitions and the fact that
$$a \wedge Z \wedge b = {(-1)}^{|b||Z|} a \wedge b \wedge Z \text{ . }$$

\noindent To verify part 4, choose a basis $\{x_1,...,x_m\}$ of $V$. Let $\{y_1,...,y_m\}$ be a basis of $V^*$ dual to $\{x_1,..,x_m\}$. For an
ordered subset $S = \{i_1,..,i_k\}$ of $\{1,...,m\}$ let $y_S$ denote $y_{i_1} \wedge .... \wedge y_{i_k}$ and let $x_S$ denote $x_{i_1} \wedge
... \wedge x_{i_k}$. If $T$ is disjoint from $S$, note that
$$ \textbf{j}_{V^*}(y_S)(x_T \wedge x_S) = x_T \wedge
\textbf{j}_{V^*}(y_S)(x_S) = {(-1)}^{\frac{k(k-1)}{2}} x_T $$. Let
$T$ and $T'$ be subsets of $\{1,..,m\}$ disjoint from $S$. Then,
$$ \textbf{j}_{V^*}(y_S)(x_T \wedge x_S) \wedge (x_{T'} \wedge x_S)
= {(-1)}^{\frac{k(k-1)}{2}} x_T \wedge x_{T'} \wedge x_S =
{(-1)}^{\frac{k(k-1)}{2}} {(-1)}^{|T'||S|} x_T \wedge x_{S} \wedge
x_{T'} $$ $$ = {(-1)}^{|T'||S|} x_T \wedge x_{S} \wedge
\textbf{j}_{V^*}(y_S)(x_{T'} \wedge x_S) =
{(-1)}^{(|T'|+|S|)|S|+|S|} x_T \wedge x_{S} \wedge
\textbf{j}_{V^*}(y_S)(x_{T'} \wedge x_S) \text{ . }$$ Putting $a =
x_T \wedge x_S$, $b= x_{T'} \wedge x_{S}$ we see that $|b| =
|T'|+|S|$. Further, $|Y|= -|S|$. Part 4 now follows from the above
computation once we
recall that $\langle a,b \rangle = \pi_m(a \wedge b) $. \\

\noindent Part 5 follows from part 2,part 3 and part 4. $G_l(Y
\otimes Z) = {(-1)}^{|Z||Y|}  \textbf{j}_{V^*}(Y)  \circ
\textbf{i}_{V}(Z) $. Thus,
$$G_l(Y \otimes Z)^+ = {(-1)}^{|Z||Y|} (\textbf{j}_{V^*}(Y)  \circ
\textbf{i}_{V}(Z))^+ ={(-1)}^{|Z||Y|} {(-1)}^{|Z||Y|}
(\textbf{i}_{V}(Z))^+ \circ  (\textbf{j}_{V^*}(Y))^+$$ $$ =
{(-1)}^{|Y|}(\textbf{i}_{V}(Z))  \circ (\textbf{j}_{V^*}(Y)) =
{(-1)}^{|Y|} G_r(Y \otimes Z) \text{ . }$$

\end{proof}

\noindent Recall that $\pi_j:\sss(V[1]) \rar \wedge^jV[j]$ denotes
the natural projection. We will denote the projection
$$ \sss(V^*[-1]) \otimes \pi_j: \sss(V^*[-1])  \otimes  \sss(V[1]) \rar
 \sss(V^*[-1]) \otimes \wedge^jV[j] $$ by $\pi_j$ itself. Let $I$ denote the
endomorphism of $\sss(V^*[-1])$ that multiplies $\wedge^iV^*[-i]$ by
${(-1)}^i$. The
following Proposition is really a corollary of Proposition 11. \\

\begin{prop}
If $L \in \enn(\sss(V[1]))$, then
$$ \pi_0 (F_l(L)) = I (\pi_0 ( F_r(L^+))) \text{ . }$$
\end{prop}

\begin{proof}
This is almost immediate from Part 5 of Proposition 11. Let $Y \in \sss(V^*[-1])$ and $Z \in \sss(V[1])$ be homogenous. By part 5 of Proposition
11
$$G_l(Y \otimes Z)^+ = {(-1)}^{|Y|}G_r(Y \otimes Z) \text{ . }
$$ By definition, $F_r(G_r(Y \otimes Z)) = Y \otimes Z$. On the other hand,
$$F_l(G_r(Y \otimes Z)^+) = F_l({(-1)}^{|Y|} G_l(Y \otimes Z) =
{(-1)}^{|Y|} Y \otimes Z \text{ . } $$ Since $G_r$ and
$\textbf{A}_V$ are degree preserving isomorphisms of $\fld$-vector
spaces, any homogenous element in $\enn(\sss(V[1]))$ is of the form
$G_r(Y \otimes Z)^+$. It follows from the above computation that
$$F_l(M)= (I \otimes \sss(V[1]))F_r(M^+)$$ for any homogenous $M \in
\enn(\sss(V[1]))$. Therefore,
$$\pi_j(F_l(L)) = I(\pi_j(F_r(L^+)))
$$ for any $j$.
 When $j=0$, we get
the
desired proposition.\\

\end{proof}

\noindent \textbf{Convention} To simplify notation, we follow the following convention : If \\ $a \in \sss(V^*[-1]) \otimes \sss(V[1])$ , then
$G_r(a)$ will be denoted by $a$ itself. {\it Keeping this convention in mind },

\begin{prop} If $a,b \in \sss(V^*[-1]) \otimes \sss(V[1])$, then,
$$\pi_0(F_r(a \circ b)) = \pi_0(F_r(\pi_0(a) \circ
b)) \text{ . }$$
\end{prop}

\begin{proof}  It suffices to check that if $Z \in \wedge^iV[i]$with
$i >0$, then $$\pi_0(Y \otimes Z \circ Y' \otimes Z') =0 $$ for any $Z' \in \sss(V[1])$, $Y,Y' \in \sss(V^*[-1])$. Note that by our convention,
$$Y \otimes Z \circ Y' \otimes Z' =\textbf{i}_V(Z) \circ \textbf{j}_{V^*}(Y) \circ \textbf{i}_V(Z') \circ \textbf{j}_{V^*}(Y') \text{ . } $$
 For subsets $S$ and $T$ of
$\{1,...,m\}$, let $x_S$ and $y_T$ be as in the proof of Proposition 11, part 4. By Proposition 10,
$$ \textbf{j}_{V^*}(Y) \circ \textbf{i}_V(Z') \circ
\textbf{j}_{V^*}(Y') = \sum_{S,T \subset \{1,...,m\} } a_{S,T} y_T \otimes x_S  $$ for some $a_{S,T} \in \fld$. Then,
$$ \textbf{i}_V(Z) \circ
\textbf{j}_{V^*}(Y) \circ \textbf{i}_V(Z') \circ
\textbf{j}_{V^*}(Y') = \sum_{S,T \subset \{1,...,m\} } a_{S,T} y_T
 \otimes x_S \wedge Z \text{ . } $$ Since $Z \in \wedge^iV[i]$, $x_S \wedge Z
\in \oplus_{k \geq i} \wedge^kV[k] \subset \oplus_{k >0}
\wedge^kV[k]$. It follows that
$$ \pi_0(\sum_{S,T \subset \{1,...,m\} } a_{S,T} y_T \otimes
x_S \wedge Z ) =0 $$. This proves the desired proposition.

\end{proof}

\noindent \textbf{Remark - 3 :} Proposition 10 said that every
element of $\enn(\sss(V[1]))$ can be thought of as a differential
operator on $M_V$. The isomorphism $F_r$ makes this identification.
The map $\pi_0: \sss(V^*[-1]) \otimes \sss(V[1]) \rar \sss(V^*[-1])$
should be thought of as the map which "takes the constant term" of a
differential operator. Proposition 13 says that
$$ \text{const. term}({\mathcal D}_1 \circ {\calg D}_2) = \text{const. term}((\text{ const.
term }({\calg D}_1)) \circ {\calg D}_2) $$ for any two differential
operators ${\calg D}_1$ and
${\calg D}_2$ on $M_V$. \\

\noindent Recall that by Proposition 10, $F_r$ identifies
$\enn(\sss(V[1]))$ with $\sss(V^*[-1]) \otimes \sss(V[1])$. We also
remarked (in Remark 2) that the direct summand $\wedge^iV^*[-i]
 \otimes \sss(V[1])$ of $\sss(V^*[-1]) \otimes \sss(V[1])$ can be
thought of as the space of principal symbols of differential
operators of order $i$ on $\sss(V[1])$. Reflecting this
understanding, we denote
$\wedge^iV^*[-i] \otimes \sss(V[1])$  by $D_i$. \\

Note that the composition $\circ: \enn(\sss(V[1]))^{\otimes 2} \rar \enn(\sss(V[1]))$ equips $\enn(\sss(V[1]))$ with the structure of a graded
associative $\fld$-algebra. This algebra structure induces the structure of a Lie super-algebra on $\enn(\sss(V[1]))$. If $a,b \in
\enn(\sss(V[1]))$ are homogenous, then $$[a,b]_V = a \circ b - {(-1)}^{|a||b|} b \circ a \text{ . }$$ \\

Also note that the map $\circ^{\text{op}}: \enn(\sss(V[1]))^{\otimes 2} \rar \enn(\sss(V[1]))$ equips $\enn(\sss(V[1]))$ with the structure of a
graded associative $\fld$-algebra. We denote this $\fld$-algebra by $\enn(\sss(V[1]))^{\text{op}}$. The algebra structure of
$\enn(\sss(V[1]))^{\text{op}}$ induces the structure of a Lie super-algebra on $\enn(\sss(V[1]))^{\text{op}}$. If $a,b \in \enn(\sss(V[1]))$ are
homogenous, then
$$[a,b]^{\text{op}}_V = a \circ^{\text{op}} b - {(-1)}^{|a||b|} b \circ^{\text{op}} a = [b,a]_V
\text{ . }$$

\begin{prop} $$[D_1,D_1]_V \subset D_1 \text{ . }$$
\end{prop}

\begin{proof}
Let $H \in \sss(V[1])$ ,let $Z \in \sss(V[1])$ and let $Y
\in\sss(V^*[-1])$. Recall that $(H \bullet Y)$ denotes
$\textbf{j}_{V^*}(Y)(H)$. Note that $HZ:= H \wedge Z =
\textbf{i}_V(Z)(H)$. Keep in mind that $Y \otimes Z$ is identified
with $G_r(Y \otimes Z)$. Then, $$Y \otimes Z(H) = (H \bullet Y)Z
\text{ . }$$

Let $Z_1,Z_2 \in \sss(V[1])$ be homogenous and let $y_1,y_2 \in V^*[-1]$. Then, if $H \in \sss(V[1])$,
$$(y_1 \otimes Z_1) \circ (y_2 \otimes Z_2)(H) = ((H \bullet y_2)Z_2 \bullet y_1 )Z_1
= (H \bullet y_2)(Z_2 \bullet y_1)Z_1 + {(-1)}^{|Z_2|} ((H \bullet
y_2) \bullet y_1) Z_2Z_1 $$.

Therefore, \begin{equation} \label{d1d2} (y_1 \otimes Z_1) \circ
(y_2 \otimes Z_2) = y_2 \otimes (Z_2 \bullet y_1)Z_1 +{(-1)}^{|Z_2|}
y_2 \wedge y_1 \otimes Z_2Z_1 \text{ . }
\end{equation}

Similarly, \begin{equation} \label{d2d1} (y_2 \otimes Z_2) \circ
(y_1 \otimes Z_1) = y_1 \otimes (Z_1 \bullet y_2)Z_2 +{(-1)}^{|Z_1|}
y_1 \wedge y_2 \otimes Z_1Z_2 \text{ . }
\end{equation}

If ${\calg D}_1= y_1 \otimes Z_1$ and ${\calg D}_2 = y_2 \otimes Z_2$ then $|{\calg D}_1| = |Z_1|-1$ and $|{\calg D}_2|=|Z_2|-1$. It then
follows from \eqref{d1d2} and \eqref{d2d1} that

\begin{equation} \label{comm} {\calg D}_1 \circ {\calg D}_2 - {(-1)}^{|{\calg D}_1||{\calg D}_2|} {\calg D}_2 \circ {\calg D}_1 = y_2 \otimes (Z_2 \bullet y_1)Z_1  -
{(-1)}^{|{\calg D}_1||{\calg D}_2|} y_1 \otimes (Z_1 \bullet y_2)Z_2
\text{ . }
\end{equation}

Note that the right hand side is an element of $D_1$. The desired proposition now follows immediately.

 \end{proof}

We will assume that the Lie super-algebra $\enn(\sss(V[1]))$ is equipped with the bracket $[,]_V$ unless we explicitly state otherwise. \\

 \noindent Proposition 14 tells us that $D_1$ is a Lie sub-algebra of
 $\enn(\sss(V[1]))$. It also follows immediately from Proposition 14 that $$[D_1,D_1]^{\text{op}}_V \subset D_1$$. It follows that $D_1$ equipped with
 the bracket $[,]_V^{\text{op}}$ is a Lie subalgebra of $\enn(\sss(V[1]))^{\text{op}}$ .
 $D_m$ can be identified with the top symmetric
 power of $D_1$ over $\sss(V[1])$. In other words,
 $$D_m \simeq \textbf{S}_{\sss(V[1])}^m D_1 \text{ . }$$
 One can then study the right adjoint action of the Lie algebra $D_1$ ( equipped with the bracket $[,]_V^{\text{op}}$ ) on
 $D_m$ . For an element $L$ of $D_1$, let $\text{ad}(L)$ denote the
  right adjoint action of $L$ on $D_m$ with respect to the bracket $[,]_V^{\text{op}}$. Then
   $$ \text{ad}(L)({\calg D}_m{\calg D}_{m-1}...{\calg D}_1) =  {\calg D}_m{\calg D}_{m-1}...[{\calg D}_1,L]_V^{\text{op}} +
  {(-1)}^{|L||{\calg D}_1|}{\calg D}_m{\calg D}_{m-1}....[{\calg D}_2,L]_V^{\text{op}}{\calg D}_1 + $$ $$ ... {(-1)}^{|L|(|{\calg D}_1|+... +|{\calg D}_{m-1}|)}
  [{\calg D}_m,L]_V^{\text{op}}....{\calg D}_1 $$ $$ =  {\calg D}_m{\calg D}_{m-1}...[L,{\calg D}_1] +
  {(-1)}^{|L||{\calg D}_1|}{\calg D}_m{\calg D}_{m-1}....[L,{\calg D}_2]{\calg D}_1 + $$ $$ ... {(-1)}^{|L|(|{\calg D}_1|+... +|{\calg D}_{m-1}|)}
  [L,{\calg D}_m]....{\calg D}_1 $$ for homogenous elements ${\calg D}_i \in D_1$. In defining $\text{ad}(L)$,
  $D_m$ is treated as $\textbf{S}_{\sss(V[1])}^m D_1$. \\

One also has the right adjoint action $\bar{\text{ad}}(L)$ of $L$ on $D_1^{\otimes m}$ (with respect to the bracket $[,]_V^{\text{op}}$ on
$D_1$).
$$\bar{\text{ad}}(L)( {\calg D}_m \otimes {\calg D}_{m-1}\otimes ...\otimes {\calg D}_1) = {\calg D}_m \otimes {\calg D}_{m-1} \otimes
...\otimes [{\calg D}_1,L]_V^{\text{op}}  $$ $$+
  {(-1)}^{|L||{\calg D}_1|}{\calg D}_m \otimes {\calg D}_{m-1} \otimes .... \otimes [{\calg D}_2,L]_V^{\text{op}} \otimes {\calg D}_1 + $$ $$ ... {(-1)}^{|L|(|{\calg D}_1|+... +|{\calg D}_{m-1}|)}
  [{\calg D}_m,L]_V^{\text{op}} \otimes .... \otimes {\calg D}_1 $$ $$ = {\calg D}_m \otimes {\calg D}_{m-1} \otimes
...\otimes [L,{\calg D}_1] +
  {(-1)}^{|L||{\calg D}_1|}{\calg D}_m \otimes {\calg D}_{m-1} \otimes .... \otimes [L,{\calg D}_2] \otimes {\calg D}_1 + $$ $$ ... {(-1)}^{|L|(|{\calg D}_1|+... +|{\calg D}_{m-1}|)}
  [L,{\calg D}_m] \otimes .... \otimes {\calg D}_1 $$ for homogenous elements ${\calg D}_i \in D_1$.

  Let $p: D_1^{\otimes m} \rar D_m$ denote the map $${\calg D}_m \otimes {\calg D}_{m-1}\otimes
...\otimes {\calg D}_1 \leadsto {\calg D}_m{\calg D}_{m-1}...{\calg D}_1$$. The following proposition is clear from the definitions. \\

\begin{prop}
The following diagram commutes :\\

$$\begin{CD}
D_1^{\otimes m} @>p>> D_m \\
@VV{\bar{\text{ad}}(L)}V   @V{\text{ad}(L)}VV \\
D_1^{\otimes m} @>p>> D_m \\
\end{CD} $$
\end{prop}

 \noindent Choose a basis $\{x_1,...,x_m\}$ of $V$ and a basis $\{y_1,...,y_m\}$ of $V^*$ dual to $\{x_1,....,x_m\}$.
  Let $\textbf{1}^m: \fld \rar \wedge^mV^*[-m] \otimes
 \wedge^mV[m]$ be the map $$1 \leadsto y_m \wedge ..... \wedge y_1 \otimes x_1 \wedge ... \wedge x_m  \text{ . }$$

\noindent Let $\tau: \sss(V[1]) \otimes \wedge^mV^*[-m] \rar \wedge^m V^*[-m]
 \otimes \sss(V[1])$ denote the swap map.
 Denote the composite map $$\begin{CD} \sss(V[1]) @>(\tau \otimes \wedge^mV[m]) \circ (\sss(V[1])
 \otimes \textbf{1}^m ) >>    \wedge^m V^*[-m]
 \otimes \sss(V[1]) \otimes \wedge^mV[m] \end{CD} $$ by $\textbf{1}^m$
 itself.

 \begin{lem}
 Let $L \in D_1$ be a homogenous element. The following diagram commutes : \\
 $$\begin{CD}
 \sss(V[1]) @>\textbf{1}^m>> \wedge^m V^*[-m] \otimes \sss(V[1])
 \otimes \wedge^m V[m] \\
 @VV{-L^+}V     @V{{(-1)}^{|L|m} \text{ad}(L) \otimes \wedge^m V[m]}VV \\
  \sss(V[1]) @>\textbf{1}^m>> \wedge^m V^*[-m] \otimes \sss(V[1])
 \otimes \wedge^m V[m] \\
 \end{CD} $$
 \end{lem}

\begin{proof}
This lemma is proven by a direct computation. We once more recall that if $H,Z \in \sss(V[1])$ and if $Y \in
\sss(V^*[-1])$ then $HZ$ denotes $H \wedge Z = \textbf{i}_{V}(Z)(H)$ and $(H \bullet Y)$ denotes $\textbf{j}_{V^*}(Y)(H)$. \\

Choose a basis $\{y_1,...,y_m\}$ of $V^*$ and a basis $\{x_1,..,x_m\}$ of $V$ dual to $\{y_1,...,y_m\}$ . We may assume without loss of
generality that $L= y_1 \otimes Z$ where $Z \in \sss(V[1])$ is homogenous. Then, if $H \in \sss(V[1])$ is homogenous , by Proposition 11, Parts
2,3, and 4,

$$L^+(H) = {(-1)}^{-|Z|-1} ((HZ) \bullet y_1) = {(-1)}^{|Z|+1} ((HZ) \bullet y_1) \text{ . }$$ Thus,
\begin{equation} \label{lhs} \textbf{1}^m(L^+(H)) = {(-1)}^{(|H|+|Z|-1)m}{(-1)}^{|Z|+1} y_m \wedge ... \wedge y_1 \otimes  ((HZ) \bullet
y_1) \otimes x_1 \wedge ... \wedge x_m \text{ . } \end{equation}

On the other hand, if we treat $H$ as an element of $D_0$, then,

$$ (y_1 \otimes Z) \circ H (P) = (PH \bullet y_1)Z = P (H  \bullet y_1) Z +{(-1)}^{|H|}(P \bullet y_1)HZ $$
$$ H \circ (y_1 \otimes Z) (P) = (P \bullet y_1) ZH \text{ . } $$

Note that the degree of the operator $y_1 \otimes Z$ is $|Z|-1$. It follows that

$$[ y_1 \otimes Z,H](P) = P(H \bullet y_1)Z \text{ . }$$ Thus, $$[y_1 \otimes Z,H]  =  (H \bullet y_1)Z \text{ . }$$

By \eqref{comm} in the proof of Proposition 14, $$[y_1 \otimes Z,y_i] = {(-1)}^{|Z|} y_1 \otimes (Z \bullet y_i)  \text{ . }     $$

Therefore, $$y_m \wedge ... \wedge [y_1 \otimes Z,y_i] \wedge ...
\wedge y_1 = 0 \text{   } \forall \text{  } i \neq 1 \text{ . } $$

Thus, $$ \text{ad}(L) \otimes \wedge^mV[m] (y_m \wedge.... \wedge
y_1 \otimes H \otimes x_1 \wedge ... \wedge x_m) $$ $$= y_m \wedge
... \wedge y_1 \otimes ((H \bullet y_1)Z +
{(-1)}^{|H|(|Z|-1)}{(-1)}^{|Z|} (Z \bullet y_1)H ) \otimes x_1
\wedge .. \wedge x_m $$
$$  = y_m \wedge ... \wedge y_1 \otimes  ((H \bullet
y_1)Z + {(-1)}^{|H|(|Z|-1)}{(-1)}^{|Z|} {(-1)}^{|H|(|Z|-1)}H(Z
\bullet y_1)) \otimes x_1 \wedge .. \wedge x_m $$
$$ = y_m \wedge ... \wedge y_1 \otimes  ((H \bullet
y_1)Z +{(-1)}^{|Z|} H(Z \bullet y_1)) \otimes x_1 \wedge .. \wedge
x_m \text{ . }$$

Note that $$(HZ \bullet y_1) = H (Z \bullet y_1) + {(-1)}^{|Z|}(H
\bullet y_1)Z \text{ . }$$

Further recall that $$\textbf{1}^m(H) = {(-1)}^{|H|m} y_m \wedge....
\wedge y_1 \otimes H \otimes x_1 \wedge ... \wedge x_m \text{ . }$$

 Therefore,
$$ {(-1)}^{|L|m} \text{ad}(L) \otimes \wedge^mV[m] (\textbf{1}^m(H)) = $$
$$ {(-1)}^{(|H|+|Z|-1)m} \text{ad}(L) \otimes \wedge^mV[m] (y_m \wedge.... \wedge y_1 \otimes H \otimes x_1 \wedge ... \wedge x_m) $$
$$=  {(-1)}^{(|H|+|Z|-1)m} y_m \wedge ... \wedge y_1 \otimes  ((H \bullet
y_1)Z +{(-1)}^{|Z|} H(Z \bullet y_1)) \otimes x_1 \wedge .. \wedge x_m $$
$$ = - {(-1)}^{(|H|+|Z|-1)m} {(-1)}^{|Z|+1} y_m \wedge ... \wedge y_1 \otimes (HZ \bullet y_1) \otimes x_1 \wedge .. \wedge x_m $$
$$ = - \textbf{1}^m (L^+(H)) \text{ . } $$

This proves the desired lemma. \\

\end{proof}

\textbf{Remark 4 :} Lemma 1 seems to be a phenomenon occurring in purely odd super-geometry only. Let $M_V$ be the super-manifold whose
ring of functions is
$\sss(V[1])$. The pairing $\langle,\rangle$ on $\sss(V[1])$ is the pairing
$$\langle f,g \rangle = \int_{M_V} fg \text{  } \forall \text{ } f,g \in \sss(V[1]) \text{ . }$$
 By $\int_{M_V}$, we of course mean a {\it Berezinian integral}. We can think of the usual geometric analog of $\sss(V[1])$ to be the ring of
 compactly supported functions on a smooth oriented manifold $M$. The analog of the pairing $\langle,\rangle$ on $\sss(V[1])$ is the
 pairing $$(f,g) \leadsto \int_M fgd\mu $$ where $d\mu$ is the measure arising out of a volume form on $M$.On $M_V$ there is a constant top-order differential operator $\partial$ that is unique upto scalar. Lemma
 1 then says that if $D$ is a differential operator on $M_V$ that is purely of first order, and if $f$ is a function on $M_V$, the Lie bracket
 of $D$ with $f\partial$ is $\pm (D^+f) \partial$ where $D^+$ is the adjoint of $D$. In the usual geometric setting, the analog of $\partial$
 would be a global, nowhere vanishing section of the top wedge power of the tangent bundle of $M$. Even if such a section exists on $M$, the
 analog of Lemma 1 does not hold even in the $1$-dimensional case. For example, if $M = \mathbb R$, then the adjoint of the operator $f
 \frac{d}{dx} $ is the operator $-(\frac{df}{dx} + f\frac{d}{dx})$ . This follows from the standard integration by parts. But ,
 $$[f\frac{d}{dx},g\frac{d}{dx}] = f\frac{dg}{dx} - g\frac{df}{dx} \neq \pm (-g\frac{df}{dx} - f\frac{dg}{dx}) \text{ . } $$ \\
We get back to proving more propositions in linear algebra that we require for future use. \\

Choose a basis $\{x_1,..,x_m\}$ of $V$ and a basis $\{y_1,...,y_m\}$
of $V^*$ dual to $\{x_1,...,x_m\}$. Define $\textbf{k}_{V}:
\sss(V[1]) \rar \enn(\sss(V^*[-1]))$ by the formulae
$$ \textbf{k}_V(x_i)(y_j) = \delta_{ij}$$ $$\textbf{k}_V(x_i)(1)=0 $$
$$ \textbf{k}_V(x_i)(Y_1 \wedge Y_2) = \textbf{k}_V(x_i)(Y_1)
+{(-1)}^{|Y_1|} Y_1 \wedge \textbf{k}_V(x_i)(Y_2) $$
$$ \textbf{k}_V(X_1 \wedge X_2)(Y) =
\textbf{k}_V(X_1)(\textbf{k}_{V}(X_2)(Y)) \text{ . } $$

To simplify notation, we shall denote $\textbf{k}_{V}(Z)(Y)$ by $(Z
|Y)$ \\

\textbf{Remark - :} The map $\textbf{j}_{V^*}$ identifies an element
$Y$ of $\sss(V^*[-1])$ with the operation "differentiation on the
right by $Y$". The map $\textbf{k}_V$ identifies an element $Z$ of
$\sss(V[1])$ with "differentiation on the left by $Z$".\\

We have the following proposition : \\

\begin{prop} Let $Y \in \sss(V^*[-1])$ and let $Z \in \sss(V[1])$ . Then, $$ \pi_0( F_r(\textbf{j}_{V^*}(Y) \circ \textbf{i}_{V}(Z))) =
(Z|Y) \text{ . }$$
\end{prop}

\begin{proof} This is again proven by a direct computation. We may
assume without loss of generality that $Y$ and $Z$ are homogenous.
\\

Choose a basis $\{x_1,...,x_m\}$ of $V$ and a basis
$\{y_1,...,y_m\}$ of $V^*$ dual to $\{x_1,...,x_m\}$. Let $H \in
\sss(V[1])$. Let $Z = x_S$ and let $Y= y_T$ for some $S,T \subset
\{1,...,m\}$. If $S$ and $T$ are disjoint and if $S$ is nonempty,
$\textbf{j}_{V^*}(Y)$ and $\textbf{i}_V(Z)$ commute upto sign. It
follows that $\pi_0(F_r(\textbf{j}_{V^*}(Y) \circ \textbf{i}_V(Z)))
= 0$ . Also, $(Z|Y)=0$ if $S$ and $T$ are disjoint and $S$ is
nonempty. If $S$ is empty, $x_S=1$ by convention. Therefore,
$\pi_0(F_r(\textbf{j}_{V^*}(Y) \circ \textbf{i}_V(Z)))= Y$ and
$(Z|Y)=Y$. We therefore, have to
prove this proposition for the case when $S$ and $T$ are not disjoint. \\

Let $S$ and $T$ be arbitrary subsets of $\{1,..,m\}$. Suppose that $j \not\in S \cup T$. \\

$$ \textbf{j}_{V^*}(y_j \wedge Y) \circ \textbf{i}_{V}(Z \wedge x_j)
(H) = (HZ \wedge x_j \bullet y_j \wedge Y) = (HZ \bullet Y) -  ((HZ
\bullet y_j) \wedge x_j \bullet Y) \text{ . }$$

Since $j$ is not in $T$. It follows that $$Z' \wedge x_j \bullet Y =
\pm (Z' \bullet Y) \wedge x_j \text{ . }$$ Therefore $$
\textbf{j}_{V^*}(y_j \wedge Y) \circ \textbf{i}_{V}(Z \wedge x_j)
(H)  = (HZ \bullet Y)  \pm ((HZ \bullet y_j) \bullet Y) \wedge x_j
\text{ . } $$

It follows that $$\pi_0(F_r(\textbf{j}_{V^*}(y_j \wedge Y) \circ
\textbf{i}_{V}(Z \wedge x_j))) = \pi_0(F_r(\textbf{j}_{V^*}(Y) \circ
\textbf{i}_{V}(Z))) $$ . Note that $$(Z \wedge x_j | y_j \wedge Y) =
(Z | Y) $$ since $j \not\in S \cup T$. The desired proposition
follows for homogenous $Z=x_S$ and $Y=y_T$ by induction on $|S \cap
T|$ . For general $Y$ and $Z$ the proposition follows from the fact
that the maps
$$\sss(V^*[-1]) \times \sss(V[1]) \rar \sss(V^*[-1]) $$
$$(Y,Z) \leadsto \pi_0(F_r(\textbf{j}_{V^*}(Y) \circ
\textbf{i}_{V}(Z))) $$ and $$\sss(V^*[-1]) \times \sss(V[1]) \rar
\sss(V^*[-1]) $$ $$ (Y,Z) \leadsto (Z|Y) $$ are both $\fld$-bi-linear. \\

\end{proof}

Let $(\text{ - }||\text{ - })_V:\sss(V[1]) \otimes \sss(V^*[-1])
\rar \fld$ denote map $$Z \otimes Y \leadsto p_0(Z | Y)$$ where
$p_0:\sss(V^*[-1]) \rar \fld$ denotes the projection to the degree
$0$ direct summand. \\

 Let $\gamma_V: \sss(V[1]) \rar \sss(V^*[-1]) \otimes \wedge^mV[m]$
be the isomorphism such that
$$\pi_m(\text{ - } \wedge \text{ - })_V = [(\text{ - }|| \text{ - })_V
\otimes \wedge^mV[m]] \circ (\sss(V[1]) \otimes \gamma_V) \text{ .
}$$ Let $\zeta_V$ denote $\gamma_V^{-1}$.

\begin{prop}
If $Z,W \in \sss(V[1])$, then $$\zeta_V(\{(Z | \text{ - }) \otimes
\wedge^mV[m]\}(\gamma(W))) = Z \wedge W \text{ . }$$
\end{prop}

\begin{proof}
This proposition is verified by a direct computation. Choose a basis
$\{x_1,...,x_m\}$ of $V$ and a basis $\{y_1,...,y_m\}$ of $V^*$ dual
to $\{x_1,...,x_m\} $. Without loss of generality, $Z = x_1 \wedge
... \wedge x_k$ and $W = x_{l+1} \wedge .. \wedge x_m$. \\

Then, $$\gamma_V(W) = y_l \wedge... \wedge y_1 \otimes x_1 \wedge ..
\wedge x_m $$
$$((Z | \text{ - }) \otimes
\wedge^mV[m])(\gamma_V(W)) = (x_1 \wedge ... \wedge x_k | y_l \wedge
...\wedge y_1) \otimes x_1 \wedge .... \wedge x_m \text{ . }$$

If $k > l$ then $(x_1 \wedge ... \wedge x_k | y_1 \wedge ...\wedge y_l) = 0$ and $x_1 \wedge ... \wedge x_k \wedge x_{l+1} \wedge ... \wedge x_m
= 0$, proving this proposition. We may thus assume that $k \leq l$. Then,

$$((Z | \text{ - }) \otimes
\wedge^mV[m])(\gamma_V(W)) =  (x_1 \wedge ... \wedge x_k |y_l \wedge ... \wedge y_1) \otimes x_1 \wedge ... \wedge x_m $$
$$ = {(-1)}^{k(l-k)} y_l \wedge ... \wedge y_{k+1} \otimes x_1
\wedge  .. \wedge x_m $$
$$ = \gamma_V(x_1 \wedge .. \wedge x_k \wedge x_{l+1} \wedge .. \wedge
x_m) \text{ . } $$

This proves the desired proposition. \\

\end{proof}

Let $J_V$ denote the endomorphism of $\sss(V[1])$ multiplying $\wedge^iV[i]$ by ${(-1)}^i$. We also have the following proposition \\
Let $\tau: \wedge^mV[m] \otimes \sss(V^*[-1]) \otimes \wedge^mV[m] \rar  \sss(V^*[-1]) \otimes \wedge^mV[m] \otimes \wedge^mV[m]$ denote the map
that interchanges $\wedge^mV[m]$ and $\sss(V^*[-1]) \otimes \wedge^mV[m]$. Let \\ $\simeq:\wedge^mV[m] \otimes \wedge^m V^*[-m] \rar \fld$
denote the map $x_1 \wedge ... \wedge x_m \otimes
y_m \wedge ... \wedge y_1 \leadsto 1 $. \\

\begin{prop}
The following diagram commutes : \\

$$\begin{CD}
\wedge^mV[m] \otimes \sss(V^*[-1]) \otimes \wedge^mV[m] \otimes
\wedge^mV^*[-m]  @>{(\text{ - } \bullet \text{ - }) \otimes
\simeq}>>  \sss(V[1])\\
 @VV{\tau \otimes \wedge^mV^*[-m]}V   @V{J_V}VV \\
 \sss(V^*[-1]) \otimes \wedge^mV[m] \otimes \wedge^mV[m] \otimes \wedge^mV^*[m]
 @>{\zeta_V
\otimes \simeq}>>
  \sss(V[1]) \\
\end{CD} $$
\end{prop}

\begin{proof}
This is verified by a direct computation as well. Let $\{x_i\}$, $\{y_i\}$  be as in the
proof of the previous proposition. \\

$$(x_1 \wedge .... \wedge x_m \bullet y_m \wedge ... \wedge y_{k+1} )
= x_1 \wedge ... \wedge x_k $$
$$\implies ((\text{ - }\bullet\text{ - }) \otimes \simeq)(x_1 \wedge .... \wedge x_m \otimes y_m \wedge ... \wedge
y_{k+1} \otimes x_1 \wedge .... \wedge x_m \otimes y_m \wedge ...
\wedge y_{1}) $$ $$ = x_1 \wedge ... \wedge x_k \text{ . } $$
$$ \tau ( x_1 \wedge .... \wedge x_m \otimes y_m \wedge ... \wedge
y_{k+1} \otimes x_1 \wedge ... \wedge x_m ) = {(-1)}^{mk}(y_m \wedge
... \wedge y_{k+1} \otimes x_1 \wedge ... \wedge x_m \otimes x_1
\wedge .... \wedge x_m) $$
$$\zeta_V(y_m \wedge ... \wedge
y_{k+1} \otimes x_1 \wedge ... \wedge x_m ) = {(-1)}^{k(m-k)} x_1 \wedge ... \wedge x_k $$
$$ \implies (\zeta_V \otimes \simeq) \circ (\tau \otimes \wedge^mV^*[-m]) ( x_1 \wedge .... \wedge x_m \otimes y_m \wedge ... \wedge
y_{k+1} \otimes x_1 \wedge ... \wedge x_m \otimes y_m \wedge...
\wedge y_1 ) $$ $$= {(-1)}^{-k^2} x_1 \wedge ... \wedge x_k =
J_V(x_1 \wedge ... \wedge x_k)  \text{ . }$$

This proves the desired proposition. \\

\end{proof}

Let $\textbf{C}_V$ denote the co-multiplication on $\sss(V[1])$.
Think of $\textbf{C}_V$ as an element of  \\ $\enn(\sss(V[1]))
\otimes \sss(V[1])$.

\begin{prop}

If $Y \in \sss(V^*[-1])$ then $$(\enn(\sss(V[1])) \otimes ( \text{ - } || Y) \circ \textbf{C}_V) = \textbf{j}_{V^*}(Y) \text{ . }$$

\end{prop}

\begin{proof}
This is yet another proposition that is verified by a direct computation. Let $\{x_i\}$,$\{y_i\}$ be as in the proof of the previous
proposition. Without loss of generality, $Y = y_1 \wedge .... \wedge y_k$. Then,

$$\textbf{C}_V(x_l \wedge .... \wedge x_1) = x_l \wedge .... \wedge
x_{k+1} \otimes x_k \wedge ... \wedge x_1 + \sum_{S \neq
\{l,...,k+1\}}  \pm x_S \otimes x_{\bar{S}} $$
$$\implies (\enn(\sss(V[1])) \otimes ( \text{ -
} || Y) \circ \textbf{C}_V)(x_l \wedge .... \wedge x_1) = x_{l}
\otimes .... \otimes x_{k+1} = \textbf{j}_{V^*}(Y)(x_l \wedge ..
\wedge x_1) \text{ . } $$

The desired proposition follows immediately from the above
computation.\\

\end{proof}

\subsection{Applying the linear algebra to $\sss(\Omega[1])$}

Let $\calg M$ be a locally free coherent $\strc$-module. We begin
with the remark that every proposition in the previous subsection
holds in $\bcc$ (and hence in $\dcat$ ) with $V$ replaced by $\calg
M$, $V^*$ replaced by ${\calg M}^* = {\calg
H}\text{om}_{\strc}(\calg M, \strc)$. We are interested in the case
when $\calg M = \Omega$. All graded $\strc$-modules in this subsection are
to be thought of as complexes of $\strc$-modules with $0$-differential. \\

{\it All $\fld$-vector spaces in this section are finite dimensional. In this subsection and in future sections, $n$ shall denote the dimension
of $X$.} \\

\textbf{An important point for the reader to note :}Each proposition
or lemma in this subsection is proven by proving it for an arbitrary
open subscheme $U$ of $X$ such that $\Omega$ is trivial over $U$.
Then, $\Omega |_U = \strcc{U} \otimes_{\fld} V$ and $T |_U =
\strcc{U} \otimes_{\fld} V^* $ for some finite dimensional
$\fld$-vector space $V$. After observing that every map involved in
the proposition/Lemma is $\strcc{U}$-linear, proving the
proposition/Lemma reduces to proving the corresponding
proposition/Lemma in Section 3.1. The Proposition/Lemma in Section
3.1 corresponding to each proposition here can be thought of as the
"local computation" required to prove the corresponding
proposition/lemma
in this subsection. \\

With the above announcement we can just state the propositions and lemma that we wish to state. \\

Note that there is a morphism $$\textbf{i}_{\Omega}: \sss(\Omega[1]) \rar \ennnn(\sss(\Omega[1]))$$ such that whenever $U$ is an open subscheme
on $X$ such that $\Omega \simeq \strcc{U} \otimes_\fld V$ for some $\fld$-vector space $V$, then
$$\textbf{i}_{\Omega} |_U = \textbf{i}_{V} \otimes_\fld \strcc{U} \text{ . }$$

We denote $\textbf{i}_{\Omega}$ by $\textbf{i}$ to simplify
notation. \\

Similarly, we have a morphism $$\textbf{j}_{T}:\sss(T[-1]) \rar
\ennnn(\sss(\Omega[1]))$$ such that whenever $U$ is an open
subscheme on $X$ such that $\Omega \simeq \strcc{U} \otimes_\fld V$
for some $\fld$-vector space $V$, then $$\textbf{j}_{T} |_U =
\textbf{j}_{V^*} \otimes_\fld \strcc{U} \text{ . } $$

We denote $\textbf{j}_{T}$ by $\textbf{j}$ to simplify notation. \\

The following proposition corresponds to Proposition 10. \\

\begin{prop}
The composite
$$\begin{CD} \sss(T[-1]) \otimes \sss(\Omega[1]) @>\textbf{j} \otimes
\textbf{i} >> \ennnn(\sss(\Omega[1])) \otimes \ennnn(\sss(\Omega[1])
@> \circ^{\text{op}} >> \ennnn(\sss(\Omega[1])) \end{CD} $$ is an
isomorphism in $\bcc$.
\end{prop}

\subsection*{Notation}  $G_r$ shall denote
the isomorphism in Proposition 20 and $F_r$ shall denote its
inverse. In addition, if $\tau: \sss(T[-1]) \otimes \sss(\Omega[1])
\rar \sss(\Omega[1]) \otimes \sss(T[-1])$ denotes the swap map the
composite
$$
\begin{CD}\sss(T[-1]) \otimes \sss(\Omega[1]) @> \tau>> \sss(\Omega[1])
\otimes \sss(T[-1]) @>\textbf{i} \circ^{\text{op}} \textbf{j}>> \ennnn(\sss(\Omega[1]))
\end{CD} $$ shall be denoted by $G_l$ .
 The inverse of $G_l$ will be denoted by $F_l$. \\

Let $\pi_k:\sss(\Omega[1]) \rar \Omega^k[k]$ denote the natural projection. Note that we have a non-degenerate pairing $\langle,\rangle$ on
 $\sss(\Omega[1])$. This is given by the composite
 $$\begin{CD} \sss(\Omega[1]) \otimes \sss(\Omega[1]) @>(\text{ - } \wedge \text{ - })>> \sss(\Omega[1])
 @>\pi_n>> \Omega^n[n] \\ \end{CD} \text{ . }$$

Recall the definition of the Adjoint of an element of
$\enn(\sss(V[1]))$ from Section 3.1. Let
$\textbf{A}_{V}:\enn(\sss(V[1])) \rar \enn(\sss(V[1])) $ denote the
map taking an element of $\enn(\sss(V[1]))$ to its adjoint. We have
a morphism $\textbf{A}_{\Omega}:\ennnn(\sss(\Omega[1])) \rar
\ennnn(\sss(\Omega[1])) $ in $\bcc$ such that if $U$ is an open
subscheme of $X$ such that $\Omega \simeq \strcc{U} \otimes_{\fld}
V$ for some $\fld$-vector space $V$, then
$$\textbf{A}_{\Omega} = \textbf{A}_{V} \otimes_{\fld} \strcc{U} \text{ . } $$

Similarly, the map $$\text{ev}_{V}: \sss(V[1]) \otimes
\enn(\sss(V[1])) $$ $$ Z \otimes \Phi \leadsto \Phi(Z) $$ yields a
map $$ \text{ev}:\sss(\Omega[1]) \otimes \ennnn(\sss(\Omega[1]))
\rar \sss(\Omega[1]) $$ such that $$\text{ev} |_U = \text{ev}_V
\otimes_{\fld} \strcc{U}$$ on any
open subscheme $U$ of $X$ such that $\Omega |_U = V \otimes_{\fld} \strcc{U}$ for some $\fld$-vector space $V$. \\

We denote $\textbf{A}_{\Omega}$ by $\textbf{A}$ to simplify notation. Let $\langle \text{ },\text{ev}^+ \rangle$ denote the composite
$$\begin{CD}
\sss(\Omega[1]) \otimes \sss(\Omega[1]) \otimes \ennnn(\sss(\Omega[1])\\
@VV{\sss(\Omega[1]) \otimes \sss(\Omega[1]) \otimes \textbf{A}}V\\
\sss(\Omega[1]) \otimes \sss(\Omega[1]) \otimes \ennnn(\sss(\Omega[1]) @>\sss(\Omega[1]) \otimes \text{ev}>> \sss(\Omega[1]) \otimes
\sss(\Omega[1])
 @>\langle \text{ },\text{ } \rangle>> \Omega^n[n] \end{CD} $$.
 By the corresponding fact for a finite dimensional $\fld$-vector space
$V$, the following diagram commutes in $\bcc$ (and hence in $\dcat$): \\

\begin{equation} \label{adjdef}
\begin{CD}
\sss(\Omega[1]) \otimes \ennnn(\sss(\Omega[1])) \otimes \sss(\Omega[1]) @> \sss(\Omega[1]) \otimes \tau >> \sss(\Omega[1]) \otimes \sss(\Omega[1]) \otimes \ennnn(\sss(\Omega[1])) \\
@VV{\text{ev} \otimes \sss(\Omega[1])}V     @V{\langle \text{ },\text{ev}^+ \rangle}VV\\
\sss(\Omega[1]) \otimes \sss(\Omega[1])    @> \langle,\rangle>> \Omega^n[n] \\
\end{CD}
\end{equation}

The map $\tau: \ennnn(\sss(\Omega[1])) \otimes \sss(\Omega[1]) \rar \sss(\Omega[1]) \otimes \ennnn(\sss(\Omega[1]))$ in the top row of the
above diagram swaps factors.\\

The following proposition corresponds to Proposition 11 in Section
3.1. \\

\begin{prop}1. To begin with,
$$ \textbf{A} \circ \textbf{A} =
\textbf{1}_{\ennnn(\sss(\Omega[1]))} \text{ . } \\ $$

2. If $\tau: \ennnn(\sss(\Omega[1])) \otimes \ennnn(\sss(\Omega[1]))
\rar  \ennnn(\sss(\Omega[1])) \otimes \ennnn(\sss(\Omega[1]))$
denotes the swap map, then the following diagram commutes in $\bcc$:

$$\begin{CD}
\ennnn(\sss(\Omega[1])) \otimes \ennnn(\sss(\Omega[1])) @>
(\textbf{A} \otimes \textbf{A}) \circ \tau >>
\ennnn(\sss(\Omega[1])) \otimes \ennnn(\sss(\Omega[1])) \\
@VV{ \circ}V       @V{\circ}VV \\
\ennnn(\sss(\Omega[1])) @> \textbf{A}>> \ennnn(\sss(\Omega[1])) \\
\end{CD} $$

3. Also ,

$$\textbf{A} \circ \textbf{i} = \textbf{i} \text{ . }$$

4. If $I:\sss(T[-1]) \rar \sss(T[-1])$ denotes the endomorphism
multiplying $\wedge^i T[-i]$ by ${(-1)}^i$ then

$$ \textbf{A} \circ \textbf{j} = \textbf{j} \circ I \text{ . } $$

5. Finally,

$$ \textbf{A} \circ G_l = G_r \circ (I \otimes \sss(\Omega[1])) \text{ . } $$
\end{prop}

\noindent Recall that $\pi_j:\sss(\Omega[1]) \rar \Omega^j[j]$
denotes the natural projection. We will denote the projection
$$ \sss(T[-1]) \otimes \pi_j: \sss(T[-1])  \otimes  \sss(\Omega[1]) \rar
 \sss(T[-1]) \otimes \Omega^j[j] $$ by $\pi_j$ itself. Let $I$ be as in Part 4, Proposition 21 . The
following Proposition corresponds to Proposition 12 of Section 3.1.
\\

\begin{prop}

$$\pi_0 \circ F_l = I \circ \pi_0 \circ F_r \circ \textbf{A} \text{ . }$$

\end{prop}

Denote the composite
$$ \begin{CD} {(\sss(T[-1]) \otimes \sss(\Omega[1]))}^{\otimes 2} @> G_r \circ G_r >> \ennnn(\sss(\Omega[1])) @>F_r>> \sss(T[-1]) \otimes \sss(\Omega[1])
\end{CD} $$ by $\text{m}$ to simplify notation.
The proposition below corresponds to Proposition 13 of Section 3.1. \\

\begin{prop}
As morphisms in $\bcc$, $$\pi_0 \circ \text{m} = \pi_0 \circ
\text{m} \circ (\pi_0 \otimes \sss(T[-1]) \otimes \sss(\Omega[1]))
\text{ . }$$
\end{prop}

Denote the direct summand $\wedge^iT[-i] \otimes \sss(\Omega[1])$ of $\ennnn(\sss(\Omega[1]))$ by $D_i$.  Note that the map $[,]_V
:\enn(\sss(V[1]))^{\otimes 2} \rar \enn(\sss(V[1]))$ extends to a morphism \\ $[,]: \ennnn(\sss(\Omega[1]))^{\otimes 2} \rar
\ennnn(\sss(\Omega[1]))$. If $U$ is an open subscheme on $X$ such that $\Omega \simeq \strcc{U} \otimes_{\fld} V$ for some $\fld$-vector space
$V$, then  $$[,]|_U = [,]_V \otimes_{\fld} \strcc{U} \text{ . } $$ The following proposition corresponds to Proposition 14 of Section 3.1.

\begin{prop}
The composite $$\begin{CD} D_1 \otimes D_1 @>>> \ennnn(\sss(\Omega[1]))^{\otimes 2} @>[,]>> \ennnn(\sss(\Omega[1])) \end{CD} $$ factors through
the inclusion of $D_1$ in $\ennnn(\sss(\Omega[1]))$.
\end{prop}

Let $V$ be a $\fld$-vector space of dimension $n$. Denote the direct
summands $V^*[-1] \otimes \sss(V[1])$ and $\wedge^n V^*[-n] \otimes
\sss(V[1])$ of $\sss(V^*[-1]) \otimes \sss(V[1])$ by $\bar{D_1}$ and
$\bar{D_n}$ respectively, {\it unlike} in Section 3.1. Recall that
in Section 3.1 we defined maps
$$\text{ad}(L): \wedge^n V^*[-n] \otimes \sss(V[1]) \rar \wedge^n V^*[-n] \otimes \sss(V[1])
$$ $$\bar{\text{ad}}(L): \bar{D_1}^{\otimes n} \rar \bar{D_1}^{\otimes n}$$ for any $L \in \bar{D_1}$. These yield us maps $$\text{ad}_V: \wedge^n V^*[-n] \otimes \sss(V[1]) \otimes \bar{D_1} \rar \wedge^n V^*[-n] \otimes
\sss(V[1])$$
$$ Y \otimes L \leadsto \text{ad}(L)(Y) $$
$$\bar{\text{ad}}_V: \bar{D_1}^{\otimes n} \otimes \bar{D_1} \rar \bar{D_1}^{\otimes n} $$ $$ W \otimes L \leadsto \bar{\text{ad}}(L)(W) \text{ . } $$
These yield morphisms $$\text{ad}: \wedge^nT[-n] \otimes
\sss(\Omega[1]) \otimes D_1 \rar \wedge^nT[-n] \otimes
\sss(\Omega[1]) $$
$$\bar{\text{ad}}:D_1^{\otimes n} \otimes D_1 \rar D_1^{\otimes n}$$ in $\bcc$. If $U$ is an open subscheme on $X$ such that
 $\Omega \simeq \strcc{U} \otimes_{\fld} V$ for some $\fld$-vector
space $V$, then
$$ \text{ad}|_U = \text{ad}_V \otimes_{\fld} \strcc{U} $$
$$ \bar{\text{ad}} |_U = \bar{\text{ad}}_V \otimes_{\fld} \strcc{U} \text{ . } $$

The map $p_V:\bar{D_1}^{\otimes n} \rar \bar{D_n} $ also yields a map $$p:D_1^{\otimes n} \rar \wedge^nT[-n] \otimes \sss(\Omega[1]) $$ such
that if $\Omega \simeq V \otimes_{\fld} \strcc{U}$ for some open subscheme $U$ of $X$ and some $\fld$-vector space $V$, then,
$$p = p_V  \otimes_{\fld} \strcc{U} \text{ . } $$

Note that the dimension $n$ of $X$ is the rank of $\Omega$ as well. Let $\textbf{1}^n: \strc \rar \wedge^n T[-n] \otimes \Omega^n[n]$ denote the
map dual to the evaluation map. There are maps \\ $\tau: \wedge^n T[-n] \otimes \Omega^n[n] \otimes \sss(\Omega[1]) \otimes D_1 \rar \wedge^n
T[-n] \otimes \sss(\Omega[1]) \otimes D_1  \otimes \Omega^n[n] $ and  \\ $\tau':\wedge^n T[-n] \otimes \Omega^n[n] \otimes \sss(\Omega[1]) \rar
\wedge^n T[-n] \otimes \sss(\Omega[1]) \otimes \Omega^n[n]$. $\tau$ is obtained by swapping $\Omega^n[n]$ and $\sss(\Omega[1]) \otimes D_1$.
$\tau'$ is obtained by swapping $\Omega^n[n]$ and $\sss(\Omega[1])$. Denote the composites
$$\tau \circ (\textbf{1}^n \otimes \sss(\Omega[1]) \otimes D_1 )$$ and $$\tau' \circ (\textbf{1}^n \otimes \sss(\Omega[1]))$$
 by $\textbf{1}^n$. The following proposition corresponds to Proposition 15 of Section 3.1. \\

\begin{prop}
The following diagram commutes in $\bcc$: \\

$$\begin{CD}
D_1^{\otimes n} \otimes D_1 @>p \otimes D_1>> D_n \otimes D_1 \\
@VV{\bar{\text{ad}}}V    @V{\text{ad}}VV \\
D_1^{\otimes n} @>p>> D_n \\
\end{CD} $$

\end{prop}

  The following Lemma corresponds
to Lemma 1 of Section 3.1. In the following Lemma, $\textbf{A}$
denotes the composite $$\begin{CD} D_1 @>>> \ennnn(\sss(\Omega[1]))
@>\textbf{A}>> \ennnn(\sss(\Omega[1])) \end{CD} \text{ . }$$

\begin{lem}
The following diagram commutes in $\bcc$: \\

$$ \begin{CD}
\sss(\Omega[1]) \otimes D_1 @> \textbf{1}^n >> \wedge^nT[-n] \otimes \sss(\Omega[1]) \otimes D_1 \otimes
\Omega^n[n] \\
@VV{-\text{ev} \circ (\sss(\Omega[1]) \otimes \textbf{A})}V     @V{\text{ad} \otimes \Omega^n[n]}VV \\
\sss(\Omega[1]) @> \textbf{1}^n >>    \wedge^nT[-n] \otimes \sss(\Omega[1]) \otimes \Omega^n[n] \\
\end{CD} $$

\end{lem}

The map $$(\text{ - }|\text{ - })_V :\sss(V[1]) \otimes
\sss(V^*[-1]) \rar \sss(V^*[-1]) $$ $$Z \otimes Y \leadsto (Z|Y) $$
yields a  map
$$(\text{ - }|\text{ - }):\sss(\Omega[1]) \otimes \sss(T[-1]) \rar \sss(T[-1]) $$. As usual, if $U$ is an open subscheme on $X$ such
that $\Omega \simeq \strcc{U} \otimes_\fld V$ for some $\fld$-vector space $V$, then
$$( \text{ - }|\text{ - }) = (\text{ - }|\text{ - })_V \otimes_\fld \strcc{U} \text{ . }$$

Recall that the composition product on $\ennnn(\sss(\Omega[1]))$ was denoted by $\circ$. Let $\circ^{\text{op}}$ denote the product on
$\ennnn(\sss(\Omega[1]))^{\text{op}}$. The following proposition corresponds to Proposition 16 of Section 3.1.\\

\begin{prop}

$$\pi_0 \circ F_r \circ (\textbf{i} \circ^{\text{op}} \textbf{j}) = (\text{ - }|\text{ - }): \sss(\Omega[1]) \otimes \sss(T[-1]) \rar \sss(T[-1]) \text{ . } $$

\end{prop}

Let $(\text{ - }||\text{ - }): \sss(\Omega[1]) \otimes \sss(T[-1]) \rar \strc$ denote the map such that for any open $U \subset X$ such that
$\Omega \simeq V \otimes_{\fld} \strcc{U}$ for some $\fld$-vector space $V$,

$$ (\text{ - }||\text{ - }) = (\text{ - }||\text{ - })_V
\otimes_{\fld} \strcc{U} \text{ . } $$

Let $\gamma:\sss(\Omega[1]) \rar \sss(T[-1]) \otimes S_X$ be the
isomorphism such that
$$\pi_n \circ (\text{ - } \wedge \text{ - }) = ((\text{ - } || \text{
- }) \otimes S_X) \circ (\sss(\Omega[1]) \otimes \gamma) \text{ .
}$$ Let $\zeta$ denote the inverse of $\gamma$.

The following proposition corresponds to Proposition 17 of Section
3.1. \\

\begin{prop}
$$\zeta([(\text{ - }|\text{ - }) \otimes S_X] \circ[ \sss(\Omega[1]) \otimes \gamma]) =
(\text{ - } \wedge \text{ - }): \sss(\Omega[1]) \otimes
\sss(\Omega[1]) \rar \sss(\Omega[1]) \text{ . }$$
\end{prop}

Let $\tau: S_X \otimes \sss(T[-1]) \otimes S_X \rar \sss(T[-1]) \otimes S_X \otimes S_X$ denote the map swapping $S_X$ with $\sss(T[-1]) \otimes
S_X$. Let $\simeq$ denote the identification of $S_X \otimes S_X^{-1}$ with $\strc$. The following proposition corresponds to Proposition 18 of
Section 3.1. Recall that $J$ is the endomorphism of $\sss(\Omega[1])$ that
multiplies $\Omega^j[j]$ by ${(-1)}^j$. \\

\begin{prop}
The following diagram commutes in $\bcc$: \\

$$\begin{CD}
S_X \otimes \sss(T[-1]) \otimes S_X \otimes S_X^{-1} @>\tau \otimes
S_X^{-1} >> \sss(T[-1]) \otimes S_X \otimes S_X \otimes S_X^{-1} \\
@VV{(\text{ - } \bullet \text{ - }) \otimes \simeq}V @V{\zeta
\otimes \simeq}VV \\
\sss(\Omega[1]) @>J>> \sss(\Omega[1])\\
\end{CD} $$

\end{prop}

Let $\textbf{C}_{\Omega}$ denote the co-product on $\sss(\Omega[1])$. Think of $\textbf{C}_{\Omega}$ as a morphism in $\dcat$ from $\strc$ to
$\ennnn(\sss(\Omega[1])) \otimes \sss(\Omega[1])$. The following proposition corresponds to
Proposition 19 of Section 3.1: \\

\begin{prop} $$\ennnn(\sss(\Omega[1])) \otimes (\text{ - } || \text{
-}) \circ \textbf{C}_{\Omega} = \textbf{j}: \sss(T[-1]) \rar
\ennnn(\sss(\Omega[1])) \text{ . } $$ \end{prop}

\section{The adjoint of $\Phi_R$}

This section is a continuation of Section 2. All maps in this
section are in $\dcat$ unless explicitly stated otherwise. Both in
this section and the next, we use results from Section 3.2. All
morphisms in $\bcc$ described in Section 3.2 induce morphisms in
$\dcat$. The diagrams that were shown to commute in $\bcc$ in Section 3.2 also commute in $\dcat$. \\

\subsection{Stating the main lemma of this section} Let $\Phi_L$ and $\Phi_R$ be as in Section 2.By Theorem 2' , the following diagrams commute
in $\dcat$ :

$$\begin{CD}
\hc^{\bullet}(X) @>\alpha_L>> \hc^{\bullet}(X) \otimes \Omega[1] \\
@VV\hkr V   @V \hkr \otimes \Omega[1] VV \\
\sss(\Omega[1]) @> \Phi_L >> \sss(\Omega[1]) \otimes \Omega[1] \\
\end{CD} $$

$$\begin{CD}
\hc^{\bullet}(X) @>\alpha_R>> \hc^{\bullet}(X) \otimes \Omega[1] \\
@VV\hkr V   @V \hkr \otimes \Omega[1] VV \\
\sss(\Omega[1]) @> \Phi_R >> \sss(\Omega[1]) \otimes \Omega[1] \\
\end{CD} $$

$\Phi_R$ can be thought of as an element of $\text{Hom}_{\dcat}(\strc, \ennnn(\sss(\Omega[1])) \otimes \Omega[1])$. Recall the
endomorphism $\textbf{A}$ of $\ennnn(\sss(\Omega[1]))$ (defined in Section 3.2) which takes a section of $\ennnn(\sss(\Omega[1]))$ to its adjoint.\\

 \textbf{Definition} The adjoint $\Phi_R^+$ of $\Phi_R$ is the element $(\textbf{A} \otimes \Omega[1]) \circ \Phi_R$
 of \\ $\text{Hom}_{\dcat}(\strc,\ennnn(\sss(\Omega[1])) \otimes \Omega[1])$. \\

 This section is devoted to finding an explicit formula for $\Phi_R^+$. \\

 Note that the Atiyah class of the tangent bundle $T$ of $X$ yields a morphism \\ $\text{At}_T: \Omega[1] \rar \Omega[1] \otimes \Omega[1]$.
 Let $p: \Omega[1]^{\otimes i} \rar \Omega^i[i]$ denote the standard projection. Let $$\text{At}_T^i: \Omega[1] \rar \sss(\Omega[1]) \otimes
 \Omega[1]$$ denote the composite $$ (p \otimes \Omega[1]) \circ (\Omega[1]^{\otimes i-1} \otimes \text{At}_T) \circ ... \circ \text{At}_T \text{ . }$$
 Then if $\frac{z}{\text{e}^z -1} = 1 + \sum_i c_i Z^i $ the map
 $$\frac{\text{At}_T}{\text{exp}(\text{At}_T)-1} := \textbf{1} + \sum_i c_i \text{At}_T^i : \Omega[1] \rar \sss(\Omega[1]) \otimes \Omega[1] $$
 makes sense.
 Denote the element $\text{det}(\frac{\text{At}_T}{\text{exp}(\text{At}_T)-1}) \in \text{Hom}_{\dcat}(\strc,\sss(\Omega[1]))$ by
 $\textbf{f}$.\\
Note that $\textbf{i}(\textbf{f})$ and $\textbf{i}(\textbf{f}^{-1})$ are elements of $\text{Hom}_{\dcat}(\strc,\ennnn(\sss(\Omega[1])))$.\\

Let $ \textbf{i}(\textbf{f})  \circ \Phi_R \circ \textbf{i}(\textbf{f}^{-1})$ denote the composite
$$\begin{CD} \strc \\
@V{\Phi_R}VV\\
 \ennnn(\sss(\Omega[1])) \otimes \Omega[1] \simeq \strc \otimes \ennnn(\sss(\Omega[1])) \otimes \strc \otimes \Omega[1] \\
 @V{\textbf{i}(\textbf{f}) \otimes \ennnn(\sss(\Omega[1])) \otimes \textbf{i}(\textbf{f}^{-1}) \otimes \Omega[1]}VV \\
 \ennnn(\sss(\Omega[1])) \otimes \ennnn(\sss(\Omega[1])) \otimes \ennnn(\sss(\Omega[1])) \otimes \Omega[1] @>\circ \otimes \Omega[1]>>
 \ennnn(\sss(\Omega[1])) \otimes \Omega[1] \end{CD} $$

The following lemma is the main lemma of this section. \\

\begin{lem}

$$\Phi_R^+ = \textbf{i}(\textbf{f}) \circ \Phi_R \circ \textbf{i}(\textbf{f}^{-1}) \text{ . }$$

\end{lem}

The proof of this lemma requires further preparation.  The following subsection is devoted to a key lemma (Lemma 4) used to prove Lemma 3.
The proof of Lemma 3 itself is at the end of this section (in Section 4.4). \\

\subsection{ Comparing two "sections" of $D_n \otimes \Omega^n[n]$}

Let $p_k: \sss(T[-1]) \rar \wedge^k T[-k]$ be the standard
projection. View $\Phi_L$ as an element of $\text{Hom}_{\dcat}(\strc, \ennnn(\sss(\Omega[1])) \otimes \Omega[1])$. We have the following
proposition. \\

\begin{prop}
$$(p_k \otimes \sss(\Omega[1]) \otimes \Omega[1]) \circ (F_r \otimes \Omega[1])(\Phi_L) = 0 \text{  }
\forall \text{   }
 k \neq 1 \text{ . }$$
\end{prop}

\textbf{Remark :} The above proposition just states that $\Phi_L$
can be thought of as an element of $\text{Hom}_{\dcat}(\strc, D_1
\otimes \Omega[1])$. Recall that $D_1$ is like the space of "purely
first order differential operators on $\sss(\Omega[1])$". In
principle, this proposition should follow from Proposition 7 and
Theorem 2'. I however, give a concrete proof below that uses the
definition of $\Phi_L$ from Section 2 since I cant see
how the above proposition follows immediately from Proposition 7. \\

\begin{proof}

Let $\tilde{\textbf{C}}: \sss(T[-1]) \rar \sss(T[-1]) \otimes
\sss(T[-1])$ denote co-product on $\sss(T[-1]$. Denote the composite
$(\sss(T[-1]) \otimes p_1) \circ \tilde{\textbf{C}} $ by
$\hat{\textbf{C}}$. Denote the wedge product on $\sss(\Omega[1])$ by
$\mu$ {\it in this proof only}. Note that the dual of the map
$$\bar{\omega}:\sss(\Omega[1]) \otimes \Omega[1] \rar \sss(\Omega[1])
\otimes \Omega[1] $$ defined in Section 2 is the composite
$$ (\sss(T[-1]) \otimes \text{At}_{T}) \circ (\hat{\textbf{C}}
\otimes T[-1]) \text{ . } $$ It follows that $$\bar{\omega} = (\mu
\otimes \Omega[1]) \circ (\sss(\Omega[1]) \otimes \text{At}_{T})
\text{ . } $$ In the latter composite, $\text{At}_{T}$ is thought of
as a morphism in $\dcat$ from $\Omega[1]$ to $\Omega[1] \otimes
\Omega[1]$. Therefore, \begin{equation} \label{omega} \bar{\omega}^i
= (\mu \otimes \Omega[1]) \circ (\sss(\Omega[1]) \otimes
\text{At}_T^i)
\end{equation} where $\text{At}_{T}^i$ is
as in the previous subsection. \\

It follows from \eqref{omega} that $ (p_k \otimes \sss(\Omega[1])
\otimes \Omega[1]) \circ (F_r \otimes \Omega[1])(\bar{\omega}^i
\circ \bar{\textbf{C}}) $ is given by the composite

$$\begin{CD} \strc @> \bar{\textbf{C}} >> \ennnn(\sss(\Omega[1])) \otimes
\Omega[1] @> ((p_k \otimes \sss(\Omega[1])) \circ F_r) \otimes
\text{At}_{T}^i
>> \wedge^k T[-k] \otimes \sss(\Omega[1]) \otimes \Omega^i[i]
\otimes \Omega[1]
\end{CD} $$ $$ \begin{CD} \wedge^k T[-k] \otimes \sss(\Omega[1]) \otimes
\Omega^i[i] \otimes \Omega[1] @>\wedge^k T[-k] \otimes \mu \otimes \Omega[1] >> \wedge^kT[-k] \otimes \sss(\Omega[1]) \otimes \Omega[1]
\end{CD} $$

Note that from the above description, proving that $(((p_k \otimes
\sss(\Omega[1])) \circ F_r)
 \otimes \Omega[1]) \circ \bar{\textbf{C}} = 0$ will imply that
$$(p_k \otimes \sss(\Omega[1]) \otimes \Omega[1]) \circ (F_r \otimes
\Omega[1])(\bar{\omega}^i \circ \bar{\textbf{C}}) = 0 $$. The desired proposition will then follow the fact that
$$\Phi_L = \sum_i c_i \bar{\omega}^i \circ \bar{\textbf{C}} $$ where
$\sum_i c_iz^i = \frac{z}{\text{e}^z -1} $. \\

Proving this proposition has therefore been reduced to proving that
\begin{equation} \label{calc} (((p_k \otimes \sss(\Omega[1])) \circ F_r)
 \otimes \Omega[1]) \circ \bar{\textbf{C}} = 0 \text{      } \forall k \neq 1 \text{ . }\end{equation}
 Let $U$ be an affine open subscheme of $X$ such that $\Omega \simeq
 V \otimes_K \strcc{U}$ for some $K$-vector space $V$. Proving
 \eqref{calc} reduces to proving that
 \begin{equation} \label{calcred} (((p_k \otimes \sss(V[1])) \circ F_r)
 \otimes V[1]) \circ \textbf{C}_V = 0 \text{     } \forall k \neq 1 \end{equation} where
 $\textbf{C}_V \in \enn(\sss(V[1])) \otimes V[1] \simeq \text{Hom}_{\fld}(\sss(V[1]),\sss(V[1]) \otimes V[1])$
 is the map
 $$v_1 \wedge ... \wedge v_j \leadsto \sum_i {(-1)}^{j-i}
 \widehat{v_1 \wedge ... i... \wedge v_j} \otimes v_i $$
 , $F_r$ is as in Section 3.1 and $p_k:\sss(V^*[-1]) \rar \wedge^k
 V^*[-k]$ denotes the standard projection. \eqref{calcred} however,
 follows from the fact that $$\textbf{C}_V = \sum_{j=1}^{j=n}
 \textbf{j}_{V}(y_j) \otimes x_j$$ $$ \implies  (F_r
 \otimes V[1]) \circ  \textbf{C}_V = \sum_{j=1}^{j=n} y_j \otimes x_j $$  for any bases
 $\{x_1,...,x_n\}$ of $V$ and $\{y_1,...,y_n\}$ of $V^*$ dual to
 each  other.

\end{proof}

The proof of Proposition 30 also helps us understand the map
$\Phi_L$ better: Let \\ $\text{id}:\strc \rar T[-1] \otimes
\Omega[1]$ denote the dual of the evaluation map from $\Omega[1]
\otimes T[-1]$ to $\strc$.

\begin{prop} As an element of $\text{Hom}_{\dcat}(\strc, D_1 \otimes \Omega[1])$, $\Phi_L$ is given by the composite
$$\begin{CD} \strc @>\text{id}>> T[-1] \otimes \Omega[1] @>T[-1] \otimes \frac{\text{At}_T}{\text{exp}(\text{At}_T)-1}>> T[-1] \otimes \sss(\Omega[1])
\otimes \Omega[1] \end{CD} \text{ . } $$
\end{prop}

\begin{proof} By the proof of Proposition 30, $\bar{\textbf{C}}:\sss(\Omega[1]) \rar \sss(\Omega[1]) \otimes \Omega[1]$ is given by the composite
$$\begin{CD} \sss(\Omega[1]) \otimes \strc @>\sss(\Omega[1]) \otimes \text{id}>> \sss(\Omega[1]) \otimes T[-1] \otimes \Omega[1] @> (\text{ - }
\bullet \text{ - }) \otimes \Omega[1]>> \sss(\Omega[1]) \otimes
\Omega[1] \end{CD} \text{ . }$$

Note that by definition, $\Phi_L =
\frac{\bar{\omega}}{\text{e}^{\bar{\omega}}-1} \circ
\bar{\textbf{C}}$. Also, equation \eqref{omega} in the proof of
Proposition 30 says that $\bar{\omega}^i$ is given by the composite
$$\begin{CD} \sss(\Omega[1]) \otimes \Omega[1] @>\sss(\Omega[1]) \otimes \text{At}_T^i>> \sss(\Omega[1]) \otimes \sss(\Omega[1]) \otimes
\Omega[1] @>(\text{ - } \wedge \text{ - }) \otimes \Omega[1]>>
\sss(\Omega[1]) \otimes \Omega[1] \end{CD} $$.

Therefore $\frac{\bar{\omega}}{\text{e}^{\bar{\omega}}-1}$ is given by the composite
$$\begin{CD} \sss(\Omega[1]) \otimes \Omega[1] @>\sss(\Omega[1]) \otimes \frac{\text{At}_T}{\text{exp}(\text{At}_T)-1} >> \sss(\Omega[1]) \otimes \sss(\Omega[1]) \otimes
\Omega[1] @>(\text{ - } \wedge \text{ - }) \otimes \Omega[1]>>
\sss(\Omega[1]) \otimes \Omega[1] \end{CD}\text{ . }$$

It follows that as a morphism in $\dcat$ from $\sss(\Omega[1])$ to $\sss(\Omega[1]) \otimes \Omega[1]$, $\Phi_L$ is given by the composite
$$\begin{CD}
\sss(\Omega[1]) \otimes \strc \\
@V{\sss(\Omega[1]) \otimes \text{id}}VV\\
\sss(\Omega[1]) \otimes T[-1] \otimes \Omega[1] \\
 @V{\sss(\Omega[1]) \otimes T[-1] \otimes  \frac{\text{At}_T}{\text{exp}(\text{At}_T)-1}}VV \\
\sss(\Omega[1]) \otimes T[-1] \otimes \sss(\Omega[1]) \otimes
\Omega[1] @>((\text{ - } \bullet \text{ - }) \wedge \text{ - })
\otimes \Omega[1]>>\sss(\Omega[1]) \otimes \Omega[1] \end{CD} $$.

This proves the desired proposition. \\

\end{proof}

By Proposition 30, $\Phi_L$ can be thought of as an element of \\ $\text{Hom}_{\dcat}(\strc,D_1 \otimes \Omega[1]) \simeq
\text{Hom}_{\dcat}(\strc, T[-1] \otimes \sss(\Omega[1]) \otimes \Omega[1])$. Let $\Phi_L^n$ denote the composite

$$\begin{CD} \strc \\
@V{\Phi_L^{\otimes n}}VV\\
 (T[-1] \otimes \sss(\Omega[1]) \otimes \Omega[1])^{\otimes n}\\
  @V{\tau}VV \\
  T^{\otimes n}[-n] \otimes \sss(\Omega[1])^{\otimes
n} \otimes \Omega^{\otimes n}[n] @>p' \otimes \text{m} \otimes p
>> \wedge^nT[-n] \otimes \sss(\Omega[1]) \otimes \Omega^n[n]
\end{CD}
$$.

The map $\tau$ in the above diagram is a rearrangement of factors. The map \\$\text{m}:\sss(\Omega[1])^{\otimes n} \rar \sss(\Omega[1])$ is
$n$-fold multiplication. $p'$ is the standard projection from $T^{\otimes n}[-n]$ to $\wedge^nT[-n]$ and $p$ is the projection from
$\Omega^{\otimes n}[n]$ to $\Omega^n[n]$. Note that $\Phi_L^n$ is an element of $\text{Hom}_{\dcat}(\strc,D_n \otimes \Omega^n[n])$. Let
$\textbf{1}^n:\sss(\Omega[1]) \rar D_n \otimes \Omega^n[n]$ be as in Section 3.2. Then,

\begin{lem}
The following diagram commutes in $\dcat$: \\

$$\begin{CD}
\strc @>\textbf{1}>> \strc \\
@VV{\textbf{f}}V   @V{\Phi_L^n}VV \\
\sss(\Omega[1]) @> \textbf{1}^n>> D_n \otimes \Omega^n[n]\\
\end{CD} $$
\end{lem}

\begin{proof}

{\it Step 1: The inverse of $\textbf{1}^n$.}\\

 {\it For this proof}, let $\textbf{e}:\strc \rar \Omega^n[n] \otimes \wedge^nT[-n]$ denote the natural isomorphism dual to the
 evaluation map from $\wedge^nT[-n] \otimes \Omega^n[n]$ to $\strc$. Denote the evaluation map from $\Omega^n[n] \otimes \wedge^n T[-n]$
 to $\strc$ by $\textbf{b}$ {\it for this proof}.
 We claim that the inverse to $\textbf{1}^n$ is given by the following
 composite:

\begin{equation} \label{inv}
\begin{CD}
\strc \otimes D_n \otimes \Omega^n[n] \\
@V{\textbf{e} \otimes D_n \otimes \Omega^n[n]}VV\\
\Omega^n[n] \otimes \wedge^n T[-n] \otimes D_n \otimes \Omega^n[n] \\
@V{\tau}VV \\
 \Omega^n[n] \otimes D_n \otimes \Omega^n[n] \otimes \wedge^n T[-n] \simeq \Omega^n[n] \otimes \wedge^n T[-n] \otimes \sss(\Omega[1]) \otimes
 \Omega^n[n] \otimes \wedge^n T[-n] \\
  @V{\textbf{b} \otimes \sss(\Omega[1]) \otimes \textbf{b}}VV \\
   \sss(\Omega[1]) \\
    \end{CD} \end{equation}.
The map $\tau$ in \eqref{inv} swaps $\wedge^nT[-n]$ and $D_n \otimes \Omega^n[n]$. \\

 To see this, note that if $U$ is an open subscheme of $X$ such that $\Omega |_U \simeq V \otimes_{\fld} \strcc{U}$ for some $n$-dimensional
 $\fld$-vector space $V$, the maps involved in \eqref{inv} can be described explicitly. Choose a basis $\{x_1,...,x_n\}$ of $V$ and a basis
 $\{y_1,...,y_n\}$ of $V^*$ dual to $\{x_1,..,x_n\}$. Let $\textbf{e}_V: \fld \rar \wedge^n V[n] \otimes \wedge^n V^*[-n]$ be the map taking $1$
 to $x_1 \wedge ... \wedge x_n \otimes y_n \wedge ... \wedge y_1$. Let $\textbf{b}_V: \wedge^n V[n] \otimes \wedge^n V^*[-n]$ be the map taking
 $x_1 \wedge .... \wedge x_n \otimes y_n \wedge ... \wedge y_1$ to $1$. Let $\textbf{1}^n_V$ be as in Section 3.1. The map denoted by $\textbf{1}^n_V$ was
 denoted by $\textbf{1}^m$ in Section 3.1. Then, if $H \in \sss(V[1])$ is
 a homogenous element,

 $$ \textbf{1}^n(H) = {(-1)}^{n|H|} y_n \wedge ..... \wedge y_1 \otimes H \otimes x_1 \wedge ..... \wedge x_n \text{ . } $$
 $$ (\textbf{e}_V \otimes \wedge^n V^*[-n] \otimes \sss(V[1]) \otimes \wedge^nV[n]) \circ \textbf{1}^n(H) $$ $$ = {(-1)}^{n|H|} x_1 \wedge .... \wedge x_n \otimes y_n \wedge ... \wedge y_1 \otimes  y_n \wedge ... \wedge y_1
 \otimes H \otimes x_1 \wedge .... \wedge x_n  \text{ . }$$
$$ \tau({(-1)}^{n|H|}x_1 \wedge .... \wedge x_n \otimes y_n \wedge ... \wedge y_1 \otimes  y_n \wedge ... \wedge y_1
 \otimes H \otimes x_1 \wedge .... \wedge x_n) $$ $$ = {(-1)}^{n|H|} {(-1)}^{n|H|} x_1 \wedge .... \wedge x_n \otimes  y_n \wedge ... \wedge y_1 \otimes  H \otimes
  x_1 \wedge .... \wedge x_n \otimes y_n \wedge .... \wedge y_1 \text{ . }$$
 $$(\textbf{b}_V \otimes \sss(V[1]) \otimes \textbf{b}_V) ({(-1)}^{n|H|} {(-1)}^{n|H|} x_1 \wedge .... \wedge x_n \otimes  y_n \wedge ... \wedge y_1 \otimes  H \otimes
  x_1 \wedge .... \wedge x_n \otimes y_n \wedge .... \wedge y_1 ) $$ $$ = H \text{ . } $$

  The fact that the composite given in \eqref{inv} is the inverse to $\textbf{1}^n$ follows from the facts that $\textbf{b} |_U = \textbf{b}_V
  \otimes_{\fld} \strcc{U}$, $\textbf{e} |_U = \textbf{e}_V \otimes_{\fld} \strcc{U}$ and $\textbf{1}^n |_U = \textbf{1}^n_V \otimes_{\fld}
  \strcc{U}$. \\

{\it Step 2:} \\

Let us look at the composition of the composite map in \eqref{inv} with $\Phi_L^n$. $\Phi_L^n$ can also be identified with the composite
$$(\textbf{b} \otimes \sss(\Omega[1]) \otimes \Omega^n[n]) \circ (\Omega^n[n] \otimes \Phi_L^n):\Omega^n[n] \rar \sss(\Omega[1]) \otimes \Omega^n[n] \text{ . }$$

Thinking of $\Phi_L^n$ as the above composite, it follows from
$\eqref{inv}$ that the inverse of $\textbf{1}^n$ composed with
$\Phi_L^n: \strc \rar D_n \otimes \Omega^n[n]$ is given by the
composite
\begin{equation} \label{finans} \begin{CD} \strc @>\textbf{e} >>
\Omega^n[n] \otimes \wedge^n T[-n] @> \Phi_L^n \otimes \wedge^n
T[-n]>> \sss(\Omega[1]) \otimes \Omega^n[n] \otimes \wedge^n T[-n]
@>\sss(\Omega[1]) \otimes \textbf{b}>> \sss(\Omega[1]) \end{CD}
\end{equation}.

Let $\textbf{b}_1: \Omega[1] \otimes T[-1] \rar \strc$ denote the evaluation map. Let $\text{m}:\sss(\Omega[1])^{\otimes n} \rar
\sss(\Omega[1])$ denote the $n$-fold
product. Let $p:\Omega[1]^{\otimes n} \rar \Omega^n[n]$ denote the natural projection. Note that
 $\Phi_L^n:\Omega^n[n] \rar \sss(\Omega[1]) \otimes \Omega^n[n]$ is also given by the following composite : \\

\begin{equation} \label{step2} \begin{CD}
\Omega^n[n] \\
@VVV \\
{(\Omega[1] \otimes \strc)}^{\otimes n}\\
 @VV{{(\Omega[1] \otimes \Phi_L)}^{\otimes n}}V \\
 {(\Omega[1] \otimes T[-1] \otimes \sss(\Omega[1]) \otimes \Omega[1])}^{\otimes n} \\
 @VV{{(\textbf{b}_1  \otimes \sss(\Omega[1]) \otimes \Omega[1])}^{\otimes n}}V \\
 {(\sss(\Omega[1]) \otimes \Omega[1])}^{\otimes n} \\
 @VV{\tau}V \\
 \sss(\Omega[1])^{\otimes n} \otimes \Omega[1]^{\otimes n} \\
 @VV{\text{m} \otimes p}V \\
 \sss(\Omega[1]) \otimes \Omega^n[n] \\
 \end{CD} \end{equation}.

The topmost vertical arrow in the above diagram is induced by the natural inclusion from $\Omega^n$ to $\Omega^{\otimes n}$. The map $\tau$ in
the above diagram is a rearrangement of factors. \\

It follows from Proposition 31 and that the composite \\ $(\textbf{b}_1 \otimes \sss(\Omega[1]) \otimes \Omega[1] ) \circ (\Omega[1] \otimes
\Phi_L):\Omega[1] \rar
\sss(\Omega[1]) \otimes \Omega[1]$ is precisely $\frac{\text{At}_T}{\text{exp}(\text{At}_T) -1}$. \\

Therefore the composite in \eqref{step2} is equal to the composite \\

\begin{equation*}
\begin{CD}
\Omega^n[n] \\
@VVV \\
\Omega[1]^{\otimes n} \\
@VV{(\frac{\text{At}_T}{{\text{exp}(\text{At}_T)-1}})^{\otimes n}}V\\
 (\sss(\Omega[1]) \otimes \Omega[1])^{\otimes n} @>\tau>> \sss(\Omega[1])^{\otimes n} \otimes \Omega[1]^{\otimes n} @>\text{m} \otimes
 p>>\sss(\Omega[1]) \otimes \Omega^n[n] \end{CD} \end{equation*}.

It follows from equation \eqref{finans} that the inverse of $\textbf{1}^n$ composed with $\Phi_L^n$ is given by the composite \\

$$\begin{CD}
\strc \\
@VV{\textbf{e}}V \\
\Omega^n[n] \otimes \wedge^nT[-n] \\
@VVV \\
\Omega[1]^{\otimes n} \otimes \wedge^nT[-n] \\
@VV{(\frac{\text{At}_T}{\text{exp}(\text{At}_T)-1})^{\otimes n} \otimes \wedge^nT[-n]}V\\
 (\sss(\Omega[1]) \otimes \Omega[1])^{\otimes n} \otimes \wedge^nT[-n] \\
 @VV{\tau}V \\
  \sss(\Omega[1])^{\otimes n} \otimes \Omega[1]^{\otimes n} \otimes \wedge^nT[-n] \\
   @VV{\text{m} \otimes  p \otimes \wedge^nT[-n]}V \\
  \sss(\Omega[1]) \otimes \Omega^n[n] \otimes \wedge^nT[-n] @>\sss(\Omega[1]) \otimes \textbf{b}>> \sss(\Omega[1]) \end{CD}
  $$.

The above composite is $\text{det}(\frac{\text{At}_T}{\text{exp}(\text{At}_T)-1})$ by the definition of $\text{det}$. \\

\end{proof}

\subsection{Another long remark - Lemma 4 and our dictionary}

Lemma 4 is the root reason for the Todd genus, an expression having
a form similar to the Jacobian of the differential
$d(\text{exp}^{-1})$ , showing up in the Riemann-Roch theorem.
Markarian [6] remarks that a Lemma in [6] similar to Lemma 3 in this
paper is like pulling back the canonical volume form on a Lie group
via the exponential map. He makes a remark in [2] that a formula
analogous to that describing the pull back of the canonical
(left-invariant) volume form on a Lie group via the exponential map
is responsible for the Todd genus showing up in the Riemann-Roch
theorem. Lemma 4 is precisely where something like this happens. We
will attempt to make the parallel between Lemma 4 and "pulling back
the canonical left-invariant volume form by the map
$\bar{\text{exp}}$" more transparent. In lemma 4
of this paper, it is $\bar{\text{exp}}$ rather than the exponential map itself that is involved. \\

{\it We also warn the reader that $\text{ad}$ and $\bar{\text{ad}}$
have the same meaning in this Section as in Section 2.5. Their
meaning in this section is therefore , different from their meaning
in Section 3, the rest of Section 4 and Section 5.}

\subsubsection{The classical situation} Keep the dictionary developed up-to Section 2.5 in mind. Let $G$ be a Lie group and $\mathfrak{g}$ its Lie algebra. Choose a basis
$\{X_i\}$ of $\mathfrak{g}$ and a basis $\{Y_i\}$ of
$\mathfrak{g}^*$ dual to $\{X_i\}$. Let $n$ be the dimension of
$\mathfrak{g}$. Let $\textbf{1}_{\mathfrak{g}}$ denote the element
$\sum_{i=1}^{i=n} X_i \otimes Y_i$ of $\mathfrak{g} \otimes
\mathfrak{g}^*$. Let $C(G)$ and $C(\mathfrak{g})$ be as in Section
2.5. Letting an element of $\mathfrak{g}$ act as a differential
operator on $C(G)$ (as in Section 2.5) yields a connection on
$C(G)$ for each element of $\mathfrak{g} \otimes \mathfrak{g}^*$. $\textbf{1}_{\mathfrak{g}}$ yields the canonical connection $d_G$. \\

The element $\textbf{1}^n := Y_n \wedge .... \wedge Y_1 \otimes X_1
\wedge ... \wedge X_n$ of $\wedge^n {\mathfrak{g}}^* \otimes
\mathfrak{g}^n$ yields a section of the trivial line bundle $\calg V
\times \wedge^n {\mathfrak{g}}^* \otimes \mathfrak{g}^n$ over the
neighborhood $\calg V$ of $0$ in $\mathfrak{g}$. This section yields
a map
$$\textbf{1}^n_{\mathfrak{g}} : C(\mathfrak{g}) \rar C(\mathfrak{g})
\otimes \wedge^n {\mathfrak{g}}^* \otimes \mathfrak{g}^n$$ $$ f
\leadsto f \otimes Y_n \wedge .... \wedge Y_1 \otimes X_1 \wedge ...
\wedge X_n  \text{ . }$$ The map $\textbf{1}^n$ in Lemma 4 is
analogous to
$\textbf{1}^n_{\mathfrak{g}}$. \\

On the other hand, if we think of an element of $\mathfrak{g}^*$ as
a function on $\calg V$, we have an $i+1$-fold multiplication
$\mu_i:C(\mathfrak{g}) \otimes {\mathfrak{g}^* }^{\otimes i} \rar
C(\mathfrak{g})$. Let $\text{ad}^{ \circ i}$ denote the composite
$$({\mathfrak{g}^*}^{\otimes i-1} \otimes \text{ad}) \circ ... \circ
\text{ad} : \mathfrak{g}^* \rar {\mathfrak{g}^*}^{\otimes i} \otimes
\mathfrak{g}^* \text{ . }$$ Note that $(\mu_i \otimes
\mathfrak{g}^*) \circ (C(\mathfrak{g}) \otimes \text{ad}^{\circ i})=
\bar{\text{ad}}^i$ where $\bar{\text{ad}}$ is as in Section 2.5.
Given any (convergent) power series $f(z) = \sum_i c_i z^i$ let
$f(\bar{\text{ad}})$ denote the map
$$\sum_i c_i \bar{\text{ad}}^i :C(\mathfrak{g}) \otimes  \mathfrak{g}^* \rar C(\mathfrak{g}) \otimes \mathfrak{g}^*  \text{ . }$$

Let $\Psi$ be as in Section 2.5. Then $$\Psi =
\frac{-\bar{\text{ad}}}{\text{e}^{\bar{\text{ad}}}-1} \circ
d_{\mathfrak{g}}$$ as a map from $C(\mathfrak{g})$ to
$C(\mathfrak{g}) \otimes {\mathfrak{g}}^*$. Therefore, as an element
of $C(\mathfrak{g}) \otimes \mathfrak{g} \otimes \mathfrak{g}^*$,
$$ \Psi = (\mathfrak{g} \otimes
\frac{-\bar{\text{ad}}}{\text{e}^{\bar{\text{ad}}}-1}) \circ
\textbf{1}_{\mathfrak{g}} \text{ . }$$

We can therefore think of the
element$\bigwedge^n_{C(\mathfrak{g})}(-\Psi)$ of $C(\mathfrak{g})
\otimes \wedge^n \mathfrak{g}^* \otimes \mathfrak{g}^n$ . Denote
this by ${(-\Psi)}^n$. Note that ${(-\Psi)}^n$ is a section of the
trivial line bundle $\calg V \times \wedge^n\mathfrak{g}^* \otimes
\mathfrak{g}^n$ over $\calg V$ . Moreover
$\textbf{1}^n_{\mathfrak{g}}$ is an isomorphism of $\mathbb
R$-vector spaces. We can therefore ask for the function
$$f_{\mathfrak{g}} :=
{(\textbf{1}^n_{\mathfrak{g}})}^{-1}({(-\Psi)}^n) \in
C(\mathfrak{g}) \text{ . }$$
 One can check from the formula for $-\Psi$ that
\begin{equation} \label{le4} f_{\mathfrak{g}} = \text{det}(\frac{\bar{\text{ad}}}{\text{e}^{\bar{\text{ad}}}-1}) \text{ . } \end{equation}

At this stage, we remark that the map $\Phi_L^n$ in Lemma 4 is
analogous to ${(-\Psi)}^n$.
$\text{det}(\frac{\bar{\text{ad}}}{\text{e}^{\bar{\text{ad}}}-1})$
is analogous to $\textbf{f}$ in Lemma 4. Lemma 4 itself is analogous
to the equation \eqref{le4}.

\subsubsection{Pulling back the canonical left invariant volume form on $G$ via $\bar{\text{exp}}$}

Finally, we observe that an element of $C(\mathfrak{g}) \otimes
\wedge^n \mathfrak{g}^* \otimes \mathfrak{g}^n$ and a volume form on
$\calg V$ together yield another volume form on $\calg V$ by letting
$\wedge^n \mathfrak{g}^*$ contract with $\mathfrak{g}^n$. In this
manner, the canonical volume form $Y_n \wedge...
\wedge Y_1$ on $\calg V$ and $\textbf{1}^n_{\mathfrak{g}}$ yield the canonical volume form $Y_n \wedge .... \wedge Y_1$ on $\calg V$. \\

Consider the left invariant volume form $\omega_G$ on $U_G$ arising
out of the element $Y_n \wedge... \wedge Y_1$ of $\wedge^n
\mathfrak{g}^*$. In the same manner the volume form $\omega$ yielded
by $\bar{\text{exp}}^*(\omega_G)$ and $\Psi^n$ equals $Y_n \wedge
... \wedge Y_1$ . But on the other hand, $\omega$ is also equal to
$\text{det}(d(\bar{\text{exp}}^*)) {(-1)}^n f_{\mathfrak{g}}$. It
follows that \eqref{le4} is equivalent to the formula
$$\text{det}(d(\bar{\text{exp}}^*)) =
\text{det}(\frac{-\bar{\text{ad}}}{\text{e}^{\bar{\text{ad}}}-1})
\text{ . }$$ This in turn is equivalent to the formula for the pull
back of a
left invariant volume form on $G$ via $\bar{\text{exp}}$.\\

\subsection{Proof of Lemma 3}

We are now equipped to prove Lemma 3. Let \\ $\text{ad}(\Phi_R): D_n
\rar D_n \otimes \Omega[1] $ denote the composite

$$ \begin{CD} D_n \otimes \strc
@>D_n \otimes \Phi_R>> D_n \otimes D_1 \otimes \Omega[1] @>\text{ad}
\otimes \Omega[1]>> D_n \otimes \Omega[1]  \end{CD} $$ \\

where $\text{ad}:D_n \otimes D_1 \rar D_n$ is as in Section 3.2. We begin with the following proposition. \\

\begin{prop}
$$\Phi_R^+(\textbf{f}) = 0  \text{ . }$$
\end{prop}

\begin{proof}

The upper square in the commutative diagram below commutes by Lemma
4. The lower square in the diagram below commutes by Lemma 2. \\

\begin{equation} \label{keystep}
\begin{CD}
\strc    @> \textbf{1} >> \strc \\
@VV{\textbf{f}}V    @V{\Phi_L^n}VV \\
\sss(\Omega[1]) @> \textbf{1}^n >> D_n \otimes \Omega^n[n] \\
@VV{-\Phi_R^+}V      @V{\text{ad}(\Phi_R) \otimes \Omega^n[n]}VV \\
\sss(\Omega[1]) \otimes \Omega[1] @> \textbf{1}^n \otimes \Omega[1]
>> D_n \otimes \Omega[1] \otimes \Omega^n[n] \\
\end{CD}
\end{equation}

 The proof of Proposition 30 with $\Phi_R$ instead of $\Phi_L$
 would show that \\ $\Phi_R \in \text{Hom}_{\dcat}(\strc,D_1 \otimes \Omega[1])$.
 Let $\bar{\text{ad}}$ be as in Section 3.2. Let
$\bar{\text{ad}}(\Phi_R)$ denote the composite
$$\begin{CD}
(D_1 \otimes \Omega[1])^{\otimes n} \otimes \strc \\
 @VV{(D_1 \otimes \Omega[1])^{\otimes n} \otimes \Phi_R}V\\
 (D_1 \otimes \Omega[1])^{\otimes n} \otimes D_1 \otimes \Omega[1]\\
@VVV\\
D_1^{\otimes n} \otimes D_1 \otimes \Omega[1]^{\otimes n} \otimes
\Omega[1] @>\bar{\text{ad}} \otimes \Omega[1]^{\otimes n} \otimes
\Omega[1]>> D_1^{\otimes n} \otimes \Omega[1]^{\otimes n} \otimes
\Omega[1] @>>> (D_1 \otimes \Omega[1])^{\otimes n} \otimes \Omega[1]
\end{CD} $$. The unlabeled arrows in the above diagram are
rearrangements of factors.\\

Note that by the Proposition 25, the following diagram commutes :
\\

\begin{equation} \label{nextstep}
\begin{CD}
\strc @> \textbf{1} >> \strc \\
@VV{\Phi_L^{\otimes n}}V    @V{\Phi_L^n}VV \\
(D_1 \otimes \Omega[1])^{\otimes n} @>>> D_n \otimes \Omega^n[n] \\
@VV{\bar{\text{ad}}(\Phi_R)}V          @V{\text{ad}(\Phi_R) \otimes \Omega^n[n]} VV \\
(D_1 \otimes \Omega[1])^{\otimes n} \otimes \Omega[1] @>>> D_n \otimes \Omega[1] \otimes  \Omega^n[n] \\
\end{CD}
\end{equation}

By Theorem 2' and by Proposition 9, $\Phi_R$ and $\Phi_L$ are commuting operators on $\sss(\Omega[1])$. It follows that
$\bar{\text{ad}}(\Phi_R)(\Phi_L^{\otimes n})$ is $0$. Thus, $\text{ad}(\Phi_R) \otimes \Omega^n[n](\Phi_L^n)$ is $0$. The desired
proposition now follows from \eqref{keystep} and the fact that $\textbf{1}^n \otimes \Omega[1]$ is invertible in $\dcat$.\\

\end{proof}

\subsection*{Proof of Lemma 3}

\begin{proof}

Note that the pairing $\langle,\rangle : \sss(\Omega[1])^{\otimes 2}
\rar \Omega^n[n]$ induces a non-degenerate pairing
$$\langle,\rangle: \rhh_X(\strc,\sss(\Omega[1]))^{\otimes 2} \rar
\rhh_X(\strc, \Omega^n[n]) \simeq \fld  \text{ . }$$ \\

$\langle \Phi_R(\text{ - }),\text{ - } \rangle: \sss(\Omega[1]) \otimes \sss(\Omega[1]) \rar \Omega^n[n] \otimes \Omega[1]$ denotes the composite
$$\begin{CD}
\sss(\Omega[1]) \otimes \sss(\Omega[1])\\
 @V{\Phi_R \otimes \sss(\Omega[1])}VV \\
  \sss(\Omega[1]) \otimes \Omega[1] \otimes \sss(\Omega[1]) @>\sss(\Omega[1]) \otimes \tau>> \sss(\Omega[1]) \otimes \sss(\Omega[1]) \otimes \Omega[1] @>\langle,\rangle
  \otimes \Omega[1] >> \Omega^n[n] \otimes \Omega[1] \end{CD} $$ map in $\dcat$. The map $\tau$ in the above composition of maps in $\dcat$ swaps $\Omega[1]$ and
  $\sss(\Omega[1])$. \\

Similarly, $\langle \text{ - }, \Phi_R^+(\text{ - }) \rangle:\sss(\Omega[1]) \otimes \sss(\Omega[1]) \rar \Omega^n[n] \otimes \Omega[1]$ denotes the composite
$$\begin{CD} \sss(\Omega[1]) \otimes \sss(\Omega[1]) @>\sss(\Omega[1]) \otimes \Phi_R^+>> \sss(\Omega[1]) \otimes \sss(\Omega[1]) \otimes
\Omega[1] @>\langle,\rangle \otimes \Omega[1]>> \Omega^n[n] \otimes
\Omega[1] \end{CD} \text{ . }$$

Let $a,b \in \rhh_X(\strc,\sss(\Omega[1]))$. Then,

\begin{equation} \label{s0} \langle ab.\textbf{f}^{-1},\Phi_R^+(\textbf{f}) \rangle = 0
\end{equation}
by Proposition 32. But, \begin{equation} \label{s1} \langle
ab.\textbf{f}^{-1},\Phi_R^+(\textbf{f})\rangle = \langle
\Phi_R(ab.\textbf{f}^{-1}),\textbf{f} \rangle \end{equation} by the
commutative diagram \eqref{adjdef}.\\

By Proposition 7, Theorem 2' and by the fact that $\hkr: \hc^{\bullet}(X) \rar \sss(\Omega[1])$ is a homomorphism of algebra objects in $\dcat$,
\begin{equation} \label{s2} \Phi_R(ab.\textbf{f}^{-1}) =
\Phi_R(a)b\textbf{f}^{-1} + a\Phi_R(b\textbf{f}^{-1}) \text{ . }
\end{equation}

By \eqref{s0} ,\eqref{s1} and \eqref{s2}

\begin{equation} \label{s3}  0 = \langle \Phi_R(ab.\textbf{f}^{-1}),\textbf{f} \rangle
= \langle \Phi_R(a)b\textbf{f}^{-1} ,\textbf{f} \rangle + \langle
a\Phi_R(b\textbf{f}^{-1}),\textbf{f} \rangle \text{ . }
\end{equation}

But for any elements $u,v,w \in \rhh_X(\strc,\sss(\Omega[1]))$,
$$\langle uv,w \rangle = \langle u,vw \rangle $$ by the definition
of $\langle,\rangle$. It then follows from \eqref{s3} that

\begin{equation} 0 = \langle \Phi_R(a), b.\textbf{f}^{-1}\textbf{f}
\rangle + \langle a, \Phi_R(b.\textbf{f}^{-1})\textbf{f} \rangle =
\langle \Phi_R(a), b \rangle + \langle a,
\Phi_R(b.\textbf{f}^{-1})\textbf{f} \rangle \text{ . }
\end{equation}

It follows from the commutative diagram \eqref{adjdef} that

$$\Phi_R^+(b) = - \Phi_R(b\textbf{f}^{-1})\textbf{f} $$ for any $b
\in \rhh_X(\strc, \sss(\Omega[1]))$. This proves Lemma 3.
\end{proof}

\textbf{Remark :} The computation proving Lemma 3 that has been
written here can very easily be rewritten in a "canonical" manner
without choosing elements of $\rhh_X(\strc,\sss(\Omega[1]))$. Though
that would be the ideal thing to do from the point of view of rigor,
we feel that the computation we have depicted conveys the key idea
behind the computation more concretely. Computations of a similar
nature that show up in
Section 5, have however, been written down in a "canonical" manner. \\

\section{Proof of Theorem 1}

\subsection{Unwinding some definitions}
In this subsection, we shall confine ourselves to unwinding the definition of the duality map $D_{\Delta}$. This will help us focus
more on what exactly we need to compute to prove Theorem 1. Let
$$ \kappa: \strcc{\Delta} \rar \Delta_* \Delta^* \strcc{\Delta} $$
denote the unit of the adjunction $\Delta^* \dashv \Delta_*$ applied to $\strcc{\Delta} \in \text{D}^{b}(X \times X)$. Also let
$$\beta: \strc \rar \Delta^* \Delta_! \strc
$$ denote the unit of the adjunction $\Delta_! \dashv \Delta^*$ applied to $\strc \in \dcat$. \\

Let $\edcat$ denote the category whose objects are those of $\dcat$ such that $$\text{Hom}_{\edcat}(\calg F,\calg G) = \rhh_{\dcat}(\calg
F,\calg G)$$ for any pair of objects $\calg F$ and $\calg G$ of $\dcat$. Note that any diagram that commutes in $\dcat$ also does so in
$\edcat$. We perform a particular calculation in $\edcat$ instead of $\dcat$ only when absolutely necessary. This enables us to take care of the
shifts in grading that occur when an element of $\rhh_{X}(\calg F,\calg G)$ shows up instead of an element of $\text{Hom}_{\dcat}(\calg F,\calg
G)$. In [8], $\edcat$ is called the "extended derived category"
of $X$. \\

Note that $\Delta_! \strc \simeq \Delta_* S_X^{-1}$. It follows that $\Delta^* \Delta_! \strc \simeq \Delta^* \strcc{\Delta} \otimes S_X^{-1} $.
Let \\ $\simeq:S_X \otimes S_X^{-1} \rar \strc$ be as in Section 3.2. We now state the following proposition.

\begin{prop}
Let $\phi \in \rhh_X(\Delta^* \strcc{\Delta}, S_X) $. Then, as a
morphism in $\edcat$, $D_{\Delta}^{-1}(\phi)$ is given by the
composite
$$\begin{CD} \strc \\
 @V{\beta}VV \\
  \Delta^* \strcc{\Delta} \otimes
S_X^{-1} @>\Delta^*(\kappa) \otimes S_X^{-1} >> \Delta^*
\strcc{\Delta} \otimes \Delta^* \strcc{\Delta} \otimes S_X^{-1}
@>\Delta^* \strcc{\Delta} \otimes \phi \otimes S_X^{-1} >> \Delta^*
\strcc{\Delta} \otimes S_X \otimes S_X^{-1} @>\Delta^* \strcc{\Delta} \otimes \simeq >> \Delta^* \strcc{\Delta} \\
\end{CD}
$$.
\end{prop}

\begin{proof}
By the definition of $D_{\Delta}$, $D_{\Delta}^{-1} = {\calg I
}^{-1} \circ {\calg T}^{-1} \circ {\calg J}^{-1} $ where $\calg I,
\calg T$ and $\calg J$ are as in \eqref{2},\eqref{3} and \eqref{4}
respectively. \\

Now, \begin{equation} \label{a} {\calg J}^{-1}(\phi)= \Delta_* \phi
\circ \kappa \text{ . } \end{equation}  Further,
\begin{equation}\label{b} {\calg T}^{-1}(\alpha) = \alpha \otimes
p_2^* S_X^{-1}
\end{equation} and
\begin{equation}\label{c} {\calg I}^{-1}(\gamma) = \Delta^* \gamma \circ
\beta \text{ . } \end{equation} Now, $\Delta^* (\alpha \otimes p_2^*
S_X^{-1}) = \Delta^* \alpha \otimes S_X^{-1} $ since $p_2 \circ
\Delta = \text{id}$. The desired proposition now follows from
\eqref{a} \eqref{b} and \eqref{c} and the fact that $\Delta^*
\Delta_* \psi = \Delta^* \strcc{\Delta} \otimes \psi$ for any
morphism $\psi$ in $\dcat$.

\end{proof}

The following propositions help us understand $\Delta^*(\kappa)$ and $\beta$ explicitly. \\

\begin{prop}
The following diagram commutes in $\dcat$: \\

$$\begin{CD}
\Delta^* \strcc{\Delta} @> \Delta^*(\kappa) >> \Delta^* \strcc{\Delta} \otimes \Delta^* \strcc{\Delta} \\
@V{\Delta^* \strcc{\Delta}}VV    @VV{\Delta^* \strcc{\Delta} \otimes \Delta^* \strcc{\Delta}}V \\
\Delta^* \strcc{\Delta} @> \textbf{C}>> \Delta^* \strcc{\Delta} \otimes \Delta^* \strcc{\Delta} \\
\end{CD} $$
\end{prop}

\begin{proof}

{\it Part 1:} \\
 Note that if $\calg F \in \text{D}^{b}(X \times
X)$, then
$$\Delta_* \Delta^* \calg F \simeq \strcc{\Delta} \otimes \calg F \text{ . } $$

Denote the canonical quotient map $\strcc{X \times X} \rar
\strcc{\Delta}$ by $\textbf{h}$. \\

Also observe that if $\calg G \in \dcat$, then $$\Delta^* \Delta_*
\calg G \simeq \Delta^* \strcc{\Delta} \otimes \calg G \text{ . }$$

Recall that $\Delta^* \strcc{\Delta}$ is represented by the complex
$\hc^{\bullet}(X)$. Note that the projection from the graded
$\strc$-module $\hc^{\bullet}(X)$ to $\hc^{0}(X) = \strc$ is a map
of complexes of $\strc$-modules. This was denoted by $\eta$ in
Section 2. In this proof, we will
denote this projection by $\textbf{p}$. \\

We claim that tensoring with $\textbf{h}$  constitutes the unit of
the adjunction $\Delta^* \dashv \Delta_*$ and that tensoring with
$\textbf{p}$ constitutes the co-unit of the adjunction $\Delta^*
\dashv \Delta_*$. \\

To see this, note that $\Delta^*(\textbf{h})$ is just the map
$$ \epsilon: \strc \rar \hc^{\bullet}(X) $$ defined in Section 2. This was the unit of the Hopf-algebra object $\hc^{\bullet}(X)$ of $\bacc$. It follows that
$\textbf{p} \circ \Delta^*(\textbf{h}) = \strc$. Also,
$$ \Delta_*(\textbf{p}) \circ (\textbf{h} \otimes \strcc{\Delta}) = \strcc{\Delta} $$
since $\textbf{h} \otimes \strcc{\Delta}$ can be identified with the map $\Delta_* (\epsilon)$ .\\

It follows that $\kappa = \textbf{h} \otimes \strcc{\Delta}$. \\

{\it Part 2:} \\

We now show that $\Delta^*(\kappa)$ and $\textbf{C}$ yield the same
morphism in $\dcat$ from $\Delta^* \strcc{\Delta}$ to $\Delta^*
\strcc{\Delta} \otimes \Delta^*
\strcc{\Delta}$. \\

Note that $\strcc{\Delta}$ is represented by the complex
$\hb^{\bullet}(X)$ in $\text{D}^{b}(X \times X)$. It follows that
both $\hb^{\bullet}(X) \otimes \strcc{\Delta}$ and $\strcc{\Delta}
\otimes \hb^{\bullet}(X)$ represent the object $\strcc{\Delta}
\otimes \strcc{\Delta}$ of $\text{D}^{b}(X \times X)$. Let
$$\nu:
\hb^{\bullet}(X) \otimes \strcc{\Delta} \rar \hb^{\bullet}(X)
\otimes \hb^{\bullet}(X) \otimes \strcc{\Delta}$$ denote the map
such that on an open subscheme $U= \text{Spec }R \times \text{Spec
}R$ of $X \times X$ before completion, $$\nu(r_0 \otimes ....
\otimes r_{k+1} \otimes_{R \otimes R} r') = \sum_{p+q=k;p,q \geq 0}
r_0 \otimes ... \otimes r_p \otimes 1 \otimes_{R \otimes R} 1
\otimes r_{p+1} \otimes ..... \otimes r_{k+1} \otimes_{R \otimes R}
r'$$. $\nu$ is easily seen to be a map of complexes of $\strcc{X
\times X}$-modules. Similarly, let $$\bar{\nu}:\strcc{\Delta}
\otimes \hb^{\bullet}(X)  \rar \strcc{\Delta} \otimes
\hb^{\bullet}(X) \otimes \hb^{\bullet}(X) $$ denote the map such
that on an open subscheme $U= \text{Spec }R \times \text{Spec }R$ of
$X \times X$ before completion,
$$\bar{\nu}(r'  \otimes_{R \otimes R} r_0 \otimes .... \otimes
r_{k+1})
 = r'  \otimes_{R \otimes R} \sum_{p+q=k;p,q \geq 0} r_0 \otimes ... \otimes r_p \otimes 1
\otimes_{R \otimes R} 1 \otimes r_{p+1} \otimes ..... \otimes
r_{k+1} $$. $\bar{\nu}$ is easily seen to be a map of complexes of
$\strcc{X \times X}$-modules. \\

Let $\tau:\strcc{\Delta} \otimes \hb^{\bullet}(X) \rar
\hb^{\bullet}(X) \otimes \strcc{\Delta}$ denote the map swapping
factors. Let $\tau':\strcc{\Delta} \otimes \hb^{\bullet}(X) \otimes
\hb^{\bullet}(X)  \rar \hb^{\bullet}(X) \otimes \hb^{\bullet}(X)
\otimes \strcc{\Delta}$ denote the map swapping $\strcc{\Delta}$ and
$\hb^{\bullet}(X) \otimes \hb^{\bullet}(X)$. The following diagram
then commutes \\
$$\begin{CD}
\strcc{\Delta} \otimes \hb^{\bullet}(X) @>\bar{\nu}>> \strcc{\Delta}
\otimes \hb^{\bullet}(X) \otimes \hb^{\bullet}(X) \\
@VV{\tau}V @V{\tau'}VV \\
\hb^{\bullet}(X) \otimes \strcc{\Delta} @>\nu>> \hb^{\bullet}(X)
\otimes \hb^{\bullet}(X) \otimes \strcc{\Delta}\\
\end{CD} $$

Note that $\Delta^*(\bar{\nu}) = \Delta^* \strcc{\Delta} \otimes
\textbf{C}$. It follows from this and from the above commutative
diagram that \begin{equation} \label{delll} \Delta^*(\nu) =
\textbf{C} \otimes \Delta^* \strcc{\Delta} \text{ . } \end{equation}

{\it Part 3:} We use $\eqref{delll}$ to compare the morphisms $\nu$ and $\kappa \otimes \strcc{\Delta}$ in $\text{D}^{b}(X \times X)$. Recall
that $\text{Hom}_{\text{D}^{b}(X \times X)}(\strcc{\Delta} \otimes \strcc{\Delta}, \strcc{\Delta} \otimes \strcc{\Delta} \otimes
\strcc{\Delta})$ is isomorphic to \\ $\text{Hom}_{\dcat}(\Delta^* \strcc{\Delta} \otimes \Delta^* \strcc{\Delta},\Delta^* \strcc{\Delta} \otimes
\Delta^* \strcc{\Delta})$. This isomorphism takes an element $\alpha$ of $\text{Hom}_{\text{D}^{b}(X \times X)}(\strcc{\Delta} \otimes
\strcc{\Delta}, \strcc{\Delta} \otimes \strcc{\Delta} \otimes \strcc{\Delta})$ to $(\textbf{p} \otimes \Delta^* \strcc{\Delta} \otimes \Delta^*
\strcc{\Delta}) \circ \Delta^*(\alpha)$. Since $\textbf{p}$ is induced by the co-unit of the Hopf-algebra object $\hc^{\bullet}(X)$ of $\bacc$,
and $\textbf{C}$ is induced by the co-multiplication of $\hc^{\bullet}(X)$, $$(\textbf{p} \otimes \Delta^* \strcc{\Delta} ) \circ \textbf{C} =
\textbf{1}_{\Delta^* \strcc{\Delta}} \text{ . } $$ Also, since $\Delta^*(\kappa)= \epsilon \otimes \Delta^* \strcc{\Delta}$,
$$(\textbf{p} \otimes \Delta^* \strcc{\Delta}) \circ \Delta^*(\kappa)
= \textbf{1}_{\Delta^* \strcc{\Delta}} \text{ . } $$ It follows that
$$ (\textbf{p} \otimes \Delta^* \strcc{\Delta} \otimes \Delta^*
\strcc{\Delta}) \circ \Delta^*(\nu) = (\textbf{p} \otimes \Delta^*
\strcc{\Delta} \otimes \Delta^* \strcc{\Delta}) \circ
\Delta^*(\kappa \otimes \strcc{\Delta}) = \textbf{1}_{\Delta^*
\strcc{\Delta} \otimes \Delta^* \strcc{\Delta}} \text{ . } $$

This proves that $\nu$ and $\kappa \otimes \strcc{\Delta}$ represent
the same morphism in $\text{D}^{b}(X \times X)$. Therefore ,
$\textbf{C} \otimes \Delta^* \strcc{\Delta}$ and $\Delta^*(\kappa)
\otimes \Delta^* \strcc{\Delta}$ represent the same morphism in
$\dcat$. The desired proposition now follows from the fact that if
$\lambda:\calg G \rar \calg H$ is a morphism in $\dcat$, $\lambda$
is equal to the composite $$\begin{CD} \calg G @>\calg G \otimes
\epsilon >> \calg G \otimes \Delta^* \strcc{\Delta} @>\lambda
\otimes \Delta^* \strcc{\Delta}>> \calg H \otimes \Delta^*
\strcc{\Delta} @>\calg H \otimes \textbf{p}>> \calg H \end{CD}
\text{ . } $$

\end{proof}

Recall that $p:\Omega[1]^{\otimes i} \rar \Omega^i[i]$ denotes the standard projection. Let $\Phi_R^{ \circ i}$ denote the composite
$$ \Phi_R \otimes \Omega[1]^{\otimes i-1} \circ ... \circ \Phi_R : \sss(\Omega[1]) \rar \sss(\Omega[1]) \otimes \Omega[1]^{\otimes i} $$
Let $\Phi_R^i$ denote $(\sss(\Omega[1]) \otimes p) \circ \Phi_R^{\circ i}$. Then $\text{exp}(\Phi_R) := \sum_i \frac{1}{i!} \Phi_R^i $ is a
morphism in $\dcat$ from $\sss(\Omega[1])$ to $\sss(\Omega[1]) \otimes \sss(\Omega[1])$. \\

The following proposition follows immediately from Proposition 8 and
Theorem 2'.\\

\begin{prop}
The following diagram commutes in $\dcat$: \\

$$\begin{CD}
\Delta^* \strcc{\Delta} @> \textbf{C}>> \Delta^* \strcc{\Delta} \otimes \Delta^* \strcc{\Delta} \\
@V{\hkr}VV   @VV{\hkr \otimes \hkr}V \\
\sss(\Omega[1]) @>\text{exp}(\Phi_R)>>  \sss(\Omega[1]) \otimes \sss(\Omega[1]) \\
\end{CD} $$

\end{prop}

Let $\iota: S_X \rar \sss(\Omega[1])$ denote the inclusion of $S_X \simeq \Omega^n[n]$ into $\sss(\Omega[1])$ as a direct summand. Let
$\bar{\beta}$ denote the composite
$$\begin{CD} \strc @>>> S_X \otimes S_X^{-1} @> \iota \otimes S_X^{-1}>> \sss(\Omega[1]) \otimes S_X^{-1} \end{CD} \text{ . }$$

\begin{prop}
The following diagram commutes in $\dcat$ :\\

$$\begin{CD}
\strc @> \beta>> \Delta^* \strcc{\Delta} \otimes S_X^{-1} \\
@VV{\textbf{1}}V @V{\hkr \otimes S_X^{-1}}VV \\
\strc @>\bar{\beta}>> \sss(\Omega[1]) \otimes S_X^{-1} \\
\end{CD} $$
\end{prop}

\begin{proof}

{\it Part 1: }\\

Note that if $\calg F \in \text{D}^{b}(X \times X)$, then
$$\Delta_! \Delta^* \calg F \simeq \Delta_* S_X^{-1} \otimes \calg F $$

Also, if $\calg G \in \dcat$ , then $$ \Delta^* \Delta_! \calg G \simeq \Delta^* \strcc{\Delta} \otimes S_X^{-1} \otimes \calg G $$

Now, $\Delta_* S_X^{-1}$ is isomorphic in $\text{D}^{b}(X \times X)$
to $\strcc{\Delta} \otimes p_2^*S_X^{-1}$. We now refer to the
statement of the Serre duality theorem in Markarian [6]. By the
Serre duality theorem there is a canonical map in $\text{D}^{b}(X
\times X)$ from $\strcc{\Delta}$ to $p_2^*S_X$. We denote this map
by $\textbf{q}$. Tensoring $\textbf{q}$ with $p_2^*S_X^{-1}$ on the
right and making the obvious identifications gives us a morphism
from $\strcc{\Delta} \otimes
p_2^* S_X^{-1}$ to $\strcc{X \times X}$. We denote this morphism by $\textbf{p}$ in this proof. \\

Let $\beta: \strc \rar \Delta^* \strcc{\Delta} \otimes S_X^{-1}$ be the morphism in $\dcat$ such that the diagram in this proposition commutes.
$\beta$ is well defined in $\dcat$ since $\hkr$ is a quasi-isomorphism. \\

We claim that tensoring by $\beta$ and tensoring by $\textbf{p}$
constitute the unit and co-unit of the adjunction $\Delta_! \dashv
\Delta^*$ respectively. In order to verify this claim, it suffices
to verify that \begin{equation} \label{prop361} \Delta^*(\textbf{p})
\circ \beta = \strc \end{equation} as morphisms in $\dcat$
and that \begin{equation} \label{prop362} (\textbf{p} \otimes \Delta_! \strc) \circ \Delta_! \beta = \Delta_! \strc \end{equation} as morphisms in $\text{D}^{b}(X \times X)$.\\

{\it Part 2: Verifying $\eqref{prop361}$ on "good" open subschemes of $X$ }\\

We begin by verifying $\eqref{prop361}$ on an open subscheme $U=
\text{Spec }R$ of $X$ with local coordinates $y_1,..,y_n$. The
elements $\{y_i \otimes 1-1 \otimes y_i \text{  } i=1,..,m \}$ form
a regular sequence generating the ideal $I$ of $R \otimes R$
defining the diagonal on an open affine neighborhood $V = \text{Spec
} S$ of the diagonal in $U \times U$. Let $z_i = y_i \otimes 1-1
\otimes y_i $. This regular sequence gives rise to a Koszul complex
${\calg K}^{\bullet}(z_1,..,z_n)$. ${\calg
K}^{\bullet}(z_1,..,z_n)$ is a free $S$-module resolution of $\Delta_* R$. \\

The third part of the Serre duality theorem as stated in [6] says
that the map $\textbf{q}$ (restricted to $V$) is equal to the map of
complexes
$${\calg K}^{\bullet}(z_1,...,z_n) \rar p_2^* \Omega^n_{R/\fld}[n]$$
$$ z_1 \wedge .... \wedge z_n \leadsto dy_1 \wedge .... \wedge dy_n
$$
as morphisms in $\text{D}^{b}(V)$.\\

 Further, let $[z_i]$
denote the class of $z_i$ in $\text{H}^*(R \otimes_{S} {\calg
K}^{\bullet}(z_1,...,z_n))$. The $R$-linear map $dy_i \leadsto
[z_i]$ induces an isomorphism of graded algebras between
$\sss(\Omega_{R/\fld}[1])$ and $\text{Tor}_*^{S}(R,R)$ by
Proposition 3.4.7 of Loday [9]. Moreover, by the proof of the
Hochschild-Kostant-Rosenberg theorem in Section 3.4 of Loday[9],
this isomorphism coincides with the isomorphism induced on
co-homology by the anti-symmetrization map
$\varphi:\sss(\Omega_{R/\fld}[1]) \rar \hc^{\bullet}(R)$. Before
completion, $$\varphi(r_0 dr_1 \wedge .... \wedge dr_k) =
\sum_{\sigma \in S_k} \text{sgn}(\sigma) r_0 \otimes r_{\sigma(1)}
\otimes ... \otimes r_{\sigma(k)}$$ This is immediately seen to be a
right-inverse of $\hkr$. It follows from the facts recalled in this
paragraph and the description of $\textbf{q}$ in the previous
paragraph that $\Delta^*(\textbf{q}) = \pi_n \circ \hkr$ as
morphisms in
$\text{D}^{b}(U)$. \\

Therefore, in $\text{D}^{b}(U)$, $\Delta^*(\textbf{p})$ is given by
the composite $$\begin{CD} \Delta^* \strcc{\Delta} \otimes S_X^{-1}
|_U @>\pi_n \circ \hkr \otimes S_X^{-1} |_U>> S_X |_U \otimes
S_X^{-1} |_U @>>> \strcc{U}=R \end{CD} $$

The unlabeled arrow in the above diagram is just the identification
of $S_X |_U \otimes S_X^{-1} |_U$ with $\strcc{U}=R$. It follows
that $\Delta^*(\textbf{p}) \circ \beta = \strcc{U}$ as morphisms in
$\text{D}^{b}(U)$. \\

{\it Part 3: Verifying $\eqref{prop362}$ on "good" open subschemes
of $X \times X$}

Let $U$ and $V$ be as in Part 2 of this proof. Note that $$\Delta_*
\hkr \circ ( \Delta_! \beta \otimes S_{X \times X}) =
\textbf{1}_{\Delta_* S_X}$$  Further, $$\textbf{p} \otimes \Delta_!
\strc \otimes S_{X \times X} = \textbf{q} \otimes \strcc{\Delta} $$
By the discussion in Part 2 of this proof, as morphisms in
$\text{D}^{b}(V)$, $$\textbf{q} \otimes \strcc{\Delta} = \Delta_*
\Delta^* \textbf{q} = \Delta_* (\pi_n \circ \hkr) $$ Therefore, as
morphisms in $\text{D}^{b}(V)$, $$(\textbf{q} \otimes
\strcc{\Delta}) \circ (\Delta_! \beta \otimes S_{X \times X}) =
\Delta_* (\pi_n \circ \hkr) \circ (\Delta_! \beta \otimes S_{X
\times X}) = \textbf{1}_{\Delta_* S_X}$$. Tensoring the morphisms
involved in the above equation with $S_{X \times X}^{-1}$, we see
that $$(\textbf{p} \otimes \Delta_! \strc) \circ \Delta_! \beta =
\Delta_! \strc $$ as morphisms in $\text{D}^{b}(V)$. This is what we
set out to verify. \\

{\it Part 4:}\\

Now observe that $$\text{Hom}_{\dcat}(\strc,\Delta^* \strcc{\Delta}
\otimes S_X^{-1}) \simeq \oplus_i \text{Hom}_{\dcat}(\strc, \Omega^i
\otimes \wedge^n T [i-n])
$$
For $i < n$, $$ \text{Hom}_{\dcat}(\strc, \Omega^i \otimes \wedge^n
T [i-n]) = \text{Ext}^{i-n}(\strc, \Omega^i \otimes \wedge^n T
[i-n]) = 0$$. For $i=n$, $$\text{Hom}_{\dcat}(\strc, \Omega^i
\otimes \wedge^n T [i-n]) = \text{Hom}_{\dcat}(\strc, \Omega^n
\otimes \wedge^nT) \simeq \text{Hom}_{\dcat}(\strc, \strc) \simeq
\fld$$ It follows that $\text{Hom}_{\dcat}(\strc,\Delta^*
\strcc{\Delta} \otimes S_X^{-1})$ is a $1$-dimensional $\fld$-vector
space. Tensoring with $\beta$ and tensoring with $\textbf{p}$
therefore, do indeed form a valid choice of unit and counit of the
adjunction $\Delta_! \dashv \Delta^*$ upto some scalar factors. In
other words, \eqref{prop361} and \eqref{prop362} are satisfied in
$\dcat$ and in $\text{D}^{b}(X \times X)$ respectively upto some
scalar factors. That the scalar factors are indeed $1$ follows from
the local verifications in Part 2 and Part 3
of this proof. \\

{\it Part 5:} \\

Let $\calg J$ be as in \eqref{4} in Section 1. The last detail to be
checked is that this particular choice of unit and co-unit for the
adjunction $\Delta_! \dashv \Delta^*$ satisfies
\begin{equation} \label{trchk} \text{tr}_X(\calg J(\phi) \circ \beta) = \text{tr}_{X \times
X}(\phi) \end{equation} for any element $\phi$ of
$\text{Hom}_{\text{D}^{b}(X \times X)}(\Delta_*S_X^{-1},
\Delta_*S_X)$. By the arguments in Part 4 of this proof, this
"compatibility with traces" is satisfied upto a scalar factor
independent of $\phi$ . We need to verify that this
scalar factor is $1$.\\

Let $\textbf{m}:\strcc{X \times X} \rar \strcc{\Delta}$ denote the obvious morphism in this proof.Let $\calg E$ and $\calg F$ be objects in
$\dcat$. Let $\calg E^*$ denote the dual $\text{RD}(\calg E)$ of $\calg E$. We recall the second part of the Serre Duality theorem as stated in
[6]. It states that if $f \in \text{Hom}_{\dcat}(\calg E, \calg F)$, the composite
$$\begin{CD} \strcc{X \times X} @>(\Delta_* f) \circ \textbf{m}>>  p_1^*\calg F \otimes p_2^*
\calg E^{*}  \otimes  \strcc{\Delta} @> p_1^*\calg F \otimes p_2^*
\calg E^{*} \otimes \textbf{q} >> p_1^*\calg F \otimes p_2^* (\calg
E^{*} \otimes S_X) \end{CD} $$ is exactly the image of the element
$f_*$ of $\text{Hom}_{\fld}(\text{Hom}_{\dcat}(\strc,\calg
E),\text{Hom}_{\dcat}(\strc, \calg F))$ under the identification

$$\begin{CD}
\text{Hom}_{\fld}(\text{Hom}_{\dcat}(\strc,\calg
E),\text{Hom}_{\dcat}(\strc, \calg F))\\
 @VVV \\
\text{Hom}_{\dcat}(\strc,\calg F) \otimes
\text{Hom}_{\dcat}(\strc,\calg E)^* \\
@VVV\\
 \text{Hom}_{\dcat}(\strc,\calg F) \otimes
\text{Hom}_{\dcat}(\strc,\calg E^* \otimes S_X) @>p_1^* \otimes
p_2^* >> \text{Hom}_{\text{D}^{b}(X \times X)}(\strcc{X \times
X},p_1^* \calg F \otimes p_2^* (\calg E \otimes S_X))\end{CD} $$
\\

Let $\textbf{tr} \in \text{Hom}_{\dcat}(\strc,S_X)$ be the element satisfying $\text{tr}_X(\textbf{tr})=1$. Applying the above statement with
$\calg E = \calg F = \strc$ and $f = \strc$, we see that the composite $$\begin{CD} \strcc{X \times X} @>\textbf{m}>>
 \strcc{\Delta}
@> \textbf{q} >>p_2^*  S_X \end{CD} $$
 is precisely
$p_2^*(\textbf{tr})$. It follows that the composite $\textbf{t}$
given by
$$\begin{CD} \strcc{X \times X} @>\textbf{m}>>
 \strcc{X \times X} \otimes \strcc{\Delta}
@> p_1^*(\textbf{tr}) \otimes \textbf{q} >>p_1^* S_X \otimes p_2^*
S_X = S_{X \times X}
\end{CD}
$$ satisfies $$\text{tr}_{X \times X}(\textbf{t})= 1$$.

Let $\textbf{tr}_{\Delta}$ denote the following composite:

$$\begin{CD} \strcc{X \times X} \otimes \strcc{\Delta}
@>p_1^*(\textbf{tr}) \otimes \textbf{q}>> S_{X \times X} \otimes
\strcc{X \times X} @>S_{X \times X} \otimes \textbf{m}>> S_{X \times
X} \otimes \strcc{\Delta} \end{CD} $$.

By Lemma 2.2 of Caldararu's paper [4],
$$ \text{tr}_{X \times X}(\textbf{tr}_{\Delta}) = \text{tr}_{X
\times X}(\textbf{t}) =1 \text{ . }$$ Therefore,
$$\text{tr}_{X \times X}(\textbf{tr}_{\Delta} \otimes p_2^*S_X^{-1}) =1
\text{ . }$$

To verify \eqref{trchk}, it suffices to check that
$$\text{tr}_X( \calg J (\textbf{tr}_{\Delta} \otimes p_2^*S_X^{-1}) \circ \beta) = 1 $$
where $\calg J$ is as in Section 1. \\

For this, it is enough to check that
\begin{equation} \label{trcheck2} \text{tr}_X(\calg J (\textbf{tr}_{\Delta}) \circ
(\beta \otimes S_X)) = 1 \text{ .} \end{equation}

Note that $\Delta^*(\textbf{tr}_{\Delta})$ is given by the composite
$$\begin{CD} \strc \otimes \Delta^*\strcc{\Delta} @>\textbf{tr}
\otimes \Delta^*(\textbf{q}) >> S_X \otimes S_X @>S_X \otimes S_X
\otimes \epsilon >> S_X \otimes S_X \otimes \Delta^* \strcc{\Delta}
\end{CD}$$ where $\epsilon$ is the map induced in $\dcat$ by the
unit $\strc \rar \hc^{\bullet}(X)$ of the Hopf-algebra object
$\hc^{\bullet}(X)$ of $\bacc$. It follows that $\calg
J(\textbf{tr}_{\Delta})$ is given by the composite
$$\begin{CD} \strc \otimes \Delta^*\strcc{\Delta} @>\textbf{tr}
\otimes \Delta^*(\textbf{q}) >> S_X \otimes S_X \end{CD} $$

To verify \eqref{trcheck2} it suffices to check that
$\Delta^*(\textbf{q}) \circ (\beta \otimes S_X) = S_X$. This is an
immediate consequence of $\eqref{prop361}$. \\

\end{proof}

\subsection{The final (long) computation proving Theorem 1}

\text{  }\\

\textbf{Recalling some notation:} \\
 We remind the reader that as in Section 3.2, the product on $\sss(\Omega[1])$ is denoted by $(\text{ - }\wedge\text{ - }):\sss(\Omega[1]) \otimes \sss(\Omega[1]) \rar \sss(\Omega[1])$. Let $\text{ev}$
be as in section 3.2. Denote the composite
$$\begin{CD} \sss(\Omega[1]) \otimes \sss(T[-1]) @>\sss(\Omega[1])
\otimes \textbf{j}>> \sss(\Omega[1]) \otimes \ennnn(\sss(\Omega[1]))
@>\text{ev}>> \sss(\Omega[1])  \end{CD} $$ by $(\text{ - } \bullet
\text{ - })$ , as in Section 3.2 . Note that the composite
$$\begin{CD} \sss(\Omega[1]) \otimes \sss(\Omega[1]) @> \sss(\Omega[1])
\otimes \textbf{i}>> \sss(\Omega[1]) \otimes \ennnn(\sss(\Omega[1]))
@>\text{ev}>> \sss(\Omega[1]) \end{CD} $$ is $(\text{ - } \wedge
\text{ - })$. Let $(\text{ - } || \text { - })$ be as in Section
3.2. We will also denote the isomorphisms $S_X \otimes S_X^{-1} \rar
\strc$
and $S_X^{-1} \otimes S_X \rar \strc$ by $\simeq$. \\

Also recall from Section 1 that $ \langle \text{ - },\text{ - }
\rangle = \pi_n \circ (\text{ - } \wedge \text{ - }):
\sss(\Omega[1]) \otimes \sss(\Omega[1]) \rar \sss(\Omega[1])
 \text{ . } $ Let
$x \in \rhh_X(\strc,\sss(\Omega[1]))$.\\

\textbf{Outline of the final computation proving Theorem 1:} The
final computation proving Theorem 1 begins by summarizing the
results of the previous subsection to express
$\hkr(D_{\Delta}^{-1}(\widehat{\hkr}(x)))$ as a composite of
morphisms in $\edcat$. This is done in Proposition 37. After
Proposition 37, in \eqref{tr3},\eqref{tr4},\eqref{tr5} and
\eqref{tr6}, the composite yielding
$\hkr(D_{\Delta}^{-1}(\widehat{\hkr}(x)))$ is rewritten to express
it in a form that is possible to simplify. Lemma 5 is then used to
simplify this composite further, yielding the composite
\eqref{kecom4}. The simplification using Lemma 5 is a crucial step.
The composite \eqref{kecom4} is further simplified to the composite
\eqref{kecom5}. The fact that the composite \eqref{kecom5} of
morphisms in $\edcat$ yields
$\hkr(D_{\Delta}^{-1}(\widehat{\hkr}(x)))$ together with Proposition
28 and Proposition 27 of Section 3.2 yield Theorem 1. Lemma 5 is then proven in Section 5.3 .\\

With the above notation and outline in mind, view $\Phi_R$ as an
element of \\ $\text{Hom}_{\dcat}( \strc,\ennnn(\sss(\Omega[1]))
\otimes \Omega[1])$. The following proposition begins the final set
of steps towards Theorem 1.
\begin{prop}
As a morphism in $\edcat$,
$\hkr(D_{\Delta}^{-1}(\widehat{\hkr}(x)))$ is given by the composite

$$\begin{CD} \strc \\
@V\bar{\beta}VV \\
 \sss(\Omega[1]) \otimes S_X^{-1} \simeq
\sss(\Omega[1]) \otimes \strc \otimes S_X^{-1}\\
 @V \sss(\Omega[1])
\otimes \text{exp}(\Phi_R) \otimes S_X^{-1}VV\\
 \sss(\Omega[1])
\otimes \ennnn(\sss(\Omega[1])) \otimes \sss(\Omega[1]) \otimes
S_X^{-1} \simeq \sss(\Omega[1]) \otimes \ennnn(\sss(\Omega[1]))
\otimes \sss(\Omega[1]) \otimes \strc \otimes
S_X^{-1}\\
@V{ \sss(\Omega[1]) \otimes \ennnn(\sss(\Omega[1])) \otimes
\sss(\Omega[1]) \otimes x \otimes S_X^{-1}}VV \\
\sss(\Omega[1]) \otimes \ennnn(\sss(\Omega[1])) \otimes
\sss(\Omega[1]) \otimes \sss(\Omega[1]) \otimes S_X^{-1} \\
@V{\sss(\Omega[1]) \otimes \ennnn(\sss(\Omega[1]))  \otimes \langle \text{ - } , \text{ - } \rangle \otimes S_X^{-1}}VV \\
\sss(\Omega[1]) \otimes \ennnn(\sss(\Omega[1])) \otimes S_X \otimes
S_X^{-1} \\
 @V{\text{ev} \otimes {\simeq}}VV\\
 \sss(\Omega[1])
\end{CD} $$.

\end{prop}

\begin{proof}

This proposition just amounts to putting together Propositions
33,34,35 and 36 and the definition of $\widehat{\hkr}$.

\end{proof}

As in Section 3.2,let $\gamma: \sss(\Omega[1]) \simeq \sss(T[-1]) \otimes S_X$ be the isomorphism such that ,
\begin{equation} \label{tr1} \langle \text{ - } , \text{ - } \rangle = \pi_n \circ (\text{ - } \wedge \text{ - }) = ((\text{ - }|| \text{ - })
\otimes S_X) \circ (\sss(\Omega[1]) \otimes \gamma) \text{ . } \end{equation} Let $\zeta$ denote the inverse of $\gamma$. \\

Also note that by the definitions of $F_r$ and $F_l$, the following
diagrams commute : \\

\begin{equation} \label{tr2}
\begin{CD}
\sss(\Omega[1]) \otimes \ennnn(\sss(\Omega[1])) @>\sss(\Omega[1]) \otimes \tau(F_l)>> \sss(\Omega[1]) \otimes \sss(\Omega[1]) \otimes
\sss(T[-1]) \\
@V{\text{ev}}VV           @VV{((\text{ - } \wedge \text{ - }) \bullet \text{ - })}V\\
\sss(\Omega[1]) @>\sss(\Omega[1])>> \sss(\Omega[1]) \\
\end{CD} \\ \\ \end{equation}

$$\begin{CD} \sss(\Omega[1]) \otimes \ennnn(\sss(\Omega[1]))
@>\sss(\Omega[1]) \otimes F_r>> \sss(\Omega[1]) \otimes \sss(T[-1])
\otimes
\sss(\Omega[1]) \\
@V{\text{ev}}VV           @VV{((\text{ - } \bullet \text{ - }) \wedge \text{ - })}V\\
\sss(\Omega[1]) @>\sss(\Omega[1])>> \sss(\Omega[1]) \\
\end{CD} $$

In the first diagram in \eqref{tr2}, $\tau:\sss(T[-1]) \otimes \sss(\Omega[1]) \rar \sss(\Omega[1]) \otimes \sss(T[-1])$ interchanges factors.By Proposition 37,
 $\hkr(D_{\Delta}^{-1}(\widehat{\hkr}(x)))$ is
given by the composite \\

\begin{equation} \label{tr3}
\begin{CD}
\strc \\
@V{\bar{\beta}}VV \\
\sss(\Omega[1]) \otimes S_X^{-1} \simeq \sss(\Omega[1]) \otimes \strc \otimes \strc \otimes S_X^{-1} \\
@V{\sss(\Omega[1]) \otimes  \text{exp}(\Phi_R) \otimes x \otimes
S_X^{-1}}VV \\
\sss(\Omega[1]) \otimes \ennnn(\sss(\Omega[1])) \otimes
\sss(\Omega[1]) \otimes \sss(\Omega[1]) \otimes S_X^{-1} \\
@V{ \text{ev} \otimes \langle \text{ - } , \text{ - } \rangle
\otimes
S_X^{-1}}VV \\
\sss(\Omega[1]) \otimes S_X \otimes S_X^{-1} \\
@V{\sss(\Omega[1]) \otimes \simeq}VV\\
\sss(\Omega[1]) \\
\end{CD}
\end{equation}

of morphisms in $\edcat$. Denote the composite
$$\begin{CD} \strc @>\text{exp}(\Phi_R)>> \ennnn(\sss(\Omega[1]))
\otimes \sss(\Omega[1]) @>\tau(F_l) \otimes \sss(\Omega[1])>>
\sss(\Omega[1]) \otimes \sss(T[-1]) \otimes \sss(\Omega[1]) \end{CD}
$$
 by $\textbf{P}$. \\

By \eqref{tr1} and \eqref{tr2} the composite \eqref{tr3} yielding
$\hkr(D_{\Delta}^{-1}(\widehat{\hkr}(x)))$ is the same
as the composite \\

\begin{equation} \label{tr4}
\begin{CD}
\strc \\
@V{\bar{\beta}}VV \\
\sss(\Omega[1]) \otimes S_X^{-1} \simeq \sss(\Omega[1]) \otimes \strc \otimes \strc \otimes S_X^{-1} \\
@V{\sss(\Omega[1]) \otimes \textbf{P} \otimes \gamma(x) \otimes
S_X^{-1}}VV \\
\sss(\Omega[1]) \otimes \sss(\Omega[1]) \otimes \sss(T[-1]) \otimes
\sss(\Omega[1]) \otimes \sss(T[-1]) \otimes S_X \otimes S_X^{-1} \\
@V{((\text{ - } \wedge \text{ - }) \bullet \text{ - }) \otimes
(\text{ - }
|| \text{ - } ) \otimes \simeq}VV \\
\sss(\Omega[1])\\
\end{CD}
\end{equation}

of morphisms in $\edcat$. Note that the composite
$$\begin{CD} \Omega^j[j] @> "\bar{\beta} \otimes \Omega^j[j]" >>
\sss(\Omega[1]) \otimes \Omega^j[j] \otimes S_X^{-1} @> (\text{ - }
\wedge \text{ - }) \otimes S_X^{-1} >> \sss(\Omega[1]) \otimes
 S_X^{-1} \end{CD} $$ is $0$. ($"\bar{\beta} \otimes \Omega^j[j]"$ is
 a rearrangment of factors composed with $\bar{\beta} \otimes \Omega^j[j]$
 ). It follows that the composite \eqref{tr4} yielding
$\hkr(D_{\Delta}^{-1}(\widehat{\hkr}(x)))$ is equal to the
 composite

 \begin{equation} \label{tr5}
\begin{CD}
\strc \\
@V{\bar{\beta}}VV \\
\sss(\Omega[1]) \otimes S_X^{-1} \simeq \sss(\Omega[1]) \otimes \strc \otimes \strc \otimes S_X^{-1} \\
@V{\sss(\Omega[1]) \otimes \textbf{Q} \otimes \gamma(x) \otimes
S_X^{-1}}VV \\
\sss(\Omega[1]) \otimes \sss(T[-1]) \otimes
\sss(\Omega[1]) \otimes \sss(T[-1]) \otimes S_X \otimes S_X^{-1} \\
@V{( \text{ - } \bullet \text{ - }) \otimes (\text{ - }
|| \text{ - } ) \otimes \simeq}VV \\
\sss(\Omega[1])\\
\end{CD}
\end{equation}

of morphisms in $\edcat$. In \eqref{tr5},  $\textbf{Q}$ in turn
denotes the composite

$$\begin{CD} \strc @>\text{exp}(\Phi_R)>> \ennnn(\sss(\Omega[1]))
\otimes \sss(\Omega[1]) @>\pi_0 \circ F_l \otimes \sss(\Omega[1])>>
\sss(T[-1]) \otimes \sss(\Omega[1]) \end{CD} \text{ . }$$\\

Let $(\text{exp}(\Phi_R) || \text{ - }) :\sss(T[-1]) \rar
\ennnn(\sss(\Omega[1]))$ denote the composite
$$\begin{CD}
 \strc \otimes \sss(T[-1]) \\
  @V{\text{exp}(\Phi_R) \otimes \sss(T[-1])}VV \\
  \ennnn(\sss(\Omega[1])) \otimes \sss(\Omega[1]) \otimes
\sss(T[-1]) @>\ennnn(\sss(\Omega[1])) \otimes ( \text{ - }||\text{ -
})>> \ennnn(\sss(\Omega[1])) \\ \end{CD} $$.

The composite in the diagram $\eqref{tr5}$ yielding
$\hkr(D_{\Delta}^{-1}(\widehat{\hkr}(x)))$ can be rewritten as \\

\begin{equation} \label{tr6}
\begin{CD}
\strc \\
@V{\bar{\beta}}VV \\
\sss(\Omega[1]) \otimes S_X^{-1} \simeq \sss(\Omega[1])  \otimes \strc \otimes S_X^{-1} \\
@V{\sss(\Omega[1]) \otimes \textbf{R}
\otimes S_X^{-1}}VV \\
\sss(\Omega[1]) \otimes \sss(T[-1]) \otimes S_X \otimes S_X^{-1} \\
@V{(\text{ - } \bullet \text{ - }) \otimes \simeq}VV \\
\sss(\Omega[1]) \\
\end{CD}
\end{equation}

where $$\textbf{R}=[\pi_0(F_l(\text{exp}(\Phi_R) || \text{ - }))
\otimes S_X] \circ \gamma(x) :\strc \rar \sss(T[-1]) \otimes S_X
$$ in $\edcat$.

We remind the reader that since $(\text{exp}(\Phi_R) || \text{ - })$
is a morphism in $\dcat$ from $\sss(T[-1])$ to
$\ennnn(\sss(\Omega[1]))$, $\pi_0(F_l(\text{exp}(\Phi_R) || \text{ -
}))$ is a morphism in $\dcat$ from $\sss(T[-1])$ to $\sss(T[-1])$.
It follows from this and the fact that $\gamma(x) \in
\text{Hom}_{\edcat}(\strc, \sss(T[-1]) \otimes S_X)$ that
$\textbf{R}$ is indeed an element of
$\text{Hom}_{\edcat}(\strc,\sss(T[-1]) \otimes S_X)$ as mentioned
above. The following Lemma simplifies the computation of
$\textbf{R}$. Let $(\text{ - }|\text{ - }):\sss(\Omega[1]) \otimes
\sss(T[-1]) \rar \sss(T[-1])$ be as in Section 3.2. Let
$(\text{td}_X^{-1}|\text{ - })$ denote the composite $(\text{ -
}|\text{ - }) \circ (\text{td}_X^{-1} \otimes \sss(T[-1]))$. This is
a morphism
in $\dcat$ from $\sss(T[-1])$ to itself.\\

\begin{lem}
  $$ \pi_0(F_l(\text{exp}(\Phi_R) || \text{ - })) =
  (\text{td}_X^{-1}|\text{ - }): \sss(T[-1]) \rar \sss(T[-1])$$
  \end{lem}

We postpone the proof of this lemma to Section 5.3. \\

It follows from \eqref{tr6} and Lemma 5 that
$\hkr(D_{\Delta}^{-1}(\widehat{\hkr}(x)))$ is given by the composite

\begin{equation} \label{kecom4}
\begin{CD}
\strc \\
@VV{\bar{\beta}}V \\
\sss(\Omega[1]) \otimes \strc \otimes S_X^{-1} \\
@VV{\sss(\Omega[1]) \otimes \gamma(x) \otimes S_X^{-1}}V \\
\sss(\Omega[1]) \otimes \sss(T[-1]) \otimes S_X \otimes S_X^{-1} \\
@VV{(\text{ - } \bullet (\text{td}_{X}^{-1} | \text{ - } )) \otimes \simeq}V \\
\sss(\Omega[1]) \\
\end{CD}
\end{equation} of morphisms in $\edcat$.\\

Now note that the following diagram commutes in $\edcat$: \\

$$\begin{CD}
\strc @>\textbf{1}_{\strc}>> \strc \\
@VVV @VV{\bar{\beta}}V \\
S_X \otimes \strc \otimes S_X^{-1} @>\iota \otimes \strc \otimes S_X^{-1}>>\sss(\Omega[1]) \otimes \strc \otimes S_X^{-1} \\
@V{S_X \otimes \gamma(x) \otimes S_X^{-1}}VV      @VV{\sss(\Omega[1]) \otimes \gamma(x) \otimes S_X^{-1}}V \\
S_X \otimes \sss(T[-1]) \otimes S_X \otimes S_X^{-1}  @>\iota \otimes \sss(T[-1]) \otimes S_X \otimes S_X^{-1} >> \sss(\Omega[1]) \otimes \sss(T[-1]) \otimes S_X \otimes S_X^{-1} \\
@V{(\text{ - } \bullet (\text{td}_{X}^{-1} | \text{ - } )) \otimes \simeq}VV @VV{(\text{ - } \bullet (\text{td}_{X}^{-1} | \text{ - } )) \otimes \simeq}V \\
\sss(\Omega[1]) @>\textbf{1}_{\sss(\Omega[1])}>> \sss(\Omega[1]) \\
\end{CD}$$

The topmost square in the above diagram commutes by the definition
of $\bar{\beta}$. That the remaining squares in the above diagram
commute is obvious. The reader should recall that $S_X =
\Omega^n[n]$ to make sense out of the map $(\text{ - } \bullet
(\text{td}_{X}^{-1} | \text{ - } )): S_X \otimes \sss(T[-1]) \rar
\sss(\Omega[1])$ in the above diagram.\\

It follows from the above diagram and \eqref{kecom4} that
$\hkr(D_{\Delta}^{-1}(\widehat{\hkr}(x)))$ is given by the composite

\begin{equation} \label{kecom5}
\begin{CD}
\strc \\
@VVV \\
S_X \otimes \strc \otimes S_X^{-1} \\
@VV{S_X \otimes \gamma(x) \otimes S_X^{-1}}V \\
S_X \otimes \sss(T[-1]) \otimes S_X \otimes S_X^{-1} \\
@VV{(\text{ - } \bullet (\text{td}_{X}^{-1} | \text{ - } )) \otimes \simeq}V \\
\sss(\Omega[1]) \\
\end{CD}
\end{equation} of morphisms in $\edcat$. \\

Let $\tau:S_X \otimes \sss(T[-1]) \otimes S_X \rar \sss(T[-1])
\otimes S_X \otimes S_X$ denote the map swapping $S_X$ and
$\sss(T[-1]) \otimes S_X$ . Let $\tau':S_X \otimes \strc \rar \strc \otimes S_X$ swap factors. We now have the following proposition. \\

\begin{prop}
The following diagram commutes in $\edcat$ : \\
$$\begin{CD}
\strc  @> \textbf{1}>> \strc \\
@VVV      @VVV \\
S_X \otimes \strc \otimes S_X^{-1} @> \tau' \otimes S_X^{-1}>> \strc \otimes S_X \otimes S_X^{-1} \\
@V{S_X \otimes \gamma(x) \otimes S_X^{-1}}VV     @VV{ \gamma(x) \otimes S_X \otimes S_X^{-1}}V\\
S_X \otimes \sss(T[-1]) \otimes S_X \otimes S_X^{-1}  @> \tau \otimes S_X^{-1} >>  \sss(T[-1]) \otimes S_X \otimes S_X \otimes S_X^{-1}  \\
@V{(\text{ - } \bullet (\text{td}_{X}^{-1} | \text{ - } )) \otimes \simeq}VV    @VV{ \zeta((\text{td}_{X}^{-1} | \text{ - } ) \otimes S_X) \otimes \simeq}V\\
\sss(\Omega[1])     @>J >>    \sss(\Omega[1])\\
\end{CD}
$$
\end{prop}

\begin{proof} The fact that the first two squares in the above
diagram commute is clear. The third square commutes by Proposition 28
of Section 3.2. \\
\end{proof}

Recall that $\gamma:\sss(\Omega[1]) \simeq \sss(T[-1]) \otimes S_X$ and $\zeta$ denotes the inverse of $\gamma$.
By Proposition 27 of Section 3.2,
$$\zeta( ((\text{td}_{X}^{-1} | \text{ - } ) \otimes
S_X)(\gamma(x)) = \text{td}_X^{-1} \wedge x \text{ . }$$

Therefore, by Proposition 38 , the fact that the composite
\eqref{kecom5} yields $\hkr(D_{\Delta}^{-1}(\widehat{\hkr}(x)))$,
and the fact that $J^2 = \textbf{1}_{\sss(\Omega[1])}$,

$$ \hkr(D_{\Delta}^{-1}(\widehat{\hkr}(x))) = J(\text{td}_X^{-1}
\wedge x ) \text{ . } $$

It follows that $$D_{\Delta}(\hkr^{-1}(y)) =
\widehat{\hkr}(\text{td}_X \wedge Jy) \text{ . } $$

This proves Theorem 1. \\

\subsection{Proving Lemma 5} Lemma 5 is the only thing left to be
proven. The proof of Lemma 5 uses Lemma 3. We first state and prove
the following proposition. All statements in this subsection hold in $\dcat$
and hence in $\edcat$. \\

\begin{prop}
$$\pi_0(F_r(\text{exp}(\Phi_R) || \text{ - })) = \textbf{1} :
\sss(T[-1]) \rar \sss(T[-1]) \text{ . } $$
\end{prop}

\begin{proof}

The notation used here is that used while proving Proposition 30.
While proving Proposition 30, we noted that $\bar{\omega}^i \circ
\bar{\textbf{C}}$ is given by the composite

$$\begin{CD}  \sss(\Omega[1]) @> \bar{\textbf{C}}>> \sss(\Omega[1])
\otimes \Omega[1] @> \sss(\Omega[1]) \otimes \text{At}_T^i>>
\sss(\Omega[1]) \otimes \Omega^i[i] \otimes \Omega[1] @>(\text{ - }
\wedge \text{ - }) \otimes \Omega[1]>> \sss(\Omega[1]) \otimes
\Omega[1] \\ \end{CD} $$.

It follows by the definition of $\pi_0 \circ F_r$ that
$$(\pi_0 \circ F_r \otimes \Omega[1])(\bar{\omega}^i \circ \bar{\textbf{C}}) = 0$$ if $i>0$.

It follows from the fact that $\frac{z}{1-\text{e}^{-z}} = 1 +
\sum_{i \geq 1} c_i z^i$ that
$$(\pi_0 \circ F_r \otimes \Omega[1])(\Phi_R) = (\pi_0 \circ F_r \otimes
\Omega[1])(\bar{\textbf{C}}) \text{ . } $$

Therefore, by Proposition 23 , $$(\pi_0 \circ F_r \otimes
\sss(\Omega[1]))(\text{exp}(\Phi_R)) = (\pi_0 \circ F_r \otimes
\sss(\Omega[1])) (\text{exp}(\bar{\textbf{C}})) \text{ . } $$

Note that $$\text{exp}(\bar{\textbf{C}}) = \textbf{C}_{\Omega}$$
where $\textbf{C}_{\Omega}$ is the co-multiplication on
$\sss(\Omega[1])$
treated as an element of \\ $\text{Hom}_{\dcat}(\strc,\ennnn(\sss(\Omega[1])) \otimes \sss(\Omega[1]))$. \\

For the rest of this proof let $(\textbf{C}_{\Omega} || \text{ - })$
denote $[\ennnn(\sss(\Omega[1])) \otimes (\text{ - } || \text{ - })]
\circ
\textbf{C}_{\Omega}$. \\

Therefore,
$$\pi_0(F_r(\text{exp}(\Phi_R) || \text{ - })) = \pi_0(F_r(\textbf{C}_{\Omega} || \text{ - })) :\sss(T[-1]) \rar \sss(T[-1]) \text{ . }$$

But $(\textbf{C}_{\Omega} || \text{ - }) = \textbf{j} : \sss(T[-1])
\rar \ennnn(\sss(T[-1])$ by Proposition 29. This proves the desired
proposition. \\

\end{proof}

The proof of Lemma 5 now follows.\\

\begin{proof}

 Let $(\text{exp}(\Phi_R) || \text{ - })^+:\sss(T[-1]) \rar
\ennnn(\sss(\Omega[1]))$ denote the map $\textbf{A} \circ
(\text{exp}(\Phi_R)||\text{ - }):\sss(T[-1]) \rar \ennnn(\sss(\Omega[1]))$ where $\textbf{A}$ is as in Section 3.2.\\

By Proposition 22, $$ \pi_0(F_l(\text{exp}(\Phi_R) || \text{ - })) =
I (\pi_0 (F_r(\text{exp}(\Phi_R) || \text{ - } )^+)) :\sss(T[-1])
\rar \sss(T[-1]) \text{ . }$$

By Part 2 of Proposition 21,

$$ (\text{exp}(\Phi_R) || \text{ - })^+ = (\text{exp}(\Phi_R^+) || \text{ - } ) :\sss(T[-1]) \rar \ennnn(\sss(\Omega[1])) \text{ . }$$

By Lemma 3,

$$ (\text{exp}(\Phi_R^+) || \text{ - } ) = (\text{exp}(- \textbf{i}(\textbf{f}) \circ \Phi_R \circ  \textbf{i}(\textbf{f})^{-1}) || \text{ - })$$ $$= \textbf{i}(\textbf{f}) \circ ( \text{exp}(-\Phi_R) || \text{ - }) \circ
\textbf{i}(\textbf{f}^{-1}) : \sss(T[-1]) \rar
\ennnn(\sss(\Omega[1]))\text{ . } $$\\

We remark that $\textbf{i}(\textbf{f}) \circ (\text{exp}(-\Phi_R)
||\text{ - }) \circ \textbf{i}(\textbf{f}^{-1})$ is precisely the
composite
$$\begin{CD}
\strc \otimes \sss(T[-1]) \otimes \strc \\
@VV{\textbf{i}(\textbf{f}) \otimes (\text{exp}(-\Phi_R )|| \text{ - }) \otimes \textbf{i}(\textbf{f}^{-1})}V\\
\ennnn(\sss(\Omega[1])) \otimes \ennnn(\sss(\Omega[1])) \otimes
\ennnn(\sss(\Omega[1]) @>\circ>> \ennnn(\sss(\Omega[1])) \end{CD}
$$.

Thus,

 \begin{equation} \label{keycomp1} \pi_0(F_l(\text{exp}(\Phi_R)
|| \text{ - }))   = I (\pi_0 (F_r(\textbf{i}(\textbf{f}) \circ (
\text{exp}(-\Phi_R) || \text{ - }) \circ
\textbf{i}(\textbf{f}^{-1})))) : \sss(T[-1]) \rar \sss(T[-1]) \text{
. } \end{equation}

Note that $\pi_0(F_r(\textbf{i}(\textbf{f}))) = 1$ since
$\textbf{f}= \text{det}(1+\sum_{i>0} c_i \text{At}_{T}^i)$ . It
follows from Proposition 23 that

\begin{equation} \label{keycomp2} \pi_0(F_l(\text{exp}(\Phi_R)
|| \text{ - })) = I (\pi_0 (F_r(( \text{exp}(-\Phi_R) || \text{ - })
\circ \textbf{i}(\textbf{f}^{-1})))) : \sss(T[-1]) \rar \sss(T[-1])
\text{ . } \end{equation}

Another point to note is that $$(\text{exp}(-\Phi_R) || \text{ - })
= (\text{exp}(\Phi_R) || \text{ - } ) \circ I : \sss(T[-1]) \rar
\ennnn(\sss(\Omega[1]))
$$. It follows from this observation , and \eqref{keycomp2} that

$$ \pi_0(F_l(\text{exp}(\Phi_R)
|| \text{ - })) = I (\pi_0 (F_r(( \text{exp}(\Phi_R) || \text{ - })
\circ \textbf{i}(\textbf{f}^{-1})))) \circ I $$

$$ = I( \pi_0 (F_r(\textbf{j}
(\pi_0(F_r(\text{exp}(\Phi_R) || \text{ - }))) \circ
\textbf{i}(\textbf{f}^{-1})))) \circ I \text{ by Proposition 23 } $$

$$ =  I(\pi_0(F_r( \textbf{j}(I(\text{ - })) \circ
\textbf{i}(\textbf{f}^{-1})))) \text{  by Proposition 39 } $$

$$ = I( \textbf{f}^{-1} | I(\text{ - }) ) :\sss(T[-1]) \rar \sss(T[-1]) \text{ by Proposition 26 } $$

where $(\text{ - }|\text{ - }):\sss(\Omega[1]) \otimes \sss(T[-1])
\rar
\sss(T[-1])$ is as in Section 3.2. \\

We remind the reader that $(\textbf{f}^{-1}|I(\text{ -
})):\sss(T[-1]) \rar \sss(T[-1])$ is the composite
$$\begin{CD}
\sss(T[-1]) \\
@VV{I}V \\
\strc \otimes \sss(T[-1]) \\
@VV{\textbf{f}^{-1} \otimes \sss(T[-1])}V \\
\sss(\Omega[1]) \otimes \sss(T[-1]) @>(\text{ - }|\text{ - })>>
\sss(T[-1]) \\ \end{CD} $$ of morphisms in $\dcat$.\\

Note that $$I(\textbf{f}^{-1} | I(\text{ - })) = (J(\textbf{f}^{-1})
|\text{ - }) :\sss(T[-1]) \rar \sss(T[-1]) \text{ . }
$$ Also observe that $J(\textbf{f}^{-1}) = \text{td}_{X}^{-1} $.

This proves Lemma 5.\\
\end{proof}

\end{document}